\documentclass{amsart}
\usepackage{enumitem}
\usepackage{amssymb}
\usepackage{amsthm}
\usepackage{amsmath}

\usepackage[scr=boondox]{mathalfa}

\usepackage[usenames]{color}

\usepackage{hyperref}

\theoremstyle{plain}
\newtheorem{proposition}{Proposition}[section]
\newtheorem{theorem}[proposition]{Theorem}
\newtheorem{lemma}[proposition]{Lemma}
\newtheorem{corollary}[proposition]{Corollary}
\theoremstyle{definition}

\newtheorem{definition}[proposition]{Definition}
\newtheorem{observation}[proposition]{Observation}
\theoremstyle{remark}
\newtheorem{remark}[proposition]{Remark}

\DeclareMathOperator{\SL}{\mathsf{SL}}
\DeclareMathOperator{\GL}{\mathsf{GL}}
\DeclareMathOperator{\SO}{\mathsf{SO}}

\DeclareMathOperator{\PSL}{\mathsf{PSL}}

\DeclareMathOperator{\Hom}{Hom}

\DeclareMathOperator{\End}{End}

\DeclareMathOperator{\Gr}{Gr} 
\DeclareMathOperator{\id}{id} 
 
\DeclareMathOperator{\CAT}{CAT}

\DeclareMathOperator{\Span}{Span}
\DeclareMathOperator{\Ad}{Ad}
\DeclareMathOperator{\nL}{\mathfrak{n}}
\DeclareMathOperator{\gL}{\mathfrak{g}}

\DeclareMathOperator{\Ac}{\mathcal{A}}
\DeclareMathOperator{\Bc}{\mathcal{B}}
\DeclareMathOperator{\Cc}{\mathcal{C}}
\DeclareMathOperator{\Fc}{\mathcal{F}}

\DeclareMathOperator{\Oc}{\mathcal{O}}

\DeclareMathOperator{\Cb}{\mathbb{C}}

\DeclareMathOperator{\Hb}{\mathbb{H}}

\DeclareMathOperator{\Kb}{\mathbb{K}}
\DeclareMathOperator{\Nb}{\mathbb{N}}
\DeclareMathOperator{\Pb}{\mathbb{P}}
\DeclareMathOperator{\Rb}{\mathbb{R}}
\DeclareMathOperator{\Sb}{\mathbb{S}}

\DeclareMathOperator{\Zb}{\mathbb{Z}}
\DeclareMathOperator{\Qb}{\mathbb{Q}}

\DeclareMathOperator{\kL}{\mathfrak{k}}
\DeclareMathOperator{\pL}{\mathfrak{p}}

\DeclareMathOperator{\aL}{\mathfrak{a}}

\DeclareMathOperator{\Usf}{\mathsf{U}}

\DeclareMathOperator{\u2}{\mathscr{u}}
\DeclareMathOperator{\a2}{\mathscr{a}}

\newcommand{\abs}[1]{\left|#1\right|}

\newcommand{\norm}[1]{\left\|#1\right\|}
\newcommand{\wt}[1]{\widetilde{#1}}
\newcommand{\wh}[1]{\widehat{#1}}
\newcommand{\ip}[1]{\left\langle #1\right\rangle}

\newcommand{\vertiii}[1]{{\left\vert\kern-0.25ex\left\vert\kern-0.25ex\left\vert #1 
    \right\vert\kern-0.25ex\right\vert\kern-0.25ex\right\vert}}

\begin{document}


\title[Cusped Hitchin representations]{Cusped Hitchin representations and Anosov representations of geometrically finite Fuchsian groups}
\author[Canary]{Richard Canary}
\address{University of Michigan}
\author[Zhang]{Tengren Zhang}
\address{National University of Singapore}
\author[Zimmer]{Andrew Zimmer}
\address{University of Wisconsin-Madison}
\thanks{Canary was partially supported by grant  DMS-1906441 from the National Science Foundation
and grant 674990 from the Simons Foundation. Zhang was partially supported by the NUS-MOE grant R-146-000-270-133. Zimmer was partially supported by grants DMS-2105580 and DMS-2104381 from the
National Science Foundation.}

\date{\today}

\begin{abstract} 
We develop a theory of Anosov representation of geometrically finite Fuchsian groups in $\mathsf{SL}(d,\mathbb R)$
and show that cusped Hitchin representations are Borel Anosov in this sense. We establish analogues of many properties
of traditional Anosov representations. In particular, we show that our Anosov representations are stable under type-preserving
deformations and that  their limit maps vary analytically.
We also observe that our Anosov representations
fit into the previous frameworks of relatively Anosov and relatively dominated representations developed by Kapovich-Leeb and Zhu.

\end{abstract}

\maketitle

\setcounter{tocdepth}{1}
\tableofcontents

\section{Introduction}

In this paper we develop a theory of Anosov representations for geometrically finite Fuchsian groups into $\mathsf{SL}(d,\mathbb K)$ where $\Kb$ is either the field of real numbers or the field of complex numbers (and more generally, into any semisimple Lie group $G$ with finite center).
This theory shares two crucial features with Anosov representations of convex cocompact Fuchsian groups. First, if $\Gamma\subset\PSL(2,\Rb)$
is a geometrically finite Fuchsian group and $\rho:\Gamma\to\mathsf{SL}(d,\mathbb K)$ is Anosov, then there is a $\rho$-equivariant
quasi-isometric embedding of an orbit of $\Gamma$ in $\mathbb H^2$ into the Riemannian symmetric space $X_d(\mathbb K)$ associated to
$\mathsf{SL}(d,\mathbb K)$. Moreover, there is a 
$\rho$-equivariant map of the limit set of $\Gamma$ into the appropriate  (partial) flag variety. Second, small deformations
of Anosov representations which preserve the conjugacy class of the images of parabolic elements in $\Gamma$ remain Anosov.
We further show that these limit maps vary analytically in $\rho$.
We observe that our representations are  relatively Anosov representations in the sense of Kapovich and Leeb \cite{kapovich-leeb}
and relatively dominated representations as defined  by Zhu \cite{zhu-reldom}. Our concrete setting allows for simpler proofs
and  more explicit results.

Our main motivation was to study the class of cusped Hitchin representations.  A representation $\rho:\Gamma\to\mathsf{SL}(d,\mathbb R)$
of a geometrically finite Fuchsian group is said to be \emph{cusped Hitchin} if there exists a continuous $\rho$-equivariant positive map from the limit set of $\Gamma$ into the full flag
variety $\mathcal F_d$. We show that cusped Hitchin representations are irreducible and Borel Anosov (i.e. they are $P_k$-Anosov
for all $1\le k\le d-1$). This generalizes results of Labourie \cite{labourie-invent}, when $\Gamma$ is cocompact, and
Labourie-McShane \cite{labourie-mcshane} (see also Burelle-Treib \cite{burelle-treib}), when $\Gamma$ is convex cocompact. (Recent results of Sambarino \cite{sambarino-positive} also imply 
that cusped Hitchin representations are irreducible).
Fock and Goncharov \cite{fock-goncharov}
introduced the theory of positive representations. They consider the case where $\Gamma$ is not cocompact but has cofinite area, and they only require that
the limit map be defined on fixed points of peripheral elements of $\Gamma$. We  show that type-preserving positive representations,
in their sense, are in fact cusped Hitchin representations. Other examples of  cusped Anosov representations are provided by
exterior powers of cusped Hitchin representations and direct products of cusped Hitchin representations with trivial representations.

In turn, the motivation for studying cusped Hitchin representations arises from an intriguing potential analogy with the augmented
Teichm\"uller space from classical Teichm\"uller theory. The augmented Teichmuller space of a closed orientable surface $S$ is obtained
by appending to the Teichm\"uller space of $S$  points corresponding to (marked) finite area hyperbolic structures on the complement in $S$ of
any multicurve on $S$. Masur \cite{masur-wp} showed that the augmented Teichm\"uller space is the metric completion
of Teichm\"uller space with the Weil-Petersson metric. Loftin and Zhang \cite{loftin-zhang} construct an analytic model for
an  augmented Hitchin component. Bray, Canary, Kao and Martone \cite{BCKM2}
combine our work with thermodynamical results from 
\cite{BCKM} to construct pressure metrics on deformation spaces of cusped Hitchin
representations, generalizing work of Bridgeman, Canary, Labourie and Sambarino \cite{BCLS}. 
The hope is that this will allow us to investigate whether the augmented Hitchin component is the metric completion of
the Hitchin component with the pressure metric. For further discussion of the conjectural geometric picture see the survey paper
\cite{canary-hitchin}.

\medskip

We now turn to a more detailed discussion of our work. Let $\Gamma\subset \mathsf{PSL}(2,\mathbb R)$
be a geometrically finite group with limit set $\Lambda(\Gamma)\subset\partial\mathbb H^2$. Suppose that $\rho:\Gamma\to\mathsf{SL}(d,\mathbb K)$
is a representation and that we are given a continuous $\rho$-equivariant map
$$\xi_\rho=(\xi_\rho^{k},\xi_\rho^{d-k}):\Lambda(\Gamma)\to\mathrm{Gr}_k(\mathbb K^d)\times\mathrm{Gr}_{d-k}(\mathbb K^d)$$
into the Grassmanians of $k$-planes and $(d-k)$-planes in $\mathbb K^d$. We require that $\xi_\rho$ is   transverse, i.e.
$$\xi_\rho^{k}(x)\oplus\xi_\rho^{d-k}(y)=\mathbb K^d\quad  \mathrm{if}\quad  x\ne y\in\Lambda(\Gamma).$$
We obtain an associated splitting
$$E_\rho=\Usf(\Gamma) \times \Kb^d = \Theta^k\oplus\Xi^{d-k}\quad \mathrm{where}\quad \Theta^k|_v=\xi^k_\rho(v^+),\quad \Xi^{d-k}|_v=\xi^{d-k}_\rho(v^-)$$
and $\Usf(\Gamma)\subset T^1\mathbb H^2$ is the set of tangent vectors which extend to geodesics both of
whose endpoints lie in $\Lambda(\Gamma)$. This descends to a splitting
$$\wh E_\rho=\Gamma\backslash E_\rho=\wh \Theta_\rho^k\oplus\wh\Xi^{d-k}_\rho$$
of the flat bundle associated to $\rho$ over the non-wandering portion 
$\wh\Usf(\Gamma)=\Gamma \backslash \Usf(\Gamma)$ of the geodesic flow.
The geodesic flow on $\wh\Usf(\Gamma)$ lifts naturally to a flow on $\wh E_\rho$ which preserves
the splitting and is parallel to the flat connection. We say that $\rho$ is \emph{$P_k$-Anosov} if the associated flow on
$\mathrm{Hom}(\wh\Xi_\rho^{d-k},\wh\Theta_\rho^{k})$ is uniformly contracting with respect to the (operator) norm 
arising from \textbf{some} family of continuous norms on the fibers of $\wh E_\rho$ (see Definition \ref{defn:anosov_specific} for details).
In this case, we call $\xi_\rho$ the {\em $P_k$-Anosov limit map} of $\rho$. When $\Gamma$ contains a parabolic element, we will sometimes refer to our $P_k$-Anosov representations as {\em cusped $P_k$-Anosov} representations to
distinguish them from traditional Anosov representations.

We obtain generalizations of many of the classical properties of Anosov representations, see Labourie \cite{labourie-invent}
or Guichard-Wienhard \cite{guichard-wienhard}, in our setting. 
We say that a $\rho$-equivariant, continuous map $\xi=(\xi^k,\xi^{d-k}):\Lambda(\Gamma)\to\mathrm{Gr}^k(\mathbb K^d)\times\mathrm{Gr}^{d-k}(\mathbb K^d)$
is {\em strongly dynamics preserving} if whenever $\{\gamma_n\}$ is a sequence in $\Gamma$,
$\gamma_n(z) \rightarrow x \in\Lambda(\Gamma)$, and $\gamma_n^{-1}(z) \rightarrow y\in\Lambda(\Gamma)$ for some (any) $z\in\Hb^2$, then
\begin{align*}
\rho(\gamma_n)(V) \rightarrow \xi^{k}(x)
\end{align*}
for any $V\in\mathrm{Gr}(\mathbb K^d)$ which is transverse to $\xi^{d-k}(y)$. Given $g \in \SL(d, \mathbb K)$ let 
$$
\lambda_1(g) \geq \dots \geq \lambda_d(g)
$$
denote the absolute values of the (generalized) eigenvalues of $g$. Then $g$ is {\em $P_k$-proximal} if $\lambda_k(g)>\lambda_{k+1}(g)$ and $g$ is \emph{weakly unipotent} if  the (multiplicative) Jordan-Chevalley decomposition of $g$ has elliptic semisimple part and non-trivial unipotent part. 

\begin{theorem}\label{thm: intro thm 1}
If $\Gamma\subset\mathsf{PSL}(2,\mathbb R)$ is a geometrically finite group and $\rho:\Gamma\to\mathsf{SL}(d,\mathbb K)$
is $P_k$-Anosov, then
\begin{enumerate}
\item For any $z_0 \in \Hb^2$, there exists $A,a>1$ so that if $\gamma\in\Gamma$, then
$$\frac{1}{A}\exp\left(\frac{1}{a}d_{\Hb^2}(z_0,\gamma (z_0))\right)\leq \frac{\sigma_k(\rho(\gamma))}{\sigma_{k+1}(\rho(\gamma))} \leq A\exp\left(a d_{\Hb^2}(z_0,\gamma (z_0))\right).$$
\item
There exists $B,b>1$ so that if $\gamma\in\Gamma$, then
$$\frac{1}{B}\exp\left( \frac{1}{b}\ell(\gamma)\right)\leq \frac{\lambda_k(\rho(\gamma))}{\lambda_{k+1}(\rho(\gamma))} \leq B\exp\left( b\ell(\gamma)\right)$$
where  $\ell(\gamma)$ is the translation length of $\gamma$ on $\mathbb H^2$.
\item
The $P_k$-Anosov  limit map $\xi_\rho$ is strongly dynamics-preserving and unique. In particular, if $\alpha\in\Gamma$ is parabolic, then $\rho(\alpha)$ is weakly unipotent,
while if $\gamma\in\Gamma$ is hyperbolic, then $\rho(\gamma)$ is $P_k$-proximal.
 \item
If $z_0\in\mathbb H^2$ and $x_0$ is a point in the symmetric space $X_d(\mathbb K)$ associated to $\mathsf{SL}(d,\mathbb K)$, then
the orbit map $\tau_\rho:\Gamma(z_0)\to X_d(\mathbb K)$ given by $\tau_\rho(\gamma(z_0))=\rho(\gamma)(x_0)$ is a quasi-isometric embedding.
\end{enumerate}
\end{theorem}

We also give a dynamical characterization of Anosov representations in the spirit of characterizations of traditional Anosov representations by
G\'ueritaud-Guichard-Kassel-Wienhard \cite{GGKW}, Kapovich-Leeb-Porti \cite{KLP} and Tsouvalas \cite{kostas}.

\begin{theorem}\label{thm:equivalent_to_Anosov}
Suppose $\Gamma \subset \PSL(2,\Rb)$ is a geometrically finite group and $\rho : \Gamma \rightarrow \SL(d,\Kb)$ is a representation. Then 
$\rho$ is $P_k$-Anosov if and only if 
there exists a $\rho$-equivariant, transverse, continuous, strongly dynamics preserving map 
$\xi=(\xi^k,\xi^{d-k}):\Lambda(\Gamma)\to\mathrm{Gr}_k(\mathbb K^d)\times\mathrm{Gr}_{d-k}(\mathbb K^d)$.
Furthermore, $\xi$ is the $P_k$-Anosov limit map.
\end{theorem}

In general, being $P_k$-Anosov is not an open condition in the space of representations of a geometrically finite group. For instance, consider the case where $\Gamma=\ip{g_1,g_2} \subset \PSL(2,\Rb)$ is a free group, $g_2$ is parabolic, and $\rho : \Gamma \rightarrow \SL(4,\Rb)$ is $P_1$-Anosov with 
$$
\rho(g_2) = \begin{pmatrix} 1 & 1 & 0 & 0 \\ 0 & 1 & 0 & 0 \\ 0 & 0 & 1 & 0 \\ 0 & 0 & 0 & 1 \end{pmatrix}.
$$
Next define a family of representations $\rho_t$ where $\rho_t(g_1) = \rho(g_1)$ and 
$$
\rho_t(g_2) = \begin{pmatrix} 1 & 1 & 0 & 0 \\ 0 & 1 & 0 & 0 \\ 0 & 0 & 1 & t \\ 0 & 0 & 0 & 1 \end{pmatrix}.
$$
Then $\rho_t$ is not $P_1$-Anosov for any $t \neq 0$ since the sequence $\{\rho_t(g_2)^n\}$ does not converge to a rank one element of $\Pb(\End(\Rb^4))$. 
Bowditch \cite[Sec. 5]{bowditch-deform} gave  an example
of a sequence $\{\rho_n\}$  of indiscrete representations of a free group on two generators into $\mathsf{SO}(4,1)$ which converge to a geometrically finite representation $\rho$, 
so that $\rho_n(\alpha)$ is parabolic if and only if $\rho(\alpha)$ is parabolic. However, in his example $\rho_n(\alpha)$ is not conjugate to $\rho(\alpha)$ for any $n$.
Bowditch \cite[Thm. 1.5]{bowditch-deform} also established a stability theorem for deformations of geometrically finite representations which preserve the structure
of  the Jordan decomposition of parabolic elements.

To account for these examples, we introduce the following subset of the representation variety. 
If $\rho:\Gamma\to\mathsf{SL}(d,\mathbb K)$ is a  representation of a geometrically finite Fuchsian group, let
$$\mathrm{Hom}_{\rm tp}(\rho)\subset\mathrm{Hom}(\Gamma,\SL(d,\Kb))$$
be the space of representations $\sigma:\Gamma\to\SL(d,\mathbb K)$ so that  if $\alpha\in\Gamma$
is parabolic, then $\sigma(\alpha)$ is conjugate to $\rho(\alpha)$. We obtain the following stability result for
type-preserving deformations. (In their preprint, Kapovich and Leeb  \cite{kapovich-leeb} suggest that such a stability
result  holds more generally.) Furthermore, we show that limit maps vary analytically. Combined with work of
Bray-Canary-Kao-Martone \cite{BCKM} this will imply that entropy and pressure intersection vary analytically over the cusped 
Hitchin component (which will be used crucially in the construction of the pressure metric).  The proof also allows us to see
that the $P_k$-Anosov limit maps are uniformly  H\"older in a neighborhood of a $P_k$-Anosov representation.

We say that $\{\rho_u\}_{u\in M}$ is a {\em $\mathbb K$-analytic family} of representations if $M$ is a $\mathbb K$-analytic
manifold and the map $u\to\rho_u$ is a $\mathbb K$-analytic map from $M$ into $\mathrm{Hom}(\Gamma,\mathsf{SL}(d,\mathbb K))$.

\begin{theorem}\label{thm: stability intro}
If $\Gamma\subset\mathsf{PSL}(2,\mathbb R)$ is geometrically finite and $\rho_0:\Gamma\to\mathsf{SL}(d,\mathbb K)$
is $P_k$-Anosov, then there exists an open neighborhood $\Oc$ of $\rho_0$ in $\mathrm{Hom}_{\rm tp}(\rho_0)$, so that
\begin{enumerate}
\item
If $\rho\in \Oc$, then $\rho$ is $P_k$-Anosov.
\item
There exists $\alpha>0$ so that if $\rho\in\Oc$, then its $P_k$-Anosov limit map $\xi_\rho$ is $\alpha$-H\"older.
\item
If  $\{\rho_u\}_{u\in M}$ is a $\mathbb K$-analytic family of representations in $\Oc$ and $z\in\Lambda(\Gamma)$
then the map from $M$ to  $\mathrm{Gr}_{k}(\mathbb K^d)\times\mathrm{Gr}_{d-k}(\mathbb K^d)$
given by $u\to \xi_{\rho_u}(z)$ is $\mathbb K$-analytic.
\end{enumerate}
\end{theorem}

When $\wh\Usf(\Gamma)$ is compact, stability follows from standard arguments in hyperbolic dynamics and the H\"older regularity of the boundary maps is a consequence of standard results, e.g. \cite[Cor. 5.19]{shub-book}. The non-compact case is more involved. Our key idea to prove stability is to observe that if $\rho_2 \in \mathrm{Hom}_{\rm tp}(\rho_1)$, then on each cusp there is a smooth conjugacy of the flows associated to $\rho_1$ and $\rho_2$. This is made precise in Equation~\eqref{eqn:smooth_conj} below. This essentially means that the two flows differ by a compact perturbation and thus, essentially, reduces to the compact base case. Our key idea to prove H\"older regularity is to first introduce certain ``canonical families of norms'' on the flow spaces, see Section~\ref{sec:canonical_norms_and_cusp_repn}. Then we prove that if the flow is contracting with respect to any family of norms, then the flow is contracting with respect to any canonical family of norms (see Theorem \ref{dynamics implies anosov} and Theorem \ref{thm: intro thm 1 body}). Finally, the canonical family of norms are well-behaved enough that we can adapt an argument from~\cite{ZZ2019} to prove H\"older regularity of the boundary maps directly.

We now discuss the applications of our results to cusped Hitchin representations, which was the original goal of our work.
Given an ordered basis $\mathcal B$ for $\mathbb R^d$, we say that a unipotent $A\in\mathsf{SL}(d,\mathbb R)$ is {\em totally positive} with
respect to $\mathcal B$, if its matrix with respect to $\mathcal B$ is upper triangular and all its minors (which are not
forced to be 0 by the fact that the matrix is upper triangular) are strictly positive. The set  $\mathcal U_{>0}(\mathcal B)$ of unipotent, totally
positive matrices with respect to $\mathcal B$ is a semi-group (see Lusztig \cite[Section 2.12]{lusztig}). Following Fock and Goncharov \cite{fock-goncharov},
we say that an ordered  $k$-tuple $(F_1,F_2,\ldots, F_k)$ of distinct flags in $\mathcal F_d$ is positive, if there exists an ordered basis $\mathcal B=(b_1,\ldots,b_d)$
for $\mathbb R^d$ so that $b_i\in F_1^i\cap F_k^{d-i+1}$ for all $i$, and there exists $u_2,\ldots,u_{k-1}\in U_{>0}(\mathcal B)$ so that
$F_i=u_{k-1}\cdots u_{i}F_k$ for all $i=2,\ldots,k-1$. Fock and Goncharov \cite{fock-goncharov} proved that positivity of a $n$-tuple is invariant
under cyclic permutation and reversal (also see Kim-Tan-Zhang \cite[Section 3.1 -- 3.3]{KTZ}).
If $X$ is a subset of $\mathbb{S}^1$ then a map $\xi:X\to\mathcal F_d$ is {\em positive} if whenever $(x_1,\ldots, x_n)$ is a cyclically ordered subset of
distinct points in $X$, then $(\xi(x_1),\ldots,\xi(x_n))$ is a positive $n$-tuple of flags. 

We say that a representation  $\rho:\Gamma\to\mathsf{SL}(d,\mathbb R)$ of a geometrically finite Fuchsian group is a  \emph{Hitchin representation}  if there exists a continuous positive
$\rho$-equivariant map $\xi_\rho:\Lambda(\Gamma)\to \mathcal F_d$. If $\Gamma$ is cocompact and torsion-free then these are exactly the Hitchin
representations of closed surface groups introduced by Hitchin \cite{hitchin} and further studied by Labourie \cite{labourie-invent}, while if $\Gamma$ is convex cocompact they
are the same as the Hitchin representations studied by Labourie and McShane \cite{labourie-mcshane}. We distinguish the case where
$\Gamma$ contains parabolic elements by calling these {\em cusped Hitchin representations}.  If $d=3$ and $S=\Gamma\backslash\mathbb H^2$ is
a finite area hyperbolic surface, then cusped Hitchin representations of $\Gamma$ are holonomy representations of strictly convex, finite area,
projective structures on $S$ (see Marquis \cite{marquis-finite area}). Further it follows from \cite[Thm. 3.3]{fock-goncharov} that holonomy
maps of geometrically finite projective surfaces in the sense of Crampon-Marquis  \cite{crampon-marquis} are also cusped Hitchin.
Fock and Goncharov  \cite{fock-goncharov} studied the a priori more general
class of representations which admit equivariant positive maps from the set $\Lambda_p(\Gamma)$ of fixed points of peripheral elements of 
$\Gamma$ into $\mathcal F_d$ when $\Gamma$ is not cocompact but has cofinite area. We show that all such type-preserving representations are in fact cusped Hitchin representations.
Our main result here is that cusped Hitchin representations are $P_k$-Anosov for all $k$.

\begin{theorem}\label{thm: Hitchin}
If $\Gamma\subset\mathsf{PSL}(2,\mathbb R)$ is geometrically finite and $\rho:\Gamma\to\mathsf{SL}(d,\mathbb R)$
is a Hitchin representation, then $\rho$ is irreducible and $P_k$-Anosov for all $k$. Moreover:
\begin{enumerate}
\item For all $k$ the map $x\mapsto\xi_\rho^{k}(x)$ is the $P_k$-Anosov limit map.
\item If $\alpha \in \Gamma$ is parabolic, then $\rho(\alpha)=\pm u$ for some unipotent $u\in\SL(d,\Rb)$ with a single Jordan block.
\item If $\gamma \in \Gamma$ is hyperbolic, then $\rho(\gamma)$ is loxodromic. 
\end{enumerate}
\end{theorem}

If we let $\widehat{\mathcal H}_d(\Gamma)$ denote the space of all  Hitchin representations of $\Gamma$ into $\mathsf{SL}(d,\mathbb R)$,
it is easy to see that $\widehat{\mathcal H}_d(\Gamma)$ is a real analytic manifold. In fact, one may use results of Fock and Goncharov \cite{fock-goncharov}
to show that the space $\mathcal H_d(\Gamma)$  of conjugacy classes of Hitchin representations is diffeomorphic to $\mathbb R^m$ (for some $m$).
(Marquis \cite{marquis-param} gives an explicit parametrization of $\mathcal H_3(\Gamma)$ as a topological
cell  when $S=\Gamma\backslash \mathbb H^2$ is a finite area hyperbolic surface.)

\subsection*{Comparison to other results:} As mentioned above, Kapovich-Leeb \cite{kapovich-leeb} and Zhu \cite{zhu-reldom} have previously developed theories of Anosov representations for 
relatively hyperbolic groups. Their work is based on extending characterizations of Anosov representations due to Kapovich-Leeb-Porti \cite{KLP} and Bochi-Potrie-Sambarino \cite{BPS} respectively. 
Theorem \ref{thm:  intro thm 1} immediately implies that a $P_k$-Anosov representation  $\rho$ is $P_k$-relatively
dominated, in the sense of Zhu \cite{zhu-reldom}. Theorem 9.4 in \cite{zhu-reldom} then implies that $\rho$ is $P_k$-relatively asymptotically
embedded and Remark 9.10 in \cite{zhu-reldom} implies that $\rho$ is $P_k$-relatively uniform RCA in the sense of Kapovich-Leeb \cite{kapovich-leeb}.

\begin{corollary}
\label{P_k rel dom}
Suppose that  $\Gamma\subset\mathsf{PSL}(2,\mathbb R)$ is geometrically finite  and
$\rho:\Gamma\to\mathsf{SL}(d,\mathbb K)$ is a $P_k$-Anosov representation.
Then $\rho$ is $P_k$-relatively dominated, $P_k$-relatively asymptotically embedded
and $P_k$-relatively uniform RCA.
\end{corollary}


Further, it follows from  \cite[Thm. C]{zhu-reldom2}
and our Theorems \ref{thm: intro thm 1} and \ref{thm:equivalent_to_Anosov}, that a representation of a geometrically finite Fuchsian group  $\Gamma$ is $P_k$-Anosov if and only if it is $P_k$-relatively dominated with respect to $\mathcal P$.

Kapovich and Leeb \cite{kapovich-leeb} mention that they can show that a cusped Hitchin representation of a geometrically finite Fuchsian is relatively Anosov in their sense.

In a subsequent preprint, Filip \cite{filip} introduces the class of $\log$-Anosov representation of Fuchsian lattices, which also agrees with the class of cusped
Anosov representations. He introduces adapted metrics, which correspond to our canonical norms. He shows that monodromy maps of certain variations of
Hodge structures on finite area Riemann surfaces are $\log$-Anosov and uses this result to study their properties. In particular many monodromy maps
coming from hypergeometic differential equations can be analyzed in this manner, including those previously studied by Brav and Thomas \cite{brav-thomas} and by
Filip and Fougeron \cite{filip-fougeron}.

\subsection*{Acknowledgements:} The authors thank the referee for their careful reading of the original manuscript and helpful suggestions which improved the exposition. This material is based upon work supported by the National Science Foundation
under Grant No. DMS-1928930 while the first author participated in a program hosted
by the Mathematical Sciences Research Institute in Berkeley, California, during
the Fall 2020 semester.

\section{Preliminaries}
In this section, we recall some preliminary facts and introduce notation that will be used throughout this paper. 

\subsection{Hyperbolic 2-space} In this paper we will  identify $\Hb^2$ with the Poincar\'e upper half plane model. 
For any $v\in T^1\Hb^2$, let $r_v:\Rb\to\Hb^2$ denote the unit speed geodesic with $r_v'(0)=v$ and let 
$$v^+:=\lim_{t\to+\infty} r_v(t)\in\partial \Hb^2\qquad\mathrm{and}\qquad v^-:=\lim_{t\to-\infty} r_v(t)\in\partial \Hb^2$$
denote its limit points in $\partial\Hb^2 = \Rb \cup \{\infty\}$. We also let $\phi_t:T^1\Hb^2\to T^1\Hb^2$ denote the geodesic flow, i.e.
\begin{align*}
\phi_t(v) = r_v^\prime(t)
\end{align*}
for all $v \in T^1 \Hb^2$ and $t \in \Rb$. 

If $\{g_n\}$ is a sequence in $\PSL(2,\Rb)$, we say that $g_n$ {\em converges} to $x\in\partial\Hb^2$ if $g_n(z)\to x$
for some (any)  $z\in\mathbb H^2$. We often simply write $g_n\to x$.

\subsection{Geometrically finite Fuchsian groups:}We say that  $\Gamma\subset\mathsf{PSL}(2,\mathbb R)$ is {\em geometrically finite} if it is discrete, finitely generated and non-elementary 
(i.e. does not contain a cyclic subgroup of finite index). Let $\Lambda(\Gamma)\subset\partial\mathbb H^2$ denote
its limit set. 
Then let $\Usf(\Gamma)$ denote the minimal, non-empty, closed, $\phi_t$-invariant, $\Gamma$-invariant subset of $T^1\Hb^2$, i.e.
\begin{align*}
\Usf(\Gamma)
&= \left\{ v \in T^1 \Hb^2 : v^+, v^- \in \Lambda(\Gamma)\right\}. 
\end{align*}
Let $\wh{\Usf}(\Gamma):=\Gamma\backslash\Usf(\Gamma)$, and note that $\phi_t$ descends to a flow on $\wh{\Usf}(\Gamma)$, which we also denote by $\phi_t$. 
If $S=\Gamma\backslash\mathbb H^2$, then $\wh{\Usf}(\Gamma)$ is the non-wandering portion of $T^1S$ and its orbits are the complete geodesics which remain entirely in
the convex core of $S$.

\begin{definition} 
\begin{enumerate}
\item
If $p\in\Lambda(\Gamma)$ is fixed by a parabolic element of $\Gamma$, let $\mathrm{st}_\Gamma(p)=\{\gamma\in\Gamma\ |\ \gamma(p)=p\}$.
We call $\mathrm{st}_\Gamma(p)$ a {\em cuspidal subgroup} of $\Gamma$.
An open horodisk $H$ based at $p$ is a {\em precisely invariant horodisk} for $\mathrm{st}_\Gamma(p)$ if whenever $\gamma\in\Gamma$ and
$\gamma(H)\cap H$ is non-empty, then $\gamma\in\mathrm{st}_\Gamma(p)$. In this case, $C=\mathrm{st}_\Gamma(p)\backslash H$
is an {\em embedded cusp neighborhood}.
\item
A set $\mathcal C=\{C_1,\ldots,C_r\}$ of disjoint embedded cusp neighborhoods in $S=\Gamma\backslash\mathbb H^2$ is {\em full} if any curve on $S$ which is
represented by a parabolic element of $\Gamma$ is freely homotopic into some $C_i$.
\end{enumerate}
\end{definition} 

If $C=\mathrm{st}_\Gamma(p)\backslash H$ is an embedded cusp neighborhood, then we set
$${\Usf}(\Gamma)_{H}=\{v\in\Usf(\Gamma) |\ r_v(0)\in H \}\ \qquad\mathrm{and}\qquad
\wh{\Usf}(\Gamma)_{C}=\mathrm{st}_\Gamma(p)\backslash\Usf(\Gamma)_H.$$
If $\mathcal C$ is a full collection of cusp neighborhoods for $S$, then we set
$$\wh{\Usf}(\Gamma)_{\mathcal C}=\bigcup_{C\in\mathcal C} \wh{\Usf}(\Gamma)_{C}.$$
Notice that its complement $\wh{\Usf}(\Gamma)_{\mathcal C}^c$ is compact. We will sometimes informally refer
to $\wh{\Usf}(\Gamma)_{\mathcal C}$ as the thin part of the geodesic flow and its complement as the thick part.

\subsection{The representation theory of $\mathsf{SL}(2,\mathbb R)$:}\label{sec:repn_thy_of_SL2}
For $d\geq 1$,  let $\tau_d:\SL(2,\Rb)\to\SL(d,\Kb)$ denote the standard irreducible representation. Explicitly, $\tau_1 \equiv 1$ and if $d \geq 2$, then $\tau_d=i_d\circ \bar{\tau}_d$, where 
\begin{itemize}
\item $\bar{\tau}_d:\SL(2,\Rb)\to\SL(\mathrm{Sym}^{d-1}(\Rb^2))$ is the representation induced by the linear $\SL(2,\Rb)$-action on $\mathrm{Sym}^{d-1}(\Rb^2)$ given by
\[\gamma\cdot (e_1^ke_2^{d-1-k})\mapsto (\gamma\cdot e_1)^k(\gamma\cdot e_2)^{d-1-k}\]
for all $k=0,\dots,d-1$, and 
\item $i_d:\SL(\mathrm{Sym}^{d-1}(\Rb^2))\to\SL(d,\Rb)$ is the isomorphism induced by the linear isomorphism \hbox{$\mathrm{Sym}^{d-1}(\Rb^2))\simeq \Rb^d$} that identifies $e_1^ke_2^{d-1-k}\in\mathrm{Sym}^{d-1}(\Rb^2))$ with $e_{k+1}\in\Rb^d$ for all $k=0,\dots,d-1$.
\end{itemize}
One may also regard $\tau_d$ as a representation into $\mathsf{SL}(d,\mathbb C)$.
Let $\xi_d:\partial \Hb^2\simeq\mathbb{RP}^1\to\Fc(\Kb^d)$ be the map defined by $[ae_1+be_2]\mapsto F$, where $F^{(j)}=\Span_{\Kb}(f_1,\dots,f_j)$ and
\[f_j=\sum_{k=1}^{d+1-j}{d-j\choose k-1}a^{d+1-j-k}b^{k-1}e_k.\]
It is straightforward to verify that $\xi_d$ is a continuous, $\tau_d$-equivariant, strongly dynamics preserving map.  We call $\xi_d$ the {\em Veronese embedding}
associated to $\tau_d$.

One can compute that $\tau_d\begin{pmatrix}
1&a\\
0&1
\end{pmatrix}$ is the upper triangular matrix in $\SL(d,\Rb)$ given by 
\begin{equation}\label{eqn: Pascal}
\left(\tau_d\begin{pmatrix}
1&a\\
0&1
\end{pmatrix}\right)_{k,j}={j-1\choose k-1}a^{j-k}
\end{equation}
if $1\leq k\leq j\leq d$. Furthermore, it is well-known that the $d\times d$ upper triangular matrix given by (\ref{eqn: Pascal}) is unipotent and totally positive with respect to the standard basis of $\Rb^d$ when $a>0$. From this, it follows easily that $\xi_d$ is a positive map.

\subsection{Eigenvalues, singular values, and the (multiplicative) Jordan-Chevalley decomposition}
Given an element $g\in\SL(d,\Kb)$, let 
\begin{align*}
\sigma_1(g)\geq\dots\geq\sigma_d(g)
\end{align*}
denote the  singular values of $g$ and let 
\begin{align*}
\lambda_1(g)\geq\dots\geq\lambda_d(g)
\end{align*}
denote the absolute values of the (generalized) eigenvalues of $g$. Also, let $g = g_{ss} g_u = g_u g_{ss}$ denote the Jordan-Chevalley decomposition, that is $g_{ss}$ is semisimple, $g_u$ is unipotent, and $g_{ss}, g_u$ commute. We say that $g$ is \emph{elliptic} if it is semisimple and $\lambda_j(g) = 1$ for all $j$.  Notice that if $g$ is elliptic, then the cyclic group $\ip{g}$ generated by $g$ is relatively compact in $\SL(d,\Kb)$. 

If $g\in\SL(d,\Kb)$ and $\sigma_k(g)>\sigma_{k+1}(g)$, let $U_k(g)\in\Gr_k(\Kb^d)$ denote the subspace spanned by the $k$ major axes of the ellipse $g\left(\mathbb{S}^{d-1}\right)$, that is $U_k(g) = g \ip{e_1,\dots, e_k}$. 
The following lemma relates the singular values along a sequence in $\SL(d,\Kb)$ to the action of this sequence on the associated Grassmannian. We omit the proof as it is standard.

\begin{lemma}\label{lem: basic singular value} Suppose $V_0\in\Gr_k(\Kb^d)$ and $W_0\in\Gr_{d-k}(\Kb^d)$. If $\{g_n\}_{n\geq 1}$ is a sequence in $\SL(d,\Kb)$, then the following are equivalent:
\begin{enumerate}
\item There are open sets $\Oc\subset\Gr_k(\Kb^d)$ and $\Oc^\prime \subset \Gr_{d-k}(\Kb^d)$, such that $g_n(V)\to V_0$ for all $V\in\Oc$ and $g_n^{-1}(W) \to W_0$ for all $W \in \Oc^\prime$. 
\item  $g_n(V) \to V_0$ for all $V\in\Gr_k(\Kb^d)$ transverse to $W_0$. 
\item $\displaystyle\lim_{n\to\infty}\frac{\sigma_k(g_n)}{\sigma_{k+1}(g_n)}=\infty$, $\displaystyle\lim_{n \to \infty} U_k(g_n)=V_0$, and $\displaystyle\lim_{n \to \infty} U_{d-k}(g_n^{-1}) = W_0$. 
\end{enumerate}
Moreover, if $g_n=g^n$ for all $n$, then $g$ is $P_k$-proximal if $V_0\oplus W_0=\Kb^d$, and weakly unipotent if $V_0\subset W_0$ or $W_0\subset V_0$.
\end{lemma}

We will also make use of the following elementary calculations, which we recall without proof.

\begin{lemma}
\label{singular value facts}
\begin{enumerate}
\item
If $\u2=\begin{pmatrix} 1 & 1\\ 0 & 1\\ \end{pmatrix}$, then for all $d\ge 1$ there exists $c=c(d)>0$ so that
$$\frac{1}{c}\le \frac{\sigma_k(\tau_d(\u2^n))}{n^{d+1-2k}}\le c$$
for all $n\in\mathbb N$ and $k\in\{1,\ldots,d\}$.
\item
If $A,B\in\mathsf{SL}(d,\mathbb K)$ and $C=BAB^{-1}$, then
$$\sigma_k(C)\le \frac{\sigma_1(B)}{\sigma_d(B)}\sigma_k(A)$$
for all $k\in\{1,\ldots,d\}$.
\end{enumerate}
\end{lemma}

\section{Anosov representations into $\mathsf{SL}(d,\mathbb K)$}\label{sec: Anosov representations}

Following Labourie \cite{labourie-invent} and Guichard-Weinhard \cite{guichard-wienhard} we define Anosov representations for geometrically finite Fuchsian groups into
$\mathsf{SL}(d,\mathbb K)$. 
In Appendix \ref{app: general}, we will extend our definition to the setting of all semi-simple Lie groups with finite center.

If $\Gamma \subset \PSL(2,\Rb)$ is a geometrically finite group and $\rho:\Gamma\to\mathsf{SL}(d,\mathbb K)$ is a representation,  let
$$E_\rho=\Usf(\Gamma)\times \mathbb K^d\quad\mathrm{and}\quad\widehat E_\rho=\Gamma\backslash\Usf(\Gamma)\times \mathbb K^d$$
where $\gamma\in\Gamma$ acts on the first factor by $\gamma$ and on the second factor by $\rho(\gamma)$. 
The projection map $p:E_\rho\to \Usf(\Gamma)$ descends to a vector bundle
$$\hat p:\widehat E_\rho\to\wh{\Usf}(\Gamma)$$
which is called the {\em flat bundle associated to $\rho$}. The geodesic flow on $\Usf(\Gamma)$ extends to a flow on \hbox{$E_\rho = \Usf(\Gamma) \times \Kb^d$} whose action is trivial on the second factor. This in turn descends to a flow on $\wh{E}_\rho$ which covers the geodesic flow on $\wh{\Usf}(\Gamma)$. We use $\phi_t$ to denote all four of these flows.

We say that a continuous map $\xi=(\xi^k,\xi^{d-k}):\Lambda(\Gamma)\to \mathrm{Gr}_k(\mathbb K^d)\times\mathrm{Gr}_{d-k}(\mathbb K^d)$
is  
\begin{enumerate}
\item \emph{$\rho$-equivariant} if $\rho(\gamma) \circ \xi = \xi \circ \gamma$ for all $\gamma \in \Gamma$, 
\item \emph{transverse} if $\xi^k(x) \oplus \xi^{d-k}(y)=\Kb^d$ for all distinct $x,y \in \Lambda(\Gamma)$.
\end{enumerate}
Such a map induces a continuous $\phi_t$-invariant spitting
$$E_\rho=\Theta^k\oplus\Xi^{d-k}\quad \mathrm{where}\quad \Theta^k|_v=\xi^k(v^+)\quad \mathrm{and}\ \ \Xi^{d-k}|_v=\xi^{d-k}(v^-)$$
which descends to a continuous $\phi_t$-invariant splitting 
$$\wh E_\rho=\wh \Theta^k\oplus\wh\Xi^{d-k}.$$

\begin{definition}\label{defn:anosov_specific}  If $\Gamma \subset \PSL(2,\Rb)$ is a geometrically finite group and $k\in\{1,\ldots,d-1\}$, 
a representation $\rho : \Gamma \rightarrow \mathsf{SL}(d,\mathbb K)$ is \emph{$P_k$-Anosov} if:
\begin{enumerate}
\item There exists a  $\rho$-equivariant, continuous, transverse map 
$$\xi=(\xi^k,\xi^{d-k}):\Lambda(\Gamma)\to \mathrm{Gr}_k(\mathbb K^d)\times\mathrm{Gr}_{d-k}(\mathbb K^d)$$
which induces a splitting $\widehat E_\rho=\widehat\Theta^k\oplus\widehat\Xi^{d-k}$.
\item For some  family of norms $\norm{\cdot}_v$ on the fibers of $\wh E_\rho$ there exists $C>1$ and $c > 0$ so that 
\begin{align*}
\frac{\norm{\phi_{t}(Y)}_{\phi_{t}(v)}}{\norm{\phi_t(Z)}_{\phi_t(v)}} \leq Ce^{-ct} \frac{\norm{Y}_{v}}{\norm{Z}_{v}}
\end{align*}
for all $t > 0$, $v \in \wh{\Usf}(\Gamma)$,  $Y\in\wh\Theta^k|_v$ and non-zero  $Z\in\wh\Xi^{d-k}|_v$.
\end{enumerate}
We refer to $\xi$ as  a $P_k$-Anosov \emph{limit map of $\rho$}.
\end{definition}

\begin{remark} Notice that a continuous family of norms $\norm{\cdot}_{v \in \wh{\Usf}(\Gamma)}$ on the fibers of $\wh E_\rho$ lifts to a continuous family of norms $\norm{\cdot}_{v \in \Usf(\Gamma)}$ on $\Kb^d$ which is $\rho$-equivariant in the following sense: if $v \in \Usf(\Gamma)$, $Y \in \Kb^d$ and $\gamma \in \Gamma$, then 
\begin{equation}
\label{eqn:equivariance_of_norms}
\norm{\rho(\gamma)(Y)}_{\gamma (v)} = \norm{Y}_v.
\end{equation}
Conversely, any continuous family of norms satisfying Equation~\eqref{eqn:equivariance_of_norms} descends to a continuous family of norms on the fibers of $\wh{E}_\rho$. 
\end{remark} 

If one prefers a bundle-theoretic definition, one can consider the vector bundle $\mathrm{Hom}(\wh\Xi^{d-k},\wh\Theta^k)$ over $\wh{\Usf}(\Gamma)$.
Notice that since the splitting is flow-invariant, $\phi_t$ induces a flow on $\mathrm{Hom}(\wh\Xi^{d-k},\wh\Theta^k)$ given by 
\begin{align*}
f \mapsto \phi_t \circ f \circ \phi_{-t},
\end{align*}
with some abuse of notation we use $\phi_t$ to denote this flow. 
Moreover, any  norm on the fibers of  $\wh E_\rho$ induces an operator norm on  $\mathrm{Hom}(\wh\Xi^{d-k},\wh\Theta^k)$. 
We say that the flow $\phi_t$ on  $\mathrm{Hom}(\wh\Xi^{d-k},\wh\Theta^k)$ is {\em uniformly contracting}  if there exists $C,c>0$ so that
$$\norm{\phi_t(f)}_{\phi_t(v)}\le Ce^{-ct}\norm{f}_v$$ 
for all $t>0$, $v\in\wh\Usf(\Gamma)$ and $f\in\mathrm{Hom}(\wh\Xi^{d-k},\wh\Theta^k)_v$.
One may easily check that $\rho$ is $P_k$-Anosov if and only if there exists a $\rho$-equivariant, continuous transverse map 
$\xi=(\xi^k,\xi^{d-k}):\Lambda(\Gamma)\to \mathrm{Gr}_k(\Kb^d)\times\mathrm{Gr}_{d-k}(\mathbb K^d)$ so that 
the flow is uniformly contracting on $\mathrm{Hom}(\wh\Xi^{d-k},\wh\Theta^k)$ with respect to an operator norm associated to a continuous
family of norms on $\wh E_\rho$.
(The details of this equivalence are worked out carefully in the proof of \cite[Prop. 2.3]{BCLS}.) Moreover, by duality, $\phi_t$ is uniformly contracting on
$\mathrm{Hom}(\wh\Xi^{d-k},\wh\Theta^k)$ if and only if the flow, also called $\phi_t$, on $\mathrm{Hom}(\wh\Theta^{k},\wh\Xi^{d-k})$  is uniformly expanding,
i.e. the inverse flow $\phi_{-t}$ is uniformly contracting. We record these observations for future use.

\begin{proposition}\label{prop: Hom bundle Anosov}
Suppose that $\Gamma \subset \PSL(2,\Rb)$ is a geometrically finite group, $\rho : \Gamma \rightarrow \mathsf{SL}(d,\mathbb K)$ is a representation and $k\in\{1,\ldots,d-1\}$. Then the
following are equivalent:
\begin{enumerate}
\item
$\rho$ is $P_k$-Anosov,
\item
There exists a $\rho$-equivariant, continuous transverse map 
$$\xi=(\xi^k,\xi^{d-k}):\Lambda(\Gamma)\to \mathrm{Gr}_k(\mathbb R^d)\times\mathrm{Gr}_{d-k}(\mathbb K^d)$$
so that 
the flow is uniformly contracting on $\mathrm{Hom}(\wh\Xi^{d-k},\wh\Theta^k)$ with respect to an operator norm
induced by a norm on $\wh E_\rho$.
\item
There exists a $\rho$-equivariant, continuous transverse map 
$$\xi=(\xi^k,\xi^{d-k}):\Lambda(\Gamma)\to \mathrm{Gr}_k(\mathbb R^d)\times\mathrm{Gr}_{d-k}(\mathbb K^d)$$
so that 
the flow is uniformly expanding on $\mathrm{Hom}(\wh\Theta^{k},\wh\Xi^{d-k})$ with respect to an operator norm
induced by a norm on $\wh E_\rho$.
\end{enumerate}
\end{proposition}

As an immediate corollary we obtain:

\begin{corollary}
Suppose that $\Gamma \subset \PSL(2,\Rb)$ is a geometrically finite group, $\rho : \Gamma \rightarrow \mathsf{SL}(d,\mathbb K)$ is a representation and $k\in\{1,\ldots,d-1\}$. Then
$\rho$ is $P_k$-Anosov if and only if $\rho$ is $P_{d-k}$-Anosov.
\end{corollary}

\begin{remark} If we wish to allow $\Gamma$ to be a geometrically finite subgroup of $\mathsf{PGL}(2,\mathbb R)$, then we simply consider 
$\Gamma_0=\Gamma\cap\mathsf{PSL}(2,\mathbb R)$ and say that $\rho:\Gamma\to \mathsf{SL}(d,\mathbb R)$ is $P_k$-Anosov if and only if 
$\rho|_{\Gamma_0}$ is \hbox{$P_k$-Anosov}. With this definition, all of our results for geometrically finite Fuchsian groups  remain true for 
geometrically finite subgroups of $\mathsf{PGL}(2,\mathbb R)$.
\end{remark}

\subsection{Canonical norms and cusp representations}\label{sec:canonical_norms_and_cusp_repn} 
It will be useful to construct certain ``canonical'' norms on the fibers of the vector bundles $\wh E_\rho$, and hence also on $\Hom(\wh\Xi^{d-k},\wh\Theta^k)$, that are well behaved on the cusps. Later we will show that when $\rho$ is $P_k$-Anosov, the flow $\phi_t$ on $\Hom(\wh\Xi^{d-k},\wh\Theta^k)$ is uniformly contracting for any canonical norm (see Corollary \ref{cor: contraction canonical norm}).

The crucial property of our canonical norms is that they have a standard form over the thin part of the geodesic flow. In order to
describe this standard form we will use the following result about representations of $\SL(2,\Rb)$. Recall that $g = g_{ss} g_u$ denotes the multiplicative Jordan-Chevalley decomposition
of \hbox{$g \in \SL(d,\Kb)$}. 

\begin{proposition}\label{prop:building_repn}(see Appendix~\ref{appendix: cusp repn}) If $g \in \SL(d,\Kb)$ is weakly unipotent, then there exists a representation $\Psi: \SL(2,\Rb) \rightarrow \SL(d,\Kb)$ where $\Psi\left(\begin{pmatrix} 1 & 1\\ 0 & 1\\ \end{pmatrix}\right) = g_u$ and $g_{ss}$ commutes with the elements of $\Psi(\SL(2,\Rb))$. \end{proposition} 

Proposition~\ref{prop:building_repn} follows easily from the Jordan normal form of a weakly unipotent matrix and we delay the proof until Appendix~\ref{appendix: cusp repn}. 

A representation $\rho:\Gamma\to\mathsf{SL}(d,\mathbb K)$ of a geometrically finite Fuchsian group is \emph{type-preserving} if $\rho$ sends every parabolic element in $\Gamma$ to a weakly unipotent element in $\SL(d,\Kb)$. If $\rho$ is a type preserving representation and $\alpha\in \Gamma$ is parabolic, then we say that a representation $\Psi:\mathsf{SL}(2,\mathbb R)\to\SL(d,\Kb)$ is a {\em cusp representation} for $\alpha$ and $\rho(\alpha)$ if 
\begin{enumerate}
\item $\Psi(\wt{\alpha})=\rho(\alpha)_{u}$, where $\wt{\alpha}$ is the (unique) unipotent lift of  $\alpha$ to $\mathsf{SL}(2,\mathbb R)$ and 
\item $\rho(\alpha)_{ss}$ commutes with the elements of $\Psi(\mathsf{SL}(2,\Rb))$.
\end{enumerate}
Proposition~\ref{prop:building_repn} implies that cusp representations always exist.

Suppose that $\Gamma$ is a geometrically finite Fuchsian group, $\rho:\Gamma\to \mathsf{SL}(d,\mathbb K)$ is type-preserving and 
$C=\langle\alpha\rangle\backslash H$ is an embedded cusp neighborhood where $\alpha\in\Gamma$ is parabolic. Further suppose that $\Psi$ is a cusp representation for $\alpha$ and $\rho(\alpha)$, and $\norm{\cdot}_v$ is a family of norms of the fibers $T^1 \Hb^2 \times \Kb^d \rightarrow T^1 \Hb^2$ such that 
\begin{enumerate}
\item each $\norm{\cdot}_v$ is $\rho(\alpha)_{ss}$-invariant 
\item $\norm{\Psi(g) Z}_{g(v)} = \norm{Z}_v$ for all $g \in \SL(2,\Rb)$, $Z \in \Kb^d$ and $v \in T^1 \Hb^2$. 
\end{enumerate}
Such families are easy to construct: the group $K:=\overline{\ip{ \rho(\alpha)_{ss}, \Psi(-\id_2)}}$ is abelian and compact, so there exists a norm $\norm{\cdot}_0$ which is $K$-invariant, then if we fix some $v_0 \in T^1 \Hb^2$, the family of norms defined by $\norm{Z}_{g(v_0)} := \norm{\Psi(g)^{-1}Z}_0$ has the desired properties. 

Also, for such a family of norms
$$
\norm{\rho(\alpha) Z}_{\alpha(v)} = \norm{ \Psi(\wt{\alpha}) \rho(\alpha)_{ss}Z}_{\alpha(v)} = \norm{Z}_v
$$
for all $Z \in \Kb^d$ and $v \in T^1 \Hb^2$. So this family descends (and restricts) to a family of norms on the fibers of $\wh E_\rho$ over $\wh{\Usf}(\Gamma)_C$
which we call a {\em canonical family of norms on the cusp neighborhood} $C$. Observe that if $C'$ is an embedded cusp neighborhood properly contained in $C$, then $\wh\Usf(\Gamma)_C-\wh\Usf(\Gamma)_{C'}$ is precompact, and a canonical family of norms on the fibers of $\wh E_\rho$ over $\wh{\Usf}(\Gamma)_C$ is determined completely by its restriction to the fibers over $\wh\Usf(\Gamma)_C-\wh\Usf(\Gamma)_{C'}$.

\begin{definition} 
Let $\rho:\Gamma\to \SL(d,\Kb)$ be a type preserving representation of a geometrically finite Fuchsian group. A continuous family of norms $\norm{\cdot}_{v}$ on the fibers of $\widehat E_\rho$
is \emph{canonical}  if there exists a full collection $\mathcal C$ of embedded cusp neighborhoods for $\Gamma$, such that for all  $C\in\mathcal C$, the
restriction of the family of norms to the fibers over $\wh{\Usf}(\Gamma)_C$ is a canonical family of norms on $C$.
\end{definition}

It is straightforward to construct a canonical family of norms on the fibers of $\widehat E_\rho$. 
One first chooses a full collection of embedded cusp neighborhoods $\mathcal C$. For each $C\in\mathcal C$, one chooses a canonical norm on the cusp neighborhood $C$. One then chooses any continuous norm on a compact neighborhood of the fibers over the thick part
$\wh{\Usf}(\Gamma)-\wh{\Usf}(\Gamma)_{\mathcal C}$.
One may then use a cut off function to interpolate between the norms on their interface and 
obtain a family of norms on all of $\widehat E_\rho$ which is canonical with respect to a full collection of cusp neighborhoods contained in $\mathcal C$.

We  observe that  there are uniform upper and lower bounds on the growth rate of a canonical norm with respect to the (lift of the) geodesic flow $\phi_t$.

\begin{lemma}\label{obs:distortion_of_norms_holder} 
Suppose that $\Gamma \subset \PSL(2,\Rb)$ is a geometrically finite group, $\rho:\Gamma\to\mathsf{SL}(d,\mathbb K)$ is type-preserving
and $\norm{\cdot}_{v}$ is a canonical family of norms on the bundle $\wh E_\rho$. There exist $C_0>1, c_0> 0$ so that
if $v\in \wh\Usf(\Gamma)$, $t\in\Rb$ and $Z \in \wh E_\rho|_{v}$, then
\begin{align*}
\frac{1}{C_0} e^{-c_0 |t|} \norm{Z}_{v} \leq \norm{\phi_t(Z)}_{\phi_t(v)} \leq C_0 e^{c_0 |t|} \norm{Z}_{v}.  
\end{align*}
\end{lemma}

\begin{proof} For $v \in \wh\Usf(\Gamma)$, let
\begin{align*}
f(v) = \max\left\{ \norm{\phi_t(Z)}_{\phi_t(v)} : \abs{t} \leq 1, Z \in  \wh E_\rho|_{v}, \norm{Z}_v=1 \right\}.
\end{align*}
Let $\Cc$ be a full collection of embedded cusp neighborhoods such that for all $C\in\Cc$, $\norm{\cdot}_v$ restricted to the fibers over each component of $\wh\Usf(\Gamma)_C$ is a canonical family of norms on $C$. For each $C\in\Cc$, let $C'$ be an embedded cusp neighborhood such that 
$$
\bigcup_{\abs{t} \leq 1} \phi_t(C^\prime) \subset C.
$$
Then $f$ is constant on each $\wh\Usf(\Gamma)_{C^\prime}$. Further, 
\[K=\wh\Usf(\Gamma)-\bigcup_{C\in\Cc}\wh\Usf(\Gamma)_{C'}\]
is compact. Hence
\begin{align*}
C_0 :=\sup_{v\in \wh{\Usf}(\Gamma)} f(v)=\max_{v \in K} f(v)
\end{align*}
is finite. Set $c_0=\log C_0$. For any $t\in\Rb$, let $n$ be the largest integer such that $n\leq |t|$. If $t>0$, then
\begin{align*}
\norm{\phi_t(Z)}_{\phi_t(v)} \leq  C_0^n\norm{\phi_{t-n}(Z)}_{\phi_{t-n}(v)}\le e^{c_0n}C_0\norm{Z}_{v}\leq C_0 e^{c_0 |t|} \norm{Z}_{v}.
\end{align*}
On the other hand, if $t<0$, then
\begin{align*}
\norm{\phi_t(Z)}_{\phi_t(v)} \leq  C_0^n\norm{\phi_{t+n}(Z)}_{\phi_{t+n}(v)}\le e^{c_0n}C_0\norm{Z}_{v}\leq C_0 e^{c_0 |t|} \norm{Z}_{v}.
\end{align*}
This proves the required upper bound. The lower bound is similar.
\end{proof}

\begin{remark} We will not need this for our work, but with a bit more effort, one can show that any two canonical families of norms
on $\wh E_\rho$ are bilipschitz.
\end{remark}

\section{Basic properties of Anosov representations}

In this section we prove Theorem~ \ref{thm: intro thm 1} which we restate here. 

\begin{theorem}\label{thm: intro thm 1 body}
If $\Gamma\subset\mathsf{PSL}(2,\mathbb R)$ is a geometrically finite group and $\rho:\Gamma\to\mathsf{SL}(d,\mathbb K)$
is $P_k$-Anosov, then
\begin{enumerate}
\item For any $z_0 \in \Hb^2$, there exists $A,a>1$ so that if $\gamma\in\Gamma$, then
$$\frac{1}{A}\exp\left(\frac{1}{a}d_{\Hb^2}(z_0,\gamma (z_0))\right)\leq \frac{\sigma_k(\rho(\gamma))}{\sigma_{k+1}(\rho(\gamma))} \leq A\exp\left(a d_{\Hb^2}(z_0,\gamma (z_0))\right).$$
\item
There exists $B,b>1$ so that if $\gamma\in\Gamma$, then
$$\frac{1}{B}\exp\left( \frac{1}{b}\ell(\gamma)\right)\leq \frac{\lambda_k(\rho(\gamma))}{\lambda_{k+1}(\rho(\gamma))} \leq B\exp\left( b\ell(\gamma)\right)$$
where  $\ell(\gamma)$ is the translation length of $\gamma$ on $\mathbb H^2$.
\item
The $P_k$-Anosov  limit map $\xi_\rho$ is strongly dynamics-preserving and unique. In particular, if $\alpha\in\Gamma$ is parabolic, then $\rho(\alpha)$ is weakly unipotent,
while if $\gamma\in\Gamma$ is hyperbolic, then $\rho(\gamma)$ is $P_k$-proximal.
 \item
If $z_0\in\mathbb H^2$ and $x_0$ is a point in the symmetric space $X_d(\mathbb K)$ for $\mathsf{SL}(d,\mathbb K)$, then
the orbit map $\tau_\rho:\Gamma(z_0)\to X_d(\mathbb K)$ given by $\tau_\rho(\gamma(z_0))=\rho(\gamma)(x_0)$ is a quasi-isometric embedding.
\end{enumerate}
\end{theorem}

Before proving the theorem we note the following consequences which will be useful in \cite{BCKM}.

\begin{corollary}
\label{BCKM properties}
Suppose that $\Gamma\subset\mathsf{PSL}(2,\mathbb R)$ is a geometrically finite group and $\rho:\Gamma\to\mathsf{SL}(d,\mathbb K)$ is $P_k$-Anosov with 
$P_k$-Anosov limit map $\xi_\rho$.
\begin{enumerate}
\item
If $\{\gamma_n\}$ is a sequence in $\Gamma$ with $\gamma_n\to x\in\Lambda(\Gamma)$, then
$\lim U_k(\rho(\gamma_n))=\xi_\rho(x)$.
\item
If $\alpha\in\Gamma$ is parabolic and $j\in \{1,\ldots,d\}$, then there exists an integer $c(j,\alpha)$ and 
$C_j>1$ so that
$$\frac{1}{C_j}\le \frac{\sigma_j(\rho(\alpha^n))}{n^{c(j,\alpha)}}\le C_j \quad \text{for all } n \in \Nb.$$
Moreover, 
$$c(k,\alpha)-c(k+1,\alpha)>0.$$
\end{enumerate}
\end{corollary}

\begin{proof} Property (1) in Corollary~\ref{BCKM properties} is known as the $P_k$-Cartan property and is an immediate consequence of the fact that
$\xi_\rho$ is strongly dynamics-preserving and Lemma \ref{lem: basic singular value}. 

If $\alpha\in\Gamma$ is parabolic, then part (3) of Theorem \ref{thm: intro thm 1} implies that $\rho(\alpha)$ is weakly unipotent. So the group $\overline{\ip{ \rho(\alpha)_{ss}}}$ is compact. Further, there exists $\{d_1,\ldots,d_m\}$, so that $\rho(\alpha)_u$ is conjugate to $\oplus_{i=1}^m\tau_{d_i}\left(\begin{bmatrix} 1 & 1\\ 0 & 1\end{bmatrix}\right)$. Thus there exists a constant $C_0 > 1$ such that 
$$
\frac{1}{C_0} \sigma_j\left(\oplus_{i=1}^m\tau_{d_i}\left(\begin{bmatrix} 1 & n\\ 0 & 1\end{bmatrix}\right)\right) \leq \sigma_j(\rho(\alpha^n)) \leq C_0 \sigma_j\left(\oplus_{i=1}^m\tau_{d_i}\left(\begin{bmatrix} 1 & n\\ 0 & 1\end{bmatrix}\right)\right)
$$
for all $n \in \Nb$ and $j \in \{1,\dots, d\}$.

Given $j\in\Delta$ and $n \in \Nb$, there exists $i$ and $k\in\{1,\ldots,d_i\}$, so that 
$$\sigma_j\left(\oplus_{i=1}^m\tau_{d_i}\left(\begin{bmatrix} 1 & n\\ 0 & 1\end{bmatrix}\right)\right)=
\sigma_k\left(\tau_{d_i}\left(\begin{bmatrix} 1 & n\\ 0 & 1\end{bmatrix}\right)\right).$$ 
 Then Lemma~\ref{singular value facts} implies the first claim in part (2) of Corollary~\ref{BCKM properties} , while the second claim follows from part (1) of
 Theorem \ref{thm: intro thm 1 body}.

\end{proof}

\begin{proof}[Proof of Theorem~\ref{thm: intro thm 1 body}]
Suppose that $\Gamma \subset \PSL(2,\Rb)$ is a geometrically finite group and $\rho:\Gamma\to\mathsf{SL}(d,\mathbb K)$ is $P_k$-Anosov.
By definition there exist $C>1$, $c>0$, and a $\rho$-equivariant family of norms $\norm{\cdot}_v$ on the fibers of $E_\rho \rightarrow \Usf(\Gamma)$ so that 
\begin{equation}
\label{lower bound}
\frac{\norm{Y}_{\phi_{t}(v)}}{\norm{Z}_{\phi_t(v)}} \leq Ce^{-c t} \frac{\norm{Y}_{v}}{\norm{Z}_{v}}
\end{equation}
for all $t > 0$, $v \in \Usf(\Gamma)$,  $Y\in\xi^k(v^+)$ and non-zero  $Z\in\xi^{d-k}(v^-)$. 

Fix a distance $d_\infty$ on $\partial \Hb^2 \cong \Sb^1$ which is induced by a Riemannian metric. 
The following very special case of a result of Abels-Margulis-Soifer \cite[Theorem 4.1]{AMS} plays a key role in the proof.

\begin{lemma}
\label{Fuchsian AMS}
There exist $\delta>0$ and a finite subset $\mathcal B$ of $\Gamma$ such that if
$\gamma\in\Gamma$, then there exists $\beta\in\mathcal B$ so that $\gamma\beta$ is hyperbolic and $d_\infty((\gamma\beta)^+,(\gamma\beta)^-)\ge\delta$.
\end{lemma}

Lemma \ref{Fuchsian AMS} allows us to reduce much of the proof to considering hyperbolic elements $\gamma \in \Gamma$ with $d_\infty(\gamma^+,\gamma^-)\ge\delta$. Since $\xi$ is transverse and continuous these elements have the following decomposition: there exists a compact set $\mathcal A\subset\mathsf{SL}(d,\mathbb K)$ so that if $\gamma \in \Gamma$ is hyperbolic and $d_\infty(\gamma^+,\gamma^-)\ge\delta$, then there exist $g_\gamma \in \Ac$, $A_\gamma \in \GL(k,\Kb)$ and $B_\gamma \in \GL(d-k,\Kb)$ with 
\begin{equation}\label{eqn:decomp_block_diag}
\rho(\gamma) = g_\gamma \begin{pmatrix} A_\gamma & 0 \\ 0 & B_\gamma \end{pmatrix} g_\gamma^{-1}, \quad  \xi^k(\gamma^+) = g_\gamma\big(\langle e_1,\ldots,e_k\rangle\big)\quad\mathrm{and}\quad \xi^{d-k}(\gamma^-) = g_\gamma\big(\langle e_{k+1},\ldots,e_d\rangle\big).
\end{equation}
There also exists a compact set $K$ of $\mathbb H^2$ so that any bi-infinite geodesic  whose  endpoints are a distance at least $\delta$ apart intersects $K$. Let 
$$R = \max\{ d_{\Hb^2}(i, z) : z \in K\}.$$ 

The next two lemmas establish (1) for any hyperbolic element $\gamma$ with $d_\infty(\gamma^+,\gamma^-)\ge \delta$. The proof of the  first lemma makes crucial use of the 
contraction properties of the flow, while the second lemma does not depend on the contraction properties.

\begin{lemma}
\label{spread case}
There exist $C_1, C_2>0$ such that
if $\gamma\in\Gamma$ is hyperbolic and $d_\infty(\gamma^+,\gamma^-)\ge \delta$, then 
$$\frac{\sigma_{k}(\rho(\gamma))}{\sigma_{k+1}(\rho(\gamma))}\geq C_1e^{c d_{\Hb^2}(i,\gamma(i))}$$
and 
$$\sigma_{k}(A_\gamma) \geq C_2e^{c d_{\Hb^2}(i,\gamma (i))}\sigma_1(B_\gamma). $$
\end{lemma}

\begin{proof} Recall that $r_v : \Rb \rightarrow \Hb^2$ denotes the geodesic  with $r_v^\prime(0)=v\in T^1 \Hb^2$. Since the geodesic joining $\gamma^+$ to $\gamma^-$ intersects $K$, there exists $v_0\in \Usf(\Gamma)$ with $r_{v_0}(0) \in K$ and $v_0^\pm = \gamma^\pm$. Also, notice that 
\begin{align*}
\ell(\gamma) = d_{\Hb^2}(r_{v_0}(0), \gamma (r_{v_0}(0)))  \geq d_{\Hb^2}(i,\gamma( i))-2R. 
\end{align*}

Since $\phi_{\ell(\gamma)}(v_0)=\gamma (v_0)$, the $\rho$-equivariance of the norms and Equation~\eqref{lower bound} imply that
$$\frac{\norm{Y}_{v_0}}{\norm{Z}_{v_0}}=\frac{\norm{\rho(\gamma)(Y)}_{\gamma (v_0)}}{\norm{\rho(\gamma)(Z)}_{\gamma (v_0)}} \leq Ce^{-c\ell(\gamma)} \frac{\norm{\rho(\gamma)(Y)}_{v_0}}{\norm{\rho(\gamma)(Z)}_{v_0}}$$
for all $Y\in\xi^k(\gamma^+)$ and non-zero $Z\in\xi^{d-k}(\gamma^-)$.

Since $K$ is compact, there exists $L$ so that if $v\in \Usf(\Gamma)$ and $r_v(0) \in K$, then $\norm{\cdot}_v$ is $L$-bilipschitz to the standard Euclidean norm $\norm{\cdot}_2$
on $\mathbb K^d$. Therefore, 
$$\frac{\norm{\rho(\gamma)(Y)}_{2}}{\norm{\rho(\gamma)(Z)}_{2}} \geq \frac{1}{CL^4}e^{c\ell(\gamma)} \frac{\norm{Y}_{2}}{\norm{Z}_{2}}$$
for all $Y\in\xi^k(\gamma^+)$ and non-zero $Z\in\xi^{d-k}(\gamma^-)$.  So by the max-min/min-max characterization of singular values 
\begin{align*}
\sigma_k(\rho(\gamma)) \geq \min_{\substack{Y \in \xi^{k}(\gamma^+) \\ Y \neq 0}} \frac{\norm{\rho(\gamma)(Y)}_{2}}{\norm{Y}_2} \geq  \frac{1}{CL^4}e^{c\ell(\gamma)} \max_{\substack{Z \in \xi^{d-k}(\gamma^-) \\ Z \neq 0}} \frac{\norm{\rho(\gamma)(Z)}_{2}}{\norm{Z}_2} \geq \frac{1}{CL^4}e^{c\ell(\gamma)} \sigma_{k+1}(\rho(\gamma)). 
\end{align*}
Hence $C_1 : = \frac{1}{CL^4} e^{-2Rc}$ suffices. 

Since $\Ac$ is compact, 
$$S=\max\left\{\frac{\sigma_1(g)}{\sigma_d(g)} : g\in\mathcal{A}\right\}$$
is finite. So, if $Y=(Y^\prime,0) \in \Kb^k \times \{0\}$ and $Z=(0,Z^\prime) \in \{0\} \times \Kb^{d-k}$, then 
\begin{align*}
\frac{\norm{A_\gamma (Y^\prime)}_{2}}{\norm{B_\gamma (Z^\prime)}_{2}} = \frac{\norm{g_\gamma^{-1} \rho(\gamma) g_\gamma (Y)}_{2}}{\norm{g_\gamma^{-1} \rho(\gamma) g_\gamma (Z)}_{2}} \geq \frac{1}{S}  \frac{\norm{ \rho(\gamma) g_\gamma (Y)}_{2}}{\norm{ \rho(\gamma) g_\gamma (Z)}_{2}} \geq \frac{1}{CL^4S}e^{c\ell(\gamma)}\frac{\norm{g_\gamma (Y)}_{2}}{\norm{ g_\gamma (Z)}_{2}} \geq \frac{1}{CL^4S^2}e^{c\ell(\gamma)}\frac{\norm{Y^\prime}_{2}}{\norm{Z^\prime}_{2}}
\end{align*}
So, again by the max-min/min-max characterization of singular values,
\begin{align*}
\sigma_k(A_\gamma) \geq \min_{Y^\prime \neq 0} \frac{\norm{A_\gamma (Y^\prime)}_{2}}{\norm{Y^\prime}_2} \geq  \frac{1}{CL^4S^2}e^{c\ell(\gamma)}  \max_{Z^\prime \neq 0} \frac{\norm{B_\gamma (Z^\prime)}_{2}}{\norm{Z^\prime}_2} \geq \frac{1}{CL^4S^2}e^{c\ell(\gamma)} \sigma_1(B_\gamma). 
\end{align*}
Hence $C_2 : = \frac{1}{CL^4S^2} e^{-2Rc}$ suffices. 
\end{proof}

\begin{lemma}
\label{spectral radius bound}
There exist $C_3, c_0>0$ so that if $\gamma\in\Gamma$ is hyperbolic and $d_\infty(\gamma^+,\gamma^-)\ge \delta$, then
$$\sigma_1(\rho(\gamma))\le C_3e^{c_0 d_{\Hb^2}(i,\gamma ( i))}.$$
In particular, 
$$\frac{\sigma_k(\rho(\gamma))}{\sigma_{k+1}(\rho(\gamma))} \le C_3^2e^{2c_0 d_{\Hb^2}(i,\gamma (i))}.$$
\end{lemma}

\begin{proof} Fix a $\rho$-equivariant family of norms $\norm{\cdot}_{v \in \Usf(\Gamma)}^*$ which descends to a canonical family of norms on $\wh{E}_\rho$. Then by Lemma~\ref{obs:distortion_of_norms_holder} there exist $C_0, c_0>0$ so that 
\begin{equation}
\label{upper bound 2}
\frac{1}{C_0}e^{-c_0 \abs{t} }\norm{Y}_v^*\le\norm{Y}_{\phi_t(v)}^*\le C_0e^{c_0 \abs{t} }\norm{Y}_v^*
\end{equation}
for all $v \in \Usf(\Gamma)$, $Y \in \Kb^d$ and $t \in \Rb$. 

Since the geodesic joining $\gamma^+$ to $\gamma^-$ intersects $K$, there exists $v_0\in \Usf(\Gamma)$ with $r_{v_0}(0) \in K$ and $v_0^\pm = \gamma^\pm$.  Also, recall that 
\begin{align*}
\ell(\gamma) = d_{\Hb^2}(r_{v_0}(0), \gamma (r_{v_0}(0)))  \leq d_{\Hb^2}(i,\gamma( i))+2R. 
\end{align*}
Since $\phi_{\ell(\gamma)}(v_0)=\gamma (v_0)$, the $\rho$-equivariance of the norms and Equation~\eqref{upper bound 2} imply that
$$\norm{Y}_{v_0}^*=\norm{\rho(\gamma)(Y)}_{\gamma (v_0)}^* \geq \frac{1}{C_0}e^{-c_0\ell(\gamma)}\norm{\rho(\gamma) (Y)}_{v_0}^*$$
if $Y\in\Kb^d$.

Since $K$ is compact, there exists $L$ so that if $v\in \Usf(\Gamma)$ and $r_v(0) \in K$, then $\norm{\cdot}_v^*$ is $L$-bilipschitz to the standard Euclidean norm $\norm{\cdot}_2$
on $\mathbb K^d$. Therefore, if $Y \in \Kb^d$, then 
$$\norm{\rho(\gamma)(Y)}_2 \leq C_3 e^{c_0 d_{\Hb^2}(i, \gamma(i))}\norm{Y}_2$$
where $C_3 := C_0 L^2 e^{2c_0R}$.  So 
\begin{equation*}
\sigma_1(\rho(\gamma))\le  C_3 e^{c_0 d_{\Hb^2}(i, \gamma(i))}.
\end{equation*}
Finally, notice that 
\begin{equation*}
\frac{\sigma_k(\rho(\gamma))}{\sigma_{k+1}(\rho(\gamma))} \leq \frac{\sigma_1(\rho(\gamma))}{\sigma_d(\rho(\gamma))} = \sigma_1(\rho(\gamma))\sigma_1(\rho(\gamma^{-1}))
 \le C_3^2e^{2c_0 d_{\Hb^2}(i,\gamma (i))}.\qedhere
\end{equation*}
\end{proof}

We can now prove part (1). Since $\mathcal B$ is finite, both  
$$S_{\Bc}=\max\left\{\frac{\sigma_1(\rho(\beta))}{\sigma_d(\rho(\beta))}\ \Big|\ \beta\in\mathcal B\right\} \quad \text{and} \quad b = \max \left\{ d_{\Hb^2}(i,\beta (i)) : \beta \in \mathcal B\right\}$$
are finite. 

Given any $\gamma\in\Gamma$, Lemma \ref{Fuchsian AMS} implies that there exists $\beta\in\mathcal B$
so that $d_\infty((\gamma\beta)^+,(\gamma\beta)^-)\ge \delta$. Then 
\begin{align*}
\frac{1}{S_{\Bc}} \frac{\sigma_{k}(\rho(\gamma\beta))}{\sigma_{k+1}(\rho(\gamma\beta))} \leq \frac{\sigma_{k}(\rho(\gamma))}{\sigma_{k+1}(\rho(\gamma))} \leq S_{\Bc} \frac{\sigma_{k}(\rho(\gamma\beta))}{\sigma_{k+1}(\rho(\gamma\beta))}
\end{align*}
and
$$\abs{d_{\Hb^2}(i,\gamma(i))-d_{\Hb^2}(i,\gamma\beta (i))}\le d_{\Hb^2}(i,\beta (i))\leq  b.$$
Therefore,
$$\frac{1}{A}\exp\left(\frac{1}{a}d_{\Hb^2}(i,\gamma (i))\right)\leq \frac{\sigma_k(\rho(\gamma))}{\sigma_{k+1}(\rho(\gamma))} \leq A\exp\left(a \ d_{\Hb^2}(i,\gamma (i))\right)$$
where $A =\max\{ \frac{1}{C_1}S_{\Bc}e^{cb}, C_3^2 S_{\Bc}e^{2c_0b}\}$ and $a = \max\{ 1/c,2c_0\}$. This proves (1).

Recall that if $T\in\mathsf{SL}(d,\mathbb K)$, then
\begin{align*}
\lambda_j(T^n) = \lim_{n \rightarrow \infty} \left( \sigma_j(T^n) \right)^{1/n}
\end{align*}
and $\ell(\gamma)=\lim\frac{d_{\Hb}(i,\gamma^n( i))}{n}$ for all $\gamma\in\Gamma$.
Therefore, part (2) follows immediately from part (1).
(One may give a direct proof of (2) in the spirit of Lemma \ref{spread case} by noting that there is a compact subset $\hat K$ of
$\mathbb H^2$ such that every hyperbolic element of $\Gamma$ is conjugate to an element of $\Gamma$ whose axis passes
through $\hat K$.)

Part (4) is a simple consequence of part (1) and Lemma~\ref{spectral radius bound}. Recall that $X_d(\mathbb R)=\mathsf{SL}(d,\mathbb R)/\mathsf{SO}(d)$ and $X_d(\mathbb C)=\mathsf{SL}(d,\mathbb C)/\mathsf{SU}(d)$. We may choose
$x_0$ to be either  $[\mathsf{SO}(d)]$  or $[\mathsf{SU}(d)]$. Then, after possibly scaling, we have the following formula for the distance on $X_d(\Kb)$, see for instance \cite[Cor. 10.42]{bridson-haefliger}, 
\begin{align*}
d_{X_d}\left( x_0,g(x_0) \right) = \sqrt{ \sum_{j=1}^d \abs{\log \sigma_j(g)}^2 }
\end{align*}
Therefore, applying part (1), we see that 
$$d_{X_d}(x_0,\rho(\gamma)(x_0))\ge\frac{1}{\sqrt{2}} \log\left(\frac{\sigma_k(\rho(\gamma))}{\sigma_{k+1}(\rho(\gamma))}\right)\ge \frac{a}{\sqrt{2}} d_{\Hb^2}(i,\gamma(i))-\frac{\log A}{\sqrt{2}}$$
for all $\gamma\in\Gamma$.
On the other hand, by Lemma~\ref{spectral radius bound}
\begin{align*}
d_{X_d}(x_0,\rho(\gamma)(x_0))& \le \sqrt{d} \max \left\{ \abs{\log \sigma_j(\rho(\gamma))}  \right\}= \sqrt{d} \max \left\{ \log \sigma_1(\rho(\gamma)),\log \sigma_1(\rho(\gamma)^{-1})\right\} \\
&\le \sqrt{d}\log C_3+\sqrt{d}c_0 d_{\Hb^2}(i,\gamma( i)),
\end{align*}
so the orbit map is a quasi-isometry, 
which completes the proof of part (4).

We now show that $\xi$ is strongly dynamics preserving. Notice that this immediately implies  that $\xi$ is unique.
Lemma \ref{lem: basic singular value} will then imply that if $\alpha\in\Gamma$ is parabolic, then $\rho(\alpha)$ is weakly unipotent,
while if $\gamma\in\Gamma$ is hyperbolic, then $\rho(\gamma)$ is $P_k$-proximal.

Fix a  sequence $\{\gamma_n\}$ with $\gamma_n \rightarrow x \in\Lambda(\Gamma)$ and $\gamma_n^{-1} \rightarrow y\in\Lambda(\Gamma)$. 
For each $n$ there exists $\beta_n \in \Bc$ such that  $\gamma_n \beta_n$ is hyperbolic and 
\begin{align*}
d_\infty\left( (\gamma_n \beta_n)^+, (\gamma_n \beta_n)^- \right) \geq \delta.
\end{align*}
Since the set $\Bc$ is finite, we can divide $\{\gamma_n\}$ into finitely many subsequences and only consider the case when $\beta_n = \beta$ for all $n$. 
Then $(\gamma_n \beta)^+ \rightarrow x$ and $ (\gamma_n \beta)^- \rightarrow \beta^{-1}(y)$. 

Let 
\begin{align*}
\rho(\gamma_n \beta) = g_n \begin{pmatrix} A_n & 0 \\ 0 & B_n \end{pmatrix} g_n^{-1}
\end{align*}
be the block diagonal decomposition from Equation~\eqref{eqn:decomp_block_diag}. Then by Lemma~\ref{spread case}
\begin{align*}
\lim_{n \rightarrow \infty} \frac{\sigma_1(B_n)}{\sigma_k(A_n)} =0. 
\end{align*}
so
$$g_n^{-1}\rho(\gamma_n\beta)g_n(W)\to \langle e_1,\dots,e_k\rangle$$ 
for all $W\in\mathrm{Gr}_k(\mathbb K^d)$ transverse to $\langle e_{k+1},\ldots, e_d\rangle$.
Further, by construction 
$$\xi^k((\gamma_n\beta)^+)=g_n\big(\langle e_1,\ldots,e_k\rangle\big)\qquad\mathrm{and}\qquad\xi^{d-k}((\gamma_n\beta)^-)=g_n\big(\langle e_{k+1}\ldots,e_d\rangle\big).$$
So by the continuity of $\xi$, 
\begin{align*}
\rho(\gamma_{n}\beta)(V) \rightarrow \xi^{k}(x)
\end{align*}
for all $V$ transverse to $\xi^{d-k}(\beta^{-1}(y))$. This implies that
\begin{align*}
\rho(\gamma_{n})( V) \rightarrow \xi^{k}(x)
\end{align*}
for all $V$ transverse to $\xi^{d-k}(y)=\rho(\beta)\big(\xi^{d-k}(\beta^{-1}(y))\big)$ which completes the proof of (3).
\end{proof}


\section{Basic properties of cusp representations} 

In this section we establish some useful properties of the cusp representations associated to type preserving representations introduced in Section~\ref{sec:canonical_norms_and_cusp_repn}. 

We say that a representation $\Psi:\mathsf{SL}(2,\mathbb R)\to \mathsf{SL}(d,\mathbb K)$ is {\em $P_k$-proximal} if 
$\Psi(\gamma)$ is $P_k$-proximal for some (any) hyperbolic element $\gamma\in\SL(2,\Rb)$. We first observe that a cusp representation associated to a $P_k$-Anosov representation is $P_k$-proximal. 

\begin{proposition}\label{prop cusp rep are proximal} Suppose that $\Gamma\subset\mathsf{PSL}(2,\mathbb R)$ is a geometrically finite group and 
$\rho:\Gamma\to\mathsf{SL}(d,\mathbb K)$ is $P_k$-Anosov. If $\alpha \in \Gamma$ is parabolic and 
$\Psi:\mathsf{SL}(2,\mathbb R)\to\mathsf{SL}(d,\mathbb K)$ is a cusp representation associated to $\alpha$ and $\rho(\alpha)$, then $\Psi$ is $P_k$-proximal. 
\end{proposition} 

\begin{proof}  Let $\wt\alpha\in\SL(2,\Rb)$ be the unique unipotent lift of $\alpha$. Since $\rho(\alpha)$ is weakly unipotent, the group $\overline{\ip{\rho(\alpha)_{ss}}}$ is compact. Then, since $\Psi(\wt\alpha^n)=\rho(\alpha)_{ss}^{-n}\rho(\alpha^n)$, there exists $L > 1$ such that 
$$
\frac{1}{L} \sigma_j( \rho(\alpha^n)) \leq \sigma_j(\Psi(\wt\alpha^n)) \leq L\sigma_j( \rho(\alpha^n))
$$
for all $j \in \{1,\dots, d\}$ and $n \in \Zb$. Then, since $\rho$ is strongly dynamics-preserving, Lemma \ref{lem: basic singular value} implies that 
$$
\lim_{n \rightarrow \infty} \frac{\sigma_k(\Psi(\wt\alpha^n))}{\sigma_{k+1}(\Psi(\wt\alpha^n))}\geq \frac{1}{L^2} \lim_{n \rightarrow \infty}\frac{\sigma_k(\rho(\alpha^n))}{\sigma_{k+1}(\rho(\alpha^n))}=\infty.
$$
Now, write $\wt\alpha^n=\ell_n\a2_{t_n}m_n$, where $\ell_n,m_n\in\SO(2)$ and $\a2_{t}=\begin{pmatrix}e^t&0\\0&e^{-t}\end{pmatrix}$. Since $\{\Psi(\a2_t):t\in\Rb\}$ is simultaneously diagonalizable, by increasing $L>1$ we may assume that
\[\frac{1}{L}\sigma_k(\Psi(\a2_t))\leq\lambda_k(\Psi(\a2_t))\leq L\sigma_k(\Psi(\a2_t))\]
for all $t\in\Rb$. Since $\Psi(\SO(2))$ is compact, we may increase $L >1$ further and assume that 
$$
\frac{1}{L} \leq \sigma_j(\Psi(g)) \leq L
$$
for all $g \in \SO(2)$ and $j\in\{1,\dots,d\}$. 

Then
\[
\lim_{n \rightarrow \infty}\frac{\lambda_k(\Psi(\a2_{t_n}))}{\lambda_{k+1}(\Psi(\a2_{t_n}))}\geq\frac{1}{L^2}\lim_{n \rightarrow \infty}\frac{\sigma_k(\Psi(\a2_{t_n}))}{\sigma_{k+1}(\Psi(\a2_{t_n}))}\geq\frac{1}{L^6}\lim_{n \rightarrow \infty}\frac{\sigma_k(\Psi(\wt\alpha^n))}{\sigma_{k+1}(\Psi(\wt\alpha^n))}=\infty,
\]
which implies that $\Psi(\a2_t)$ is $P_k$-proximal for all $t\in\Rb-\{0\}$. 

\end{proof}

Next we observe that a $P_k$-proximal representation is itself $P_k$-Anosov and admits a $P_k$-limit map, in the following sense.

\begin{proposition}
\label{cusp rep anosov}
If $\Psi:\mathsf{SL}(2,\mathbb R)\to\mathsf{SL}(d,\mathbb K)$ is a $P_k$-proximal representation, then 
there exists a continuous, $\Psi$-equivariant transverse map 
$$\eta=(\eta^k,\eta^{d-k}):\partial \Hb^2\to\mathrm{Gr}_k(\mathbb K^d)\times \mathrm{Gr}_{d-k}(\mathbb K^d)$$
with the following properties:

\begin{enumerate}
\item
If $\norm{\cdot}$ is a $\Psi$-invariant family of norms on the fibers of $T^1\mathbb H^2\times\Kb^d$, then
there exists $B,b>0$ so that if $t>0$, $v\in T^1\mathbb H^2$, $Y\in\eta^k(v^+)$ and $Z\in\eta^{d-k}(v^-)$ is non-zero, then 
\begin{align}
\frac{\norm{Y}_{\phi_t(v)}}{\norm{Z}_{\phi_t(v)}}\le Be^{-b t} \frac{\norm{Y}_v}{\norm{Z}_v}.
\end{align}
\item
If $\{\gamma_n\}$ is a sequence in $\mathsf{SL}(2,\mathbb R)$ with $\gamma_n \rightarrow x\in\partial\Hb^2$ and $\gamma_n^{-1} \rightarrow y\in\partial\Hb^2$, 
then 
\begin{align*}
\Psi(\gamma_n)(V) \rightarrow \eta^{k}(x)
\end{align*}
locally uniformly for all $V \in \Gr_k(\Kb^d)$ transverse to $\eta^{d-k}(y)$.
\item If $g \in \SL(d,\Kb)$ commutes with the elements of $\Psi(\SL(2,\Rb))$, then $g \circ \eta = \eta$.
\end{enumerate}
\end{proposition}

\begin{proof}
By conjugating $\Psi$ we can assume that 
\begin{align}\label{eqn: decomp}
\Psi = \oplus_{i=1}^m \tau_{d_i}.
\end{align}
Let $\xi_{d_i} :\partial \Hb^2 \rightarrow \Fc(\Kb^{d_i})$ denote the $\tau_{d_i}$-equivariant boundary map described in Section~\ref{sec:repn_thy_of_SL2}.

By definition, if $\gamma\in\SL(2,\Rb)$ is a hyperbolic element, then $\tau_{d_i}(\gamma)$ is diagonalizable with eigenvalues having pairwise distinct absolute values. Furthermore, for all $k=1,\dots,d_i-1$, $\xi_{d_i}^{k}(\gamma^+)$ is the direct sum of the eigenspaces of the $k$ largest eigenvalues of $\tau_{d_i}(\gamma)$.

First, we construct the map $\eta$. Observe that since $\Psi$ is $P_k$-proximal, $\Psi(\gamma)$ is $P_k$-proximal and $P_{d-k}$-proximal for every hyperbolic element $\gamma \in \SL(2,\Rb)$. Thus, for all $i=1,\dots,m$, there are integers $k_i\in\{0,\dots,d_i\}$ such that 
\begin{itemize}
\item $\sum_{i=1}^mk_i=k$, and for all hyperbolic $\gamma\in\SL(2,\Rb)$, $\oplus_{i=1}^m\xi_{d_i}^{k_i}(\gamma^+)\in\Gr_k(\Kb^d)$ is the attracting fixed point for the action of $\Psi(\gamma)$
on $\Gr_k(\Kb^d)$, and
\item $\sum_{i=1}^md-k_i=d-k$, and for all hyperbolic $\gamma\in\SL(2,\Rb)$, $\oplus_{i=1}^m\xi_{d_i}^{d_i-k_i}(\gamma^-)\in\Gr_{d-k}(\Kb^d)$ is the attracting fixed point for the action of 
$\Psi(\gamma^{-1})$ on $\Gr_{d-k}(\Kb^d)$.
\end{itemize}
Then set 
\[\eta^{k}:=\oplus_{i=1}^m\xi_{d_i}^{k_i}:\partial\Hb^2\to\Gr_k(\Kb^d)\] 
and 
\[\eta^{d-k}:=\oplus_{i=1}^m\xi_{d_i}^{d_i-k_i}:\partial\Hb^2\to\Gr_{d-k}(\Kb^d).\]
Note that this pair of maps are continuous, $\Psi$-equivariant, and transverse.

For every $v \in T^1 \Hb^2$, there exists a one-parameter subgroup $\{\a2_t\}_{t\in\mathbb R}$ of hyperbolic elements in $\mathsf{SL}(2,\mathbb R)$
so that $\a2_t(v)=\phi_t(v)$ for all $t \in \Rb$. Since $\Psi(\a2_1)$ is $P_k$-proximal, 
$$b:=\log\left(\frac{\lambda_k(\Psi(\a2_1))}{\lambda_{k+1}(\Psi(\a2_1))}\right)>0.$$
so
$$\frac{\lambda_k(\Psi(\a2_t))}{\lambda_{k+1}(\Psi(\a2_t))}=e^{bt}$$
for all $t>0$.

If $X\in \Kb^d$, then
$$\norm{X}_{\phi_t(v)}=\norm{\Psi(\a2_t)^{-1}(X)}_v$$
so, if $Y\in\eta^k(v^+)$ and $Z\in\eta^{d-k}(v^-)$ are non-zero, then
$$\frac{\norm{Y}_{\phi_t(v)}}{\norm{Z}_{\phi_t(v)}}\le\frac{\lambda_{k+1}(\Psi(\a2_t))}{\lambda_k(\Psi(\a2_t))} \frac{\norm{Y}_v}{\norm{Z}_v}
=e^{-bt}  \frac{\norm{Y}_v}{\norm{Z}_v}$$
which completes the proof of (1).

If $\{\gamma_n\}$ is a sequence in $\mathsf{SL}(2,\mathbb R)$ with $\gamma_n \rightarrow x$ and $\gamma_n^{-1}\rightarrow y$,
one can write $\gamma_n= \ell_n\a2_{t_n} m_n$ where $\ell_n, m_n \in \mathsf{SO}(2)$, $t_n \ge 0$, and $\a2_t=\begin{pmatrix}e^t &0\\0&e^{-t}\\ \end{pmatrix}$.
By assumption, $t_n\to\infty$,  $\ell_n\to \ell$ and $m_n\to m$ with $\ell(\infty)=x$ and $m^{-1}(0)=y$.

Notice that if $W\in\mathrm{Gr}_k(\mathbb R^d)$ is transverse to  $\eta^{d-k}(0)$, then $\{\Psi(\a2_{t_n})(W)\}$ converges to
$\eta^k(\infty)$ (since, by definition, $\eta^k(\infty)$ is the attracting $k$-plane of $\Psi(\a2_1)$ and $\eta^{d-k}(0)$ is the repelling $(d-k)$-plane).
So, by equivariance,  if $V$ is transverse to $\eta^{d-k}(y)$, then $\Psi(m_n)(V)$ is transverse to 
$\eta^{d-k}(0)=\lim \Psi(m_n)(\eta^{d-k}(y))$ for all large enough $n$. Thus, $\{\Psi(\a2_{t_n}m_n)(V)\}$ converges to $\eta^k(\infty)$ locally uniformly in $V$.
Therefore, $\{\Psi(\gamma_n)(V)\}=\{\Psi(\ell_n\a2_{t_n} m_n)(V)\}$ converges to $\Psi(\ell)(\eta^k(\infty))=\eta^k(x)$ locally uniformly in $V$. This proves (2).

To prove (3), fix $x \in \partial \Hb^2$ and a hyperbolic element $\gamma \in \SL(2,\Rb)$ so that $\gamma^+ = x$. 
Then $\Psi(\gamma)$ is $P_k$-proximal and $\eta^k(x) \in \Gr_k(\Kb^d)$ is the attracting fixed point of $\Psi(\gamma)$,  so $g \eta^k(x) = \eta^k(x)$ since $g$ commutes with $\Psi(\gamma)$. 
Similar reasoning shows that $g\eta^{d-k}(x) = \eta^{d-k}(x)$. 
\end{proof}

The following technical result says that the image of the limit map of an Anosov representation is asymptotically homogeneous at a parabolic fixed point. 
In one of the arguments that follow, we need uniform control over continuous families of Anosov representations, so we introduce a parameter $u$.

\begin{proposition}\label{thm:asymptotics_of_limit_curve} 
Suppose that $\alpha\in\SL(2,\Rb)$ is a parabolic element, $X$ is a closed $\alpha$-invariant subset of $\partial\mathbb H^2$ containing the fixed point $\alpha^+=\alpha^-$ of $\alpha$, 
$\Psi:\mathsf{SL}(2,\mathbb R)\to\mathsf{SL}(d,\mathbb K)$ is a $P_k$-proximal representation with 
$P_k$-limit map $\eta$, and $\ell \in \SL(d,\Kb)$ is elliptic and commutes with the elements of $\Psi(\SL(2,\Rb))$. 

Let $U$ be a compact metric space and  $\{ g_{u} \}_{u \in U}$ be a continuous family of elements in $\SL(d,\Rb)$. Suppose 
$$\xi=(\xi^k,\xi^{d-k}):U \times X\to \mathrm{Gr}_k(\mathbb K^d)\times\mathrm{Gr}_{d-k}(\mathbb K^d)$$
is continuous and for each $u \in U$, the map $\xi_u := \xi(u,\cdot)$ is transverse, $\xi_u(\alpha^+)=g_u\eta(\alpha^+)$ and 
\begin{equation}
\label{commuting property}
\xi_u \circ \alpha=g_u\ell\Psi(\alpha)g_u^{-1} \circ \xi_u.
\end{equation}

If $\gamma\in \mathsf{SL}(2,\mathbb R)$ is a hyperbolic element with attracting fixed point $\alpha^+$,  then
\begin{align*}
\lim_{n \rightarrow \infty}\left( g_u\Psi(\gamma)^{-n}g_u^{-1}\circ \xi_u \circ \gamma^n\right) (x_n)  =  g_u\eta(x)
\end{align*}
if $\lim x_n = x\in \partial\Hb^2$ and $x_n \in \gamma^{-n}(X)$ for all $n$. Moreover, the convergence is uniform in $u \in U$. 
\end{proposition}

\begin{proof} Fix a distance $d_G$ on $\mathrm{Gr}_k(\mathbb K^d)\times\mathrm{Gr}_{d-k}(\mathbb K^d)$ induced by a Riemannian metric. 

Suppose the proposition fails for a sequence $\{x_n\}\subset X$ with $\lim x_n = x\in X$. Then there exist $\epsilon > 0$, $\{n_j\}$ converging to infinity and a sequence $\{u_j\}$ in $U$ 
such that 
\begin{align*}
d_{G} \left( \left( g_{u_j}\Psi(\gamma)^{-n_j}g_{u_j}^{-1}\circ \xi_{u_j} \circ \gamma^{n_j}\right) (x_{n_j}), g_{u_j} \eta(x) \right) > \epsilon
\end{align*}
for all $j$. Passing to subsequences we may suppose that $u_j \rightarrow u_\infty \in U$.

For notational convenience, let 
$$\eta_j= g_{u_j}\Psi(\gamma)^{-n_j}g_{u_j}^{-1}\circ \xi_{u_j} \circ \gamma^{n_j}$$
for all $j \in \Nb$. Since $\eta$ is $\Psi$-equivariant and $\xi_u(\alpha^+) = g_u \eta(\alpha^+)$ for all $u\in U$,
\begin{align*}
g_{u_j}\eta(\alpha^+) &= g_{u_j}\Psi(\gamma)^{-n_j} \eta(\alpha^+)=g_{u_j}\Psi(\gamma)^{-n_j}g_{u_j}^{-1}\circ \xi_{u_j} (\alpha^+) \\
& =g_{u_j}\Psi(\gamma)^{-n_j}g_{u_j}^{-1}\circ \xi_{u_j} \circ \gamma^{n_j}(\alpha^+)= \eta_j(\alpha^+).
 \end{align*}
So, by passing to a tail of our sequences we may assume that $x_{n_j}\ne\alpha^+$ for every $j$.

First suppose that $\{y_j=\gamma^{n_j}(x_{n_j})\}$ lies in a compact subset of $\partial\mathbb H^2-\{\alpha^+\}$. Then $x = \gamma^-$. Since $\xi_{u_j} \rightarrow \xi_{u_\infty}$, 
there exists $N > 0$ sufficiently large so that $\{ g_{u_j}^{-1}\xi_{u_j}(y_j)\}_{n \geq N}$ lies in a compact subset of flags transverse to $g_{u_\infty}^{-1}\xi_{u_\infty}(\alpha^+)=\eta(\gamma^+)$. 
Hence Proposition \ref{cusp rep anosov} part (2) implies that 
$$
g_{u_j}^{-1} \eta_j(x_{n_j}) = \Psi(\gamma)^{-n_j}(g_{u_j}^{-1}\xi_{u_j}(y_j))
$$  
converges to $\eta(\gamma^-)$. Thus
\begin{align*}
\epsilon \leq \liminf_{j \rightarrow \infty} d_G \left(  \eta_j (x_{n_j}), g_{u_j} \eta(x) \right) = d_G \left(g_{u_\infty} \eta(\gamma^-), g_{u_\infty} \eta(\gamma^-) \right) =0
\end{align*}
and we have a contradiction. 

Now suppose that $\{y_j\}$ does not lie in a compact subset of $\partial\mathbb H^2-\{\alpha^+\}$. 
We may assume without loss of generality that $\gamma^-=0$, $\gamma^+=\infty=\alpha^+$, and $\alpha  =\u2_1$ where 
\begin{align*}
\u2_t = \begin{pmatrix} 1 & t \\ 0 & 1 \end{pmatrix} \in \SL(2,\Rb).
\end{align*}
Then $y_j \in \Rb \subset \partial \Hb^2$ for all $j$. Let $z_j=\lfloor y_j\rfloor\in\mathbb Z$ and $w_j=y_j-z_j\in[0,1]$ and set
 $\delta_j=\gamma^{-n_j}\u2_{z_j}=\gamma^{-n_j}\alpha^{z_j}$.
Notice that $\alpha^{z_j}(w_j)=y_j$, 
$\delta_j(w_j)=x_{n_j}$, $\delta_j^-=\infty$ and $\delta_j^+\to x$. Passing to a subsequence, we can suppose that $w_j \rightarrow w \in [0,1]$. 
Proposition \ref{cusp rep anosov} part (2) then implies that
$\Psi(\delta_j)(V)$ converges to $\eta^k(x)$  locally uniformly for all $V\in\mathrm{Gr}_k(\mathbb K^d)$
which are transverse to $\eta^{d-k}(\infty)$. Also,
\begin{align*}
\eta_j^k(x_{n_j})&=g_{u_j}\Psi(\gamma^{-n_j})g_{u_j}^{-1}\big(\xi_{u_j}^k(\alpha^{z_j}(w_j)\big)=g_{u_j}\Psi(\gamma^{-n_j})g_{u_j}^{-1}\big(g_{u_j}\ell^{z_j}\Psi(\alpha^{z_j})g_{u_j}^{-1})\xi^k_{u_j}(w_j)\big)\\
 &=  g_{u_j}\ell^{z_j}\Psi(\gamma^{-n_j}\alpha^{z_j}) g_{u_j}^{-1}\big(\xi_{u_j}^k(w_j)\big)
=g_{u_j}\ell^{z_j}\Psi(\delta_j)\big(g_{u_j}^{-1}\xi_{u_j}^k(w_j)\big)
\end{align*}
(where in the first line we apply assumption \eqref{commuting property}).
We may pass to a subsequence so that $\ell^{z_j} \rightarrow \ell_\infty \in \SL(d,\Kb)$. Then $\ell_\infty$ is also elliptic and commutes with the elements of $\Psi(\SL(2,\Rb))$
and hence, by Proposition \ref{cusp rep anosov} part (3), fixes each element in the image of $\eta$. Then 
$$
\lim_{j \rightarrow \infty} \eta_j^k(x_{n_j}) = g_{u_\infty} \ell_\infty\eta^k(x) = g_{u_\infty} \eta^k(x)
$$
since $g_{u_j}^{-1}\xi_{u_j}^k(w_j) \rightarrow g_{u_\infty}^{-1}\xi_{u_\infty}^k(w)$, $g_{u_\infty}^{-1}\xi_{u_\infty}^k(w)$ is transverse to $\eta^{d-k}(\infty)=g_{u_\infty}^{-1}\xi_{u_\infty}^k(\infty)$,
and $\ell_\infty \circ \eta=\eta$.

Reversing the roles of $k$ and $d-k$, we may similarly show that $\eta_j^{d-k}(x_{n_j})\to g_{u_\infty}\eta^{d-k}(x)$. Hence we again have a contradiction. 
\end{proof}

\section{A dynamical characterization of linear Anosov representations} 

In this section, we prove Theorem \ref{thm:equivalent_to_Anosov}. The forward implication has already been established
as part (3) of Theorem~\ref{thm: intro thm 1}. The reverse implication follows from the following more general statement. 

First, recall from Section~\ref{sec: Anosov representations} that a transverse, $\rho$-equivariant, continuous map $\xi=(\xi^k,\xi^{d-k}):\Lambda(\Gamma)\to\mathrm{Gr}_k(\mathbb K^d)\times\mathrm{Gr}_{d-k}(\mathbb K^d)$ induces a continuous decomposition of $\wh E_\rho$ into a pair of sub-bundles $\wh E_\rho=\wh\Theta^{k}\oplus\wh\Xi^{d-k}$ of rank $k$ and $d-k$ respectively. Further, the flow $\phi_t$ induces a flow on $\Hom(\wh\Xi^{d-k},\wh\Theta^k)$, which we denote by $\phi_t$, and any canonical norm on $\wh E_\rho$ induces a canonical norm on $\Hom(\wh\Xi^{d-k},\wh\Theta^k)$ which is simply given by the associated operator norm.

\begin{theorem}\label{dynamics implies anosov}
Suppose $\Gamma \subset \PSL(2,\Rb)$ is a geometrically finite group and $\rho : \Gamma \rightarrow \SL(d,\Kb)$ is a representation. 
If there exists a $\rho$-equivariant, transverse, continuous, strongly dynamics preserving map 
$\xi=(\xi^k,\xi^{d-k}):\Lambda(\Gamma)\to\mathrm{Gr}_k(\mathbb K^d)\times\mathrm{Gr}_{d-k}(\mathbb K^d)$,
then $\rho$ is type-preserving and the flow $\phi_t$ on $\mathrm{Hom}(\wh\Xi^{d-k},\wh\Theta^k)$ is uniformly contracting with respect to
any canonical norm on $\mathrm{Hom}(\wh\Xi^{d-k},\wh\Theta^k)$.
In particular,  $\rho$ is $P_k$-Anosov and $\xi$ is its $P_k$-Anosov limit map.
\end{theorem}

The following is an immediate corollary of Theorem \ref{dynamics implies anosov} and Theorem \ref{thm: intro thm 1 body}.

\begin{corollary}\label{cor: contraction canonical norm} 
Suppose $\Gamma \subset \PSL(2,\Rb)$ is a geometrically finite group and $\rho : \Gamma \rightarrow \SL(d,\Kb)$ is a $P_k$-Anosov representation. Then the flow $\phi_t$ on $\Hom(\wh\Xi^{d-k},\wh\Theta^k)$ is uniformly contracting with respect to any canonical norm on $\mathrm{Hom}(\wh\Xi^{d-k},\wh\Theta^k)$.
\end{corollary}

As another corollary, we see that if $\rho$ is Zariski dense, then $\rho$ is $P_k$-Anosov if it admits  a transverse limit map, which generalizes
a result of Guichard and Wienhard from the uncusped Anosov setting \cite[Theorem 4.11]{guichard-wienhard}.

\begin{corollary}\label{zariski dense case}
Suppose $\Gamma \subset \PSL(2,\Rb)$ is a geometrically finite group and $\rho : \Gamma \rightarrow \SL(d,\Kb)$ is a representation. 
\begin{enumerate}
\item
If $\rho$ is irreducible, and there exists a $\rho$-equivariant, transverse, continuous, map 
$\xi=(\xi^1,\xi^{d-1}):\Lambda(\Gamma)\to\mathrm{Gr}_1(\mathbb K^d)\times\mathrm{Gr}_{d-1}(\mathbb K^d)$,
then $\rho$ is $P_1$-Anosov and $\xi$ is its $P_1$-Anosov limit map.
\item
If  $\wedge^k \rho : \Gamma \rightarrow \SL(\wedge^k \Kb^d)$ is irreducible (e.g. if $\rho(\Gamma)$ is Zariski dense in $\mathsf{SL}(d,\mathbb K)$) and there exists a $\rho$-equivariant, transverse, continuous, map 
$\xi=(\xi^k,\xi^{d-k}):\Lambda(\Gamma)\to\mathrm{Gr}_k(\mathbb K^d)\times\mathrm{Gr}_{d-k}(\mathbb K^d)$,
then $\rho$ is $P_k$-Anosov and $\xi$ is its $P_k$-Anosov limit map.
\end{enumerate}
\end{corollary}

\begin{proof} (1): It is enough to show that $\xi$ is strongly dynamics preserving. Fix an escaping sequence $\{\gamma_n\}$ in $\Gamma$ with $\gamma_n \rightarrow x$ and $\gamma_n^{-1} \rightarrow y$. Let $[\rho(\gamma_n)]$ denote the image of $\rho(\gamma_n)$ in $\Pb( \End(\Kb^{d}))$. Then it is enough to show that $[\rho(\gamma_n)]$ converges to the element $T \in \Pb(\End(\Kb^d))$ with $\ker(T) = \xi^{d-1}(y)$ and ${\rm Image}(T) = \xi^1(x)$. Since $\Pb( \End(\Rb^{d}))$ is compact it is enough to show that every convergent subsequence of $\rho(\gamma_n)$ converges to $T$. So suppose that $[\rho(\gamma_n)] \rightarrow S$ in $\Pb( \End(\Kb^{d}))$. Then
\begin{align*}
S(v) = \lim_{n \rightarrow \infty} \rho(\gamma_n)(v)
\end{align*}
for all $v \in \Pb(\Kb^d) \setminus \Pb(\ker S)$. 

We first claim that ${\rm Image}(S) = \xi^1(x)$. Since $\rho: \Gamma \rightarrow \SL(d,\Kb)$ is irreducible, there exists $x_1, \dots, x_{d} \in \Lambda(\Gamma)$ 
so that  $\xi(x_1), \dots, \xi(x_{d})$ spans $\Kb^{d}$. Since $\partial \Gamma$ is  perfect, we can perturb each $x_j$ and assume that
\begin{align*}
y \notin  \{ x_1, \dots, x_d\}.
\end{align*}
Then $\rho(\gamma_n)\big(\xi^1(x_j)\big) \rightarrow \xi^1(x)$. Since $\{\xi^1(x_1), \dots, \xi^1(x_{d})\}$ spans $\Kb^{d}$, we can relabel and suppose that 
\begin{align*}
\ker S \oplus \xi^1(x_1) \oplus \dots \oplus \xi^{1}(x_m) = \Kb^d
\end{align*}
where $m = d - \dim \ker S$. Then 
\begin{align*}
S\left(\xi^1(x_j)\right) = \lim_{n \rightarrow \infty} \rho(\gamma_n)\big(\xi^1(x_j)\big) = \xi^1(x)
\end{align*}
for all $1 \leq j \leq m$. Hence ${\rm Image}(S) = \xi^1(x)$. 

To compute the kernel, we notice that $\Gr_{d-1}(\Kb^d)$ may be identified with $\mathbb P(\Kb^{d*})$ by identifying a hyperplane $Q$ in $\Kb^d$ with the projective
class of linear functionals with kernel $Q$.  Notice that $[^t\rho(\gamma_n)]$ converges to $^tS$ in $\Pb(\End(\Kb^{d*}))$. 
Repeating the argument above shows that $\mathrm{Image}(^tS)=\xi^{d-1}(y)$, so  the kernel of $S$ is $\xi^{d-1}(y)$. 

(2):  One can argue similarly using the Pl\"ucker embeddings. 
\end{proof}

\begin{proof}[Proof of Theorem \ref{dynamics implies anosov}.]
Suppose  that $\Gamma \subset \PSL(2,\Rb)$ is a geometrically finite group, $\rho:\Gamma\to \mathsf{SL}(d,\mathbb K)$ is a representation and that 
$$\xi=(\xi^k,\xi^{d-k}):\Lambda(\Gamma)\to\mathrm{Gr}_k(\mathbb K^d)\times\mathrm{Gr}_{d-k}(\mathbb K^d)$$
is a continuous, transverse, $\rho$-equivariant, strongly dynamics-preserving map. Lemma \ref{lem: basic singular value} implies
that $\rho$ is type-preserving.

Let $\norm{\cdot}$ be a canonical family of norms on $\wh E_\rho$  and let $\mathcal C$ be a full collection of
embedded cusp neighborhoods so that the restriction to the fibers over
$\wh{\Usf}(\Gamma)_C$ is canonical for all $C\in\mathcal C$. We will also use $\norm{\cdot}$ to denote the lift of $\norm{\cdot}$ to a continuous family of norms on the fibers of $\Usf(\Gamma) \times \Kb^d$.

The proof divides into two parts. We first use properties of the canonical family of norms to control the flow over the thin part
$\wh{\Usf}(\Gamma)_{\mathcal C}$. We then use a compactness argument to control the flow on the complement.

\begin{proposition}
\label{decay on cusps}
If $C\in\mathcal C$, then there exist constants $b_C$ and $B_C$ and an embedded cusp sub-neighborhood $C^\prime\subset C$ such that if 
$v\in\wh\Usf(\Gamma)_{C^\prime}$, $t\ge 0$ and
$\phi_s(v)\in \wh\Usf(\Gamma)_{C^\prime}$ for all $s\in [0,t]$,
then 
\begin{align}\label{eqn: decay on cusps}
\frac{\norm{\phi_t(Y)}_{\phi_t(v)}}{\norm{\phi_t(Z)}_{\phi_t(v)}}\le B_C e^{-b_C t} \frac{\norm{Y}_v}{\norm{Z}_v}
\end{align}
for all $Y\in\wh\Theta^k|_v$ and non-zero $Z\in\wh\Xi^{d-k}|_v$.
\end{proposition}

\begin{proof}

Suppose that $C=\langle\alpha\rangle\backslash H$. Then it suffices to find a horodisc $H^\prime \subset H$ and constants $b_C,B_C$ such that: if 
$v\in\Usf(\Gamma)_{H'}$, $t\ge 0$ and
$\phi_s(v)\in \Usf(\Gamma)_{H'}$ for all $s\in [0,t]$,
then 
\begin{align}
\frac{\norm{Y}_{\phi_t(v)}}{\norm{Z}_{\phi_t(v)}}\le B_C e^{-b_C t} \frac{\norm{Y}_v}{\norm{Z}_v}
\end{align}
for all $Y\in\xi^k(v^+)$ and non-zero $Z\in\xi^{d-k}(v^-)$.

After possibly replacing $C$ with a subcusp, there exists a cusp representation $\Psi:\mathsf{SL}(2,\mathbb  R)\to \mathsf{SL}(d,\mathbb K)$ for $\alpha$ and $\rho(\alpha)$ 
such that $\norm{\cdot}$ on $\Usf(\Gamma)_H$ coincides with a $\rho(\alpha)_{ss}$-invariant, $\Psi$-equivariant family of norms $\norm{\cdot}^\star_{v \in T^1 \Hb^2}$. 

Let $\eta$ be the $P_k$-limit map of $\Psi$. Proposition \ref{cusp rep anosov} implies that there exists $B,b>0$ such that
\begin{equation}
\label{eta contraction}
\frac{\norm{Y}_{\phi_t(v)}^\star}{\norm{Z}_{\phi_t(v)}^\star}\le Be^{-b t} \frac{\norm{Y}_v^\star}{\norm{Z}_v^\star}
\end{equation}
for all  $t \geq 0$, $v \in T^1 \Hb^2$, $Y\in\eta^k(v^+)$ and non-zero $Z\in\eta^{d-k}(v^-)$. Choose $b_C=\frac{b}{2}$ and $T > 0$ so that 
\begin{align}
\label{eqn:choice_of_beta_T}
Be^{-bt}<e^{-b_Ct} \text{ for all\ } t > T.
\end{align}

We claim that there is a horodisk $H^\prime\subset H$ so that if $v\in \Usf(\Gamma)_{H^\prime}$, $t\in[T,2T]$ and $\phi_s(v)\in\Usf(\Gamma)_{H^\prime}$ for all $s\in[0,t]$, then 
\begin{align}\label{eqn: claim}
\frac{\norm{Y}_{\phi_t(v)}}{\norm{Z}_{\phi_t(v)}}\le e^{-b_C t} \frac{\norm{Y}_v}{\norm{Z}_v}
\end{align}
for all $Y\in\xi^{k}(v^+)$ and non-zero $Z\in\xi^{d-k}(v^-)$. If this is not the case, then there exists
\begin{itemize}
\item a sequence $\{t_n\}$ in $[T,2T]$,
\item a nested sequence $\{H_n\}$ of horodisks centered at the fixed point of $\alpha$ whose intersection is empty,
\item a sequence $\{v_n\}$ such that $\phi_s(v_n)\in \Usf(\Gamma)_{H_n}$ for all $s\in[0,t_n]$,
\item a sequence vectors $\{Y_n\}$ such that $Y_n\in\xi^k(v_n^+)$,
\item a sequence of non-zero vectors $\{Z_n\}$ such that $Z_n\in\xi^{d-k}(v_n^-)$,
\end{itemize}
such that
\begin{align}\label{eq:main_contradiction}
\frac{\norm{Y_n}_{\phi_{t_n}(v_n)}}{\norm{Z_n}_{\phi_{t_n}(v_n)}}> e^{-b_Ct_n} \frac{\norm{Y_n}_{v_n}}{\norm{Z_n}_{v_n}}.
\end{align}

Let 
$$\u2_t= \begin{pmatrix} 1 & t \\ 0 & 1\end{pmatrix} \quad \text{and} \quad \a2_s = \begin{pmatrix} e^s & 0 \\ 0 & e^{-s} \end{pmatrix}.$$
As usual, by conjugating, we can assume that $\Psi=\oplus_{i=1}^m\tau_{d_i}$, and $\alpha = [\u2_1]$. 

There is a sequence $\{s_n\}\to\infty$  of positive real numbers and a sequence  $\{m_n\}$ of integers
such that $\{w_n:=\a2_{-s_n}\u2_{m_n}(v_n)\}$ is relatively compact in $T^1 \Hb^2$. By passing to a subsequence, we may assume that $w_n$ converges to some $w_\infty\in T^1 \Hb^2$, and that 
\begin{align*}
(V_n,W_n):=\left(\frac{\Psi(\a2_{-s_n}\u2_{m_n})(Y_n)}{\norm{\Psi(\a2_{-s_n}\u2_{m_n})(Y_n)}_{w_n}^{\star}} ,\frac{\Psi(\a2_{-s_n}\u2_{m_n})(Z_n)}{\norm{\Psi(\a2_{-s_n}\u2_{m_n})(Z_n)}_{w_n}^{\star}}\right)
\end{align*}
converges to some $(V_\infty,W_\infty)\in\Kb^d\times\Kb^d$. By definition,
$$(V_n,W_n) \in\Psi(\a2_{-s_n}\u2_{m_n})\left(\xi^{k}(v_n^+)\times\xi^{d-k}(v_n^-)\right)=\xi_n^{k}(w_n^+)\times\xi_n^{d-k}(w_n^-)$$
where $\xi_n=\Psi(\a2_{-s_n})\circ \xi\circ \a2_{s_n}$. Proposition \ref{thm:asymptotics_of_limit_curve} (applied when $U$ is a singleton) implies that $\lim_{n \rightarrow \infty} \xi_n=\eta$,
so
$$(V_\infty,W_\infty)\in \eta^k(w_\infty^+)\times\eta^{d-k}(w_\infty^-).$$

Since $\norm{\cdot}_{v}^{\star}$ is $\Psi$-equivariant, Equation~\eqref{eq:main_contradiction} implies that
\begin{align}
\label{eq:main_contradiction_2}
\frac{\norm{V_n}_{\phi_{t_n}(w_n)}^\star}{\norm{W_n}_{\phi_{t_n}(w_n)}^\star}&=C_n \frac{\norm{Y_n}_{\phi_{t_n}(v_n)}}{\norm{Z_n}_{\phi_{t_n}(v_n)}}> C_ne^{-b_C t_n} \frac{\norm{Y_n}_{v_n}}{\norm{Z_n}_{v_n}}=e^{-b_C t_n} \frac{\norm{V_n}_{w_n}^\star}{\norm{W_n}_{w_n}^\star}.
\end{align}
where $C_n=\frac{\norm{\Psi(\a2_{-s_n}\u2_{m_n})(Z_n)}_{w_n}^\star} {\norm{\Psi(\a2_{-s_n}\u2_{m_n})(Y_n)}_{w_n}^\star}$.
By passing to a final subsequence, we can suppose that $t_n \rightarrow t \in [T,2T]$. Then taking $n \rightarrow \infty$ in Equation~\eqref{eq:main_contradiction_2} we obtain 
\begin{align*}
\frac{\norm{V_\infty}_{\phi_{t}(w_\infty)}^\star}{\norm{W_\infty}_{\phi_{t}(w_\infty)}^\star} \geq e^{-b_Ct} \frac{\norm{V_\infty}_{w_\infty}^\star}{\norm{W_\infty}_{w_\infty}^\star}> Be^{-bt} \frac{\norm{V_\infty}_{w_\infty}^\star}{\norm{W_\infty}_{w_\infty}^\star}.
\end{align*}
Since $V_\infty\in\eta^k(w_\infty^+)$ and $W_{\infty}\in\eta^{d-k}(w_\infty^-)$, this contradicts \eqref{eta contraction}, thus proving the claim.

Repeated applications of the claim implies that if $v\in\Usf(\Gamma)_{H^\prime}$ and $t\ge T$ is a number such that 
$\phi_s(v)\in \Usf(\Gamma)_{H^\prime}$ for all $s\in [0,t]$,
then
\[\frac{\norm{Y}_{\phi_t(v)}}{\norm{Z}_{\phi_t(v)}}\le e^{-b_Ct} \frac{\norm{Y}_v}{\norm{Z}_v}\]
for all non-zero $Y\in\xi^{k}(v^+)$ and $Z\in\xi^{d-k}(v^-)$. Hence, if we define
\begin{align*}
B_C = \max\left\{ e^{-b_C t}\frac{ \norm{Z}_v\norm{Y}_{\phi_t(v)}}{\norm{Z}_{\phi_t(v)}\norm{Y}_v} : 0\leq t \leq  T,\,\,\, v \in K,\,\,\, Y,Z \in \Kb^d-\{0\}\right\},
\end{align*}
(notice that $B_C$ is finite by Lemma~\ref{obs:distortion_of_norms_holder}) then the proposition follows.
\end{proof}

It remains to control the behavior of the flow on the thick part. The proof of the following proposition is
inspired by arguments of Tsouvalas \cite[Theorem 1.1]{kostas}.

\begin{proposition}
\label{decay on thick part}
Suppose that $\wh K$ is a compact subset of $\wh\Usf(\Gamma)$, $\{u_n\}$ is a sequence in $\wh K$,  $Y_n\in\wh\Theta^k|_{u_n}$, 
$Z_n\in\wh\Xi^{d-k}|_{u_n}$ and  $\norm{Y_n}_{u_n}=\norm{Z_n}_{u_n}=1$ for all $n$. If $\{t_n\}\subset \mathbb R$, $\lim t_n=+\infty$ and
$\phi_{t_n}(u_n)\in\wh K$ for all $n$, then
$$\lim_{n\to\infty} \frac{\norm{\phi_{t_n}(Y_n)}_{\phi_{t_n}(u_n)}}{\norm{\phi_{t_n}(Z_n)}_{\phi_{t_n}(u_n)}}=0.$$
\end{proposition}

\begin{proof}

We first choose a compact set $K\subset \Usf(\Gamma)$ so that $\wh K\subset\pi(K)$ (where $\pi:\Usf(\Gamma)\to\wh\Usf(\Gamma)$ is the
quotient map). For each $n$, choose $v_n\in K$ so that $\pi(v_n)=u_n$ and 
$\gamma_n\in\Gamma$ so that $w_n=\gamma_n^{-1}(\phi_{t_n}(v_n))\in K$. We may assume that $v_n\to v_\infty$ and $w_n\to w_\infty$ for some $v_\infty,w_\infty\in K$. 
Notice that $\gamma_n\rightarrow v_\infty^+$ and $\gamma_n^{-1}\rightarrow w_\infty^-$. 

Let $V_n \in \xi^k(v_n^+)$ and $W_n \in \xi^{d-k}(v_n^-)$ denote lifts of $Y_n$ and $Z_n$ respectively. Then 
$$\frac{\norm{\phi_{t_n}(Y_n)}_{\phi_{t_n}(u_n)}}{\norm{\phi_{t_n}(Z_n)}_{\phi_{t_n}(u_n)}}=\frac{\norm{V_n}_{\phi_{t_n}(v_n)}}{\norm{W_n}_{\phi_{t_n}(v_n)}}=\frac{\norm{\rho(\gamma_n)^{-1}(V_n)}_{w_n}}{\norm{\rho(\gamma_n)^{-1}(W_n)}_{w_n}},$$
Also, since $K$ is compact, there exists $L$ so that if $v\in K$, then 
$\norm{\cdot}_v$ is $L$-bilipschitz to the standard norm $\norm{\cdot}_2$ on $\Kb^d$. So it suffices to show that 
$$\lim_{n \rightarrow \infty}\frac{\norm{\rho(\gamma_n)^{-1}(V_n)}_2}{\norm{\rho(\gamma_n)^{-1}(W_n)}_2}=0.$$

Since $\xi$ is strongly dynamics preserving, Lemma \ref{lem: basic singular value} implies that
\begin{align}\label{eqn: sing values'}
\lim_{n\to\infty}\frac{\sigma_{d-k}(\rho(\gamma_n)^{-1})}{\sigma_{d-k+1}(\rho(\gamma_n)^{-1})}=\infty,
\end{align}
$U_k(\rho(\gamma_n))\to\xi^{k}(v_\infty^+),$ and $U_{d-k}(\rho(\gamma_n)^{-1})\to\xi^{d-k}(w_\infty^-)$. By the $\rho$-equivariance of $\xi$, $\rho(\gamma_n)^{-1}(W_n) \in\xi^{d-k}(w_n^-)$, which implies that
$$\lim_{n \rightarrow \infty} \angle(U_{d-k}(\rho(\gamma_n)^{-1}),\rho(\gamma_n)^{-1}(W_n)) = 0.$$
Therefore, 
$$\liminf_{n \rightarrow \infty} \frac{\norm{\rho(\gamma_n)^{-1}(W_n)}_2}{\sigma_{d-k}(\rho(\gamma_n)^{-1})}\ge \liminf_{n \rightarrow \infty} \norm{W_n}_2 \geq \frac{1}{L}.$$

We now write $V_n=V_n'+V_n''$, where $V_n'\in U_k(\rho(\gamma_n))$ and $V_n''\in\rho(\gamma_n)\big(U_{d-k}(\rho(\gamma_n)^{-1})\big)$ are orthogonal.
Notice that 
$$\limsup_{n \rightarrow \infty} \frac{\norm{\rho(\gamma_n)^{-1}(V_n')}_2}{\sigma_{d-k+1}(\rho(\gamma_n)^{-1})}\le L$$
so 
$$\lim_{n \rightarrow \infty} \frac{\norm{\rho(\gamma_n)^{-1}(V_n')}_2}{\norm{\rho(\gamma_n)^{-1}(W_n)}_2}= 0.$$
As such, if $\frac{\norm{\rho(\gamma_n)^{-1}(V_n)}_2}{\norm{\rho(\gamma_n)^{-1}(W_n)}_2}$ does not converge to 0 it must be the case
that $\frac{\norm{\rho(\gamma_n)^{-1}(V_n'')}_2}{\norm{\rho(\gamma_n)^{-1}(W_n)}_2}$ does not converge to 0,
and hence that
$$\limsup_{n\to\infty} \frac{\norm{\rho(\gamma_n)^{-1}(V_n'')}_2}{\norm{\rho(\gamma_n)^{-1}(V_n')}_2}=\infty.$$
We may then pass to a subsequence so  that the limits 
\[\lim_{n\to\infty}\rho(\gamma_n)^{-1}\Big(\Span_{\Kb}( V_n'')\Big)=\lim_{n\to\infty}\rho(\gamma_n)^{-1}\Big(\Span_{\Kb}(V_n)\Big)\]
are equal and exist. At the same time,
\[\lim_{n\to\infty}\rho(\gamma_n)^{-1}\Big(\Span_{\Kb}(V_n'')\Big)\subset\lim_{n\to\infty}U_{d-k}(\rho(\gamma_n)^{-1})=\xi^{d-k}(w_\infty^-).\]
and
\[\lim_{n\to\infty}\rho(\gamma_n)^{-1}\Big(\Span_{\Kb}(V_n)\Big)\subset\lim_{n\to\infty}\rho(\gamma_n)^{-1}\big(\xi^{k}(v_n^+)\big)=\lim_{n\to\infty}\xi^{k}(w_n^+)=\xi^{k}(w_\infty^+).\]
This contradicts the transversality of $\xi$, and completes the proof.
\end{proof}

We now combine Proposition \ref{decay on cusps} and Proposition \ref{decay on thick part} to finish the proof of Theorem \ref{thm:equivalent_to_Anosov}. 

For each $C\in\mathcal C$, let $C^\prime\subset C$ be the embedded cusp subneighborhood given by Proposition \ref{decay on cusps}. 
Let $\mathcal{C}^\prime=\{C^\prime : C\in\mathcal C\}$ and let $\wh K=\wh\Usf(\Gamma)-\wh\Usf(\Gamma)_{\mathcal{C}^\prime}$, which is compact.
Let $b=\min\{b_C : C\in\mathcal C\}$ and $B=\max\{B_C : C\in\mathcal C\}$.

Proposition \ref{decay on thick part} implies that there exists $T_1>0$ so that: if $t\ge T_1$, $v\in\wh K$ with $\phi_t(v)\in\wh K$, $Y\in\wh\Theta^k|_v$, and $Z\in\wh\Xi^{d-k}|_v$ is non-zero, then 
\begin{align}\label{eqn: 2B2}
\frac{\norm{\phi_t(Y)}_{\phi_t(v)}}{\norm{\phi_t(Z)}_{\phi_t(v)}}\le \frac{1}{2B^2} \frac{\norm{Y}_v}{\norm{Z}_v}.
\end{align}
Moreover, since $\wh K$ is compact, there exists $R>0$ so that: if  $v\in \wh K$, $0\leq t\leq T_1$, $Y\in\wh\Theta^k|_v$ and $Z\in\wh\Xi^{d-k}|_v$ is non-zero, then
\begin{align}\label{eqn: Ret}
\frac{\norm{\phi_t(Y)}_{\phi_t(v)}}{\norm{\phi_t(Z)}_{\phi_t(v)}}\le Re^t\frac{\norm{Y}_v}{\norm{Z}_v}.
\end{align}

Now choose $T>0$ so that 
$$B^2Re^{T_1(1+b)-bT}\le\frac{1}{2}\,\,\text{ and }\,\,Be^{-bT}\leq\frac{1}{2}.$$
We claim that if $t\ge T$,  $v\in \wh\Usf(\Gamma)$, $Y\in \wh\Theta^k|_v$ and $Z\in \wh\Xi^{d-k}|_v$ is non-zero, then
\begin{align}\label{eqn: half ineq}
\frac{\norm{\phi_t(Y)}_{\phi_t(v)}}{\norm{\phi_t(Z)}_{\phi_t(v)}}\le\frac{1}{2}  \frac{\norm{Y}_v}{\norm{Z}_v}.
\end{align}
Once we have proven our claim, we can choose
\[a:=\frac{\log(2)}{T}\,\text{ and }\,A:=\max\left\{\frac{\norm{\phi_t(Y)}_{\phi_t(v)}\norm{Z}_v}{\norm{\phi_t(Z)}_{\phi_t(v)}\norm{Y}_v}:0\leq t\leq T,\ Y\in \wh\Theta^k|_v-0,\ Z\in \wh\Xi^{d-k}|_v-0\right\},\]
(notice that $A$ is finite by Lemma~\ref{obs:distortion_of_norms_holder}) and observe that
$$\frac{\norm{\phi_t(Y)}_{\phi_t(v)}}{\norm{\phi_t(Z)}_{\phi_t(v)}}\le Ae^{-at} \frac{\norm{Y}_v}{\norm{Z}_v}$$
for all $t>0$, $v\in \wh\Usf(\Gamma)$, $Y\in \wh\Theta^k|_v$ and non-zero $Z\in \wh\Xi^{d-k}|_v$. Hence, we will
have shown that $\rho$ is $P_k$-Anosov.

We now establish \eqref{eqn: half ineq}. Fix $v\in\wh\Usf(\Gamma)$, $t \ge T$, $Y\in \wh\Theta^k|_v$ and non-zero $Z\in \wh\Xi^{d-k}|_v$. 
If $\phi_s(v) \notin \wh K$ for all $s \in [0,t]$, then Proposition \ref{decay on cusps} implies that 
\begin{align*}
\frac{\norm{\phi_t(Y)}_{\phi_t(v)}}{\norm{\phi_t(Z)}_{\phi_t(v)}}\le Be^{-bt}\frac{\norm{Y}_v}{\norm{Z}_v} \leq \frac{1}{2}\frac{\norm{Y}_v}{\norm{Z}_v}.
\end{align*}
Otherwise let
\[ s_1=s_1(t,v):=\min\left\{s\in[0,t] : \phi_s(v)\in \wh K\right\}\]
and
\[ s_2=s_2(t,v):=\max\left\{s\in[0,t] : \phi_s(v)\in \wh K\right\}.\]
If $s_2-s_1\ge T_1$, then Proposition \ref{decay on cusps} and \eqref{eqn: 2B2} give
\begin{align*}
\frac{\norm{\phi_t(Y)}_{\phi_t(v)}}{\norm{\phi_t(Z)}_{\phi_t(v)}}&\le \left(Be^{-b s_1}\right)\ \left(\frac{1}{2B^2}\right) \left( Be^{-b(t-s_2)}\right) \frac{\norm{Y}_v}{\norm{Z}_v}\\
&=\frac{1}{2}e^{b(s_2-s_1-t)} \frac{\norm{Y}_v}{\norm{Z}_v}\le\frac{1}{2}  \frac{\norm{Y}_v}{\norm{Z}_v}.
\end{align*}
On the other hand, if $s_2-s_1<T_1$, then by Proposition \ref{decay on cusps} and \eqref{eqn: Ret}, we have
\begin{align*}
\frac{\norm{Y}_{\phi_t(v)}}{\norm{Z}_{\phi_t(v)}}&\le \left(Be^{-bs_1}\right)\left(Re^{s_2-s_1} \right) \left(Be^{-b(t-s_2)}\right) \frac{\norm{Y}_v}{\norm{Z}_v}\\
&\le\left( B^2Re^{T_1(1+b)-bT}\right) \frac{\norm{Y}_v}{\norm{Z}_v}\le
\frac{1}{2}  \frac{\norm{Y}_v}{\norm{Z}_v}.\qedhere
\end{align*}
\end{proof}

\section{Hitchin representations are Borel Anosov} 

In this section, we show that Hitchin representations are irreducible and Borel Anosov, i.e. $P_k$-Anosov for all $k$.
Theorem \ref{thm:equivalent_to_Anosov} reduces the proof that Hitchin representations are Borel Anosov  to the claim that their limit
maps are strongly dynamics preserving.

\begin{theorem}\label{thm:hitchin_are_str_dyn_preserving} Suppose $\Gamma \subset \PSL(2,\Rb)$ is a geometrically finite group and $\rho : \Gamma \rightarrow \SL(d,\Rb)$ is Hitchin 
with continuous positive $\rho$-equivariant limit map $\xi:\Lambda(\Gamma)\to \mathcal F_d$.  If $\{\gamma_n\}$ is a sequence in $\Gamma$ with $\gamma_n \rightarrow x\in\Lambda(\Gamma) $
and $\gamma_n^{-1} \rightarrow y\in\Lambda(\Gamma)$, then
\begin{align*}
\rho(\gamma_n)(V) \rightarrow \xi^{k}(x)
\end{align*}
for all $V$ transverse to $\xi^{d-k}(y)$.
\end{theorem}

\begin{proof}

Our proof relies on an observation about convergence of sequences of flags.

\begin{definition}\label{defn:open_sets}
For any positive triple of flags $(F_1,F_2,F_3)$ in $\Fc(\Rb^d)$, define the open set
\begin{align}\label{eqn: o notation}
\Oc(F_1,F_2,F_3) := \left\{ F \in \Fc(\Rb^d) : ( F_1,F,F_2,F_3 ) \text{ is positive} \right\}.
\end{align}
\end{definition}

The following result is a corrected version of Lemma 3.16 in \cite{BCL} (whose statement omits
an additional assumption given here).

\begin{lemma} \label{lem: quadruple squeeze}
Let $(F_+,F_-)$ be a transverse pair of flags in $\Fc(\Rb^d)$. Suppose that $\{F_{1,n}\}$, $\{F_{2,n}\}$ and $\{F_{3,n}\}$ are sequences in $\Fc(\Rb^d)$ such that 
\begin{enumerate}
\item $F_{1,n}\to F^+$ and $F_{2,n}\to F^+$, 
\item $F_{3,n}\to F^-$, and 
\item $(F_{1,n}, F_{2,n}, F_{3,n})$ is positive for all $n$. 
\end{enumerate}
If $\{F_n\}$ is a sequence in $\Fc(\Rb^d)$ such that $F_n\in\Oc(F_{1,n},F_{2,n},F_{3,n})$ for all $n$, then $F_n\to F^+$.
\end{lemma}

We first suppose that $x \neq y$. Then by passing to the tail of the sequence $\{\gamma_n\}$, we may assume that each $\gamma_n$ is hyperbolic with attractor and repellor $\gamma_n^+$ and $\gamma_n^-$ in $\Lambda(\Gamma)$. Then $\gamma_n^+\rightarrow x$, $\gamma_n^- \rightarrow y$, and $\gamma_n(z) \rightarrow x$ for all $z \in \Lambda(\Gamma) - \{y\}$.

Since $\Lambda(\Gamma)$ is infinite, there are points $a,b \in \Lambda(\Gamma)-\{x,y\}$ such that 
\begin{itemize}
\item either $x<a<b<y$ or $y<b<a<x$, and
\item up to taking subsequences, the sequences $\{\gamma_n(a)\}$ and $\{\gamma_n(b)\}$ both converge monotonically to $x$, and from the same direction. 
\end{itemize}
Observe that
\begin{align*}
 \lim_{n \rightarrow \infty}\xi(\gamma_n(a)) =  \xi(x) =  \lim_{n \rightarrow \infty} \xi(\gamma_n(b)).
\end{align*}

Now, consider the open sets
 \begin{align*}
 \Oc_n:=\Oc\left(\xi(a), \xi(b), \xi(\gamma_n^-)\right)
 \end{align*}
 for all $n$. Since $\gamma_n^- \rightarrow y$, and either $a<b<y$ or $b<a<y$, it follows that there is some $N>0$ such that either $a<b<\gamma_n^-$ for all $n\geq N$, or $b<a<\gamma_n^-$ for all $n\geq N$. Lemma 3.15 in \cite{BCL}  then implies that $\Oc_n=\Oc_m$ for all $n,m\geq N$. Hence, if we set $\Oc:= \Oc_N$, then for all $n\geq N$, we have
 \begin{align*}
 \rho(\gamma_n)\left(\Oc\right) = \rho(\gamma_n)\left(\Oc_n\right)= \Oc(\xi(\gamma_n(a)), \xi(\gamma_n( b)), \xi(\gamma_n^-)).
 \end{align*}

Since $\xi(\gamma_n(a)), \xi(\gamma_n(b)) \rightarrow \xi(x)$ and $\xi(\gamma_n^-) \rightarrow \xi(y)$, we may apply Lemma \ref{lem: quadruple squeeze} to deduce that
 \begin{align*}
 \lim_{n \rightarrow \infty} \rho(\gamma_n) (F) = \xi(x)
 \end{align*}
 for all $F \in \Oc$. Repeating the same argument with $\gamma_n^{-1}$, we see that there exists an open set $\Oc^\prime \subset \Fc(\Rb^d)$ where
  \begin{align*}
 \lim_{n \rightarrow \infty} \rho(\gamma_n^{-1})(F) = \xi(y)
 \end{align*}
 for all $F \in \Oc^\prime$. Hence, we may apply Lemma \ref{lem: basic singular value} to deduce the proposition when $x\neq y$.

Now suppose that $x = y$. Pick $\gamma \in \Gamma$ such that $z:=\gamma^{-1} (x) \neq x$. Then $\gamma_n\gamma \rightarrow x$, 
$(\gamma_n\gamma)^{-1} \rightarrow z \neq x$. 
By the first part, $\rho(\gamma_n\gamma)(F) \rightarrow \xi(x)$ for all $F \in \Fc(\Rb^d)$ transverse to $\xi(z)$. Equivalently, $\rho(\gamma_n)( F) \rightarrow \xi(x)$ for all $F \in \Fc(\Rb^d)$ transverse to $\xi(x)$. 
\end{proof}

We recall that positive tuples of flags  are in general position in the following sense.

\begin{proposition}\label{prop: general position} {\em (Fock-Goncharov \cite[Prop. 9.4]{fock-goncharov}, Sun-Wienhard-Zhang \cite[Prop 2.21]{SWZ})}
If $(F_1,\dots,F_k)$ is a positive tuple of flags, $(n_i)_{i=1}^k\in\mathbb N^k$ and $n=n_1+\cdots +n_k\le d$, then $\oplus_{i=1}^k F_i^{n_i}$
has dimension $n$.
\end{proposition}

We also use the following equivalent formulation of the positivity of a quadruple of flags.

\begin{lemma}\label{lem: posequiv}
A quadruple of flags $(F_1,F_2,F_3,F_4)$ is positive if and only if there is a basis $(b_1,\dots,b_d)$ of $\Rb^d$ such that $b_i\in F_1^{i}\cap F_3^{d-i+1}$ for all $i\in\{1,\dots,d\}$, and some $u,v\in U_{>0}(b_1,\dots,b_d)$ 
such that $u(F_3)=F_2$ and $v^{-1}(F_3)=F_4$.
\end{lemma}

\begin{proof}
Suppose first that $(F_1,F_2,F_3,F_4)=(F_1,u(F_3),F_3,v^{-1}(F_3))$ for some basis $(b_1,\dots,b_d)$ of $\Rb^d$ such that $b_i\in F_1^{i}\cap F_3^{d-i+1}$ for all $i\in\{1,\dots,d\}$, and some $u,v\in U_{>0}(b_1,\dots,b_d)$. Then $v(F_1,F_2,F_3,F_4)=(F_1,vu(F_3),v(F_3),F_3)$, which implies that $v(F_1,F_2,F_3,F_4)$ is positive. Thus, $(F_1,F_2,F_3,F_4)$ is positive.

Conversely, suppose that $(F_1,F_2,F_3,F_4)$ is positive. By Proposition \ref{prop: general position}, $F_4$ and $F_3$ are both transverse to $F_1$, there is a unique unipotent $w\in\SL(d,\Rb)$ that fixes $F_1$ and sends $F_4$ to $F_3$. Then $w(F_1,F_2,F_3,F_4)$ is positive, which implies that there is a basis $(b_1,\dots,b_d)$ of $\Rb^d$ such that $b_i\in F_1^{i}\cap F_3^{d-i+1}$ for all $i\in\{1,\dots,d\}$, and some $u,v\in U_{>0}(b_1,\dots,b_d)$ such that 
\[(F_1,w(F_2),w(F_3),F_3)=w(F_1,F_2,F_3,F_4)=(F_1,vu(F_3),v(F_3),F_3).\]
Since $v$ is unipotent, fixes $F_1$, and sends $F_3$ to $w(F_3)$, it follows that $v=w$. Therefore, 
\[(F_1,F_2,F_3,F_4)=v^{-1}(F_1,vu(F_3),v(F_3),F_3)=(F_1,u(F_3),F_3,v^{-1}(F_4)).\qedhere\]
\end{proof}

\begin{proof}[Proof of Theorem \ref{thm: Hitchin}]
The fact that $\rho$ is $P_k$-Anosov, and that $x\mapsto\xi^{k}(x)$ is the Anosov limit map for $k=1,\dots,d-1$, follows from Theorems
\ref{thm:equivalent_to_Anosov} and  \ref{thm:hitchin_are_str_dyn_preserving}. Further, Theorem~\ref{thm:hitchin_are_str_dyn_preserving} and
Lemma \ref{lem: basic singular value}  imply that:
\begin{enumerate}
\item If $\alpha \in \Gamma$ is parabolic, then $\rho(\alpha)$ is weakly unipotent.
\item If $\gamma \in \Gamma$ is hyperbolic, then $\rho(\gamma)$ is loxodromic. 
\end{enumerate}

If  $\alpha$ is parabolic, let $x\in\Lambda(\Gamma)-\{\alpha^+\}$, and note that $(\alpha^+,\alpha^{-1}(x),x,\alpha(x))$ is a cyclically ordered set of distinct points in $\Lambda(\Gamma)$. Then we may
apply the following lemma to $(\xi(\alpha^+),\xi(\alpha^{-1}(x)),\xi(x),\xi(\alpha(x))$ to further conclude that $\rho(\alpha)=\pm u$ for some unipotent $u\in\SL(d,\Rb)$ with a single Jordan block.

\begin{lemma}\label{lem: peripheral parabolic}
Let $g\in\SL(d,\Rb)$ be weakly unipotent, and suppose that there are flags $F_1,F_2\in\Fc(\Rb^d)$ such that $F_1$ is fixed by $g$ and $(F_1,g^{-1}(F_2),F_2,g(F_2))$ is positive. Then $F_1$ is the unique fixed flag of $g$, or equivalently, $g=\pm u$ for some unipotent $u\in\SL(d,\Rb)$ with a single Jordan block.
\end{lemma}

\begin{proof}
By Lemma \ref{lem: posequiv}, there is a basis $(b_1,\dots,b_d)$ of $\Rb^d$ such that $b_i\in F_1^{i}\cap F_2^{d-i+1}$ for all $i\in\{1,\dots,d\}$, and some $u,v\in U_{>0}(b_1,\dots,b_d)$ 
such that $u(F_2)=g(F_2)$ and $v^{-1}(F_2)=g^{-1}(F_2)$. Then $a=u^{-1}g$ and $b=gv^{-1}$ both fix $F_1$ and $F_2$, so they are diagonal in
the basis $(b_1,\dots,b_d)$. Furthermore, since $g$ is weakly unipotent, the diagonal entries of $a$ and $b$ are either $1$ or $-1$. 

Assume for contradiction that there is some $i,j\in\{1,\dots,d\}$ such that the $i$-th diagonal entry of $b$ is $1$, while the $j$-th diagonal entry of $b$ is $-1$. Since all the upper triangular entries of $v$ are positive, this implies that the upper triangular entries of $bv$ along the $i$-th row are positive, while the upper triangular entries of $bv$ along the $j$-th row are negative. But this is impossible since $ua=g=bv$, and for every column of $ua$, the upper triangular entries in that column must have the same sign. As such, $b=\pm\id$. This implies that $g_{ss}=\pm\id$ and $g_u=v=u$.

It now suffices to show that $u$ has a unique fixed flag in $\Fc(\Rb^d)$. Observe the following linear algebra facts:
\begin{enumerate}
\item If $w$ is a unipotent element that is represented in a basis $(e_1,\dots,e_d)$ by upper triangular matrix where all the upper triangular entries are positive, then the line spanned by $e_1$ is the unique fixed line of $w$. 
\item If $w\in U_{>0}(b_1,\dots,b_d)$, then for all $k\in\{1,\dots,d-1\}$, the linear action of $w$ on $\bigwedge^{k}\Rb^d$ is represented in the basis $(b_{i_1}\wedge\dots \wedge b_{i_k})_{1\leq i_1<\dots<i_k\leq d}$ by an upper triangular matrix where all the upper triangular entries are positive.
\end{enumerate}
These two observations imply that the unique fixed flag of $u$ is the flag $F$ given by $F^k=\Span_{\Rb}(b_1,\dots,b_k)$ for all $k\in\{1,\dots,d-1\}$. \end{proof}

It only remains to show that $\rho$ is irreducible. Suppose that $\rho$ is not irreducible. Then there is a proper subspace $W\subset \Rb^d$ which is invariant under $\rho(\Gamma)$. 
By Theorem ~\ref{thm:hitchin_are_str_dyn_preserving}, $\rho(\gamma)$ is loxodromic for any hyperbolic $\gamma\in\Gamma$, so either $W$ contains the attracting fixed point (in $\Pb(\Rb^d)$) of $\rho(\gamma)$, or $W$ lies in the repelling hyperplane (in $\Pb(\Rb^d)$) of $\rho(\gamma)$. Since $\xi$ is transverse, this implies that either $W$ lies in the repelling hyperplanes of $\rho(\gamma)$ for all hyperbolic $\gamma\in\Gamma$, or $W$ contains the attracting fixed point of $\rho(\gamma)$ for all hyperbolic $\gamma\in\Gamma$. However, this contradicts Proposition \ref{prop: general position} in either case.
Therefore, $\rho$ is irreducible.
\end{proof}


\section{Stability of Anosov representations}

In this section, we prove Theorem \ref{thm: stability intro}, which we restate here. 

\begin{theorem}\label{thm: stability}
If $\Gamma \subset \PSL(2,\Rb)$ is a geometrically finite group and $\rho_0:\Gamma\to \mathsf{SL}(d,\mathbb K)$ is $P_k$-Anosov, then
there exists an open neighborhood $\Oc$ of $\rho_0$ in $\mathrm{Hom}_{\rm tp}(\rho_0)$, so that 
\begin{enumerate}
\item
If $\rho\in \Oc$, then $\rho$ is $P_k$-Anosov.
\item
There exists $\alpha>0$ so that if  $\rho\in \Oc$, then $\xi_\rho$ is $\alpha$-H\"older (with respect to any distance on $\Lambda(\Gamma)$ induced by a Riemannian metric on $\partial\Hb^2$ and any distance on $\mathrm{Gr}_k(\mathbb K^d)\times\mathrm{Gr}_{d-k}(\mathbb K^d)$ induced by a Riemannian metric).
\item
If  $\{\rho_u\}_{u\in M}$ is a $\mathbb K$-analytic family of representations in $\Oc$ and $z\in\Lambda(\Gamma)$
then the map from $M$ to  $\mathrm{Gr}_{k}(\mathbb K^d)\times\mathrm{Gr}_{d-k}(\mathbb K^d)$
given by $u\mapsto \xi_{\rho_u}(z)$ is $\mathbb K$-analytic.
\end{enumerate}
\end{theorem}

The proof of (1) is based on the proof of stability for Anosov diffeomorphisms on compact manifolds given in Shub's book \cite[Cor. 5.19]{shub-book}. The two key features which allow 
us to overcome the non-compactness of the base space are the smooth conjugacy of the flows on the cusps, see Equation~\eqref{eqn:smooth_conj}, 
and uniform estimates for families of canonical norms, see Lemma~\ref{obs:distortion_of_norms_holder}. 

\begin{proof}[Proof of Theorem \ref{thm: stability}:]
We define
$$E(\Oc)=\Oc\times\Usf(\Gamma)\times\mathbb K^d\,\,\text{ and }\,\,
\wh E(\Oc)=\Gamma\backslash(\Oc\times\Usf(\Gamma)\times \mathbb K^d).$$
The geodesic flow on $\Usf(\Gamma)$ extends to a flow on $E(\Oc)$ whose action is trivial on the first and third factor. This in turn descends to a flow on $\wh E(\Oc)$. As usual, we use $\phi_t$ to denote these flows. Also, notice that $\wh E(\Oc)|_{\rho}$ naturally identifies with $\wh E_\rho$.

Let $\norm{\cdot}^0$ be a canonical family of norms for $\wh E_{\rho_0}$ and let $\mathcal C$ be a full collection of embedded cusp neighborhoods for
$\Gamma$ so that $\norm{\cdot}^0$ is canonical with respect to $\mathcal C$. Suppose that $\wh E_{\rho_0}=\wh\Theta_{\rho_0}^k\oplus\wh\Xi_{\rho_0}^{d-k}$ is the $P_k$-Anosov splitting of $\wh E_{\rho_0}$. 

\subsection*{A $\phi_t$-invariant splitting of $\wh E(\Oc)$.}
First, we prove that (after possibly shrinking $\Oc$) there exists a continuous $\phi_t$-invariant splitting
\[\wh E(\Oc)=\wh\Theta^k\oplus\wh\Xi^{d-k}\]
that restricts to the splitting $\wh E_{\rho_0}=\wh\Theta_{\rho_0}^k\oplus\wh\Xi_{\rho_0}^{d-k}$ over $\rho_0$. 

If $C\in\mathcal C$ and $C= \langle\alpha\rangle\backslash H$, then by shrinking $\Oc$ if necessary, we may assume that there is a continuous map $g_C:\Oc\to\SL(d,\Kb)$ such that 
\[g_C(\rho)\rho_0(\alpha)g_C(\rho)^{-1}=\rho(\alpha)\] 
for all $\rho\in\Oc$.
Moreover, if $\rho\in\Oc$, the bundle isomorphism
$$\Phi_\rho^H: E_{\rho_0}|_{\Usf(\Gamma)_H}\to E_\rho|_{\Usf(\Gamma)_H} \ \ \text{given by}\ \  (\rho_0,v,Z)\to (\rho,v,g_C(\rho)(Z))$$
descends to a bundle isomorphism
$$\wh \Phi_\rho^C:\wh E_{\rho_0}|_{\wh\Usf(\Gamma)_C}\to \wh E_\rho|_{\wh\Usf(\Gamma)_C}$$
so that if $\phi_s(Z)\in \wh E_{\rho_0}|_{\wh\Usf(\Gamma)_C}$ for all $s\in[0,t]$, then
\begin{equation}
\label{eqn:smooth_conj}
\wh\Phi_\rho^C( \phi_t(Z)) = \phi_t (\wh\Phi_\rho^C(Z)).
\end{equation}
With this, we may extend the splitting $\wh E_{\rho_0}=\wh\Theta_{\rho_0}^k\oplus\wh\Xi_{\rho_0}^{d-k}$ to a global splitting
$$\wh E(\Oc)=\wh F^k\oplus\wh G^{d-k}$$
by first setting
$$\wh F^k|_{(\rho,v)}=\wh\Phi_\rho^C(\wh\Theta^k_{\rho_0}|_v)\ \ \text{and}\ \  \wh G^{d-k}|_{(\rho,v)}=\wh\Phi_\rho^C(\wh\Xi^{d-k}_{\rho_0}|_v)$$
for all $\rho\in\Oc$, $C\in\Cc$ and $v\in\Usf(\Gamma)_C$, and then extending this globally after perhaps shrinking $\Oc$ and each $C$. 

The flow $\phi_t$ does not necessarily preserve the splitting $\wh E(\Oc)=\wh F^k\oplus\wh G^{d-k}$. To find the required $\phi_t$-invariant splitting, we will use the contraction mapping theorem. For that purpose, we extend $\norm{\cdot}^0$ to a canonical family of norms
$\norm{\cdot}$ on  the fibers of $\wh E(\Oc)$ over $\Oc\times\wh\Usf(\Gamma)$ as follows.
If $C\in\mathcal C$, we define
$$\norm{\wh\Phi_\rho^C(Z)}_{(\rho,v)}=\norm{Z}^0_{(\rho_0,v)}\ \ 
\text{for all}\ \  \rho\in\Oc,\ \ v\in\wh\Usf(\Gamma)_C\  \text{and}\   Z\in\wh E_{\rho_0}|_v.$$
This gives us canonical norms over all $C\in\mathcal C$. Then, perhaps after once more shrinking $\Oc$ and each $C$, we may extend this to a continuous family of 
norms $\norm{\cdot}$ for the fibers of $\wh E(\Oc)$ over $\Oc\times\wh\Usf(\Gamma)$ such that the restriction to $\wh E_\rho$ is canonical
for all $\rho\in\Oc$.
 
Suppose that $V$ is a subspace of $\wh E(\Oc)|_{(\rho,v)}$, $W$ is a subspace of $\wh E(\Oc)|_{(\rho,w)}$ and $T\in \mathrm{Hom}(V,W)$. We define the operator norm
\begin{align*}
\norm{T}_{(\rho,v)} := \max\left\{ \norm{T(Z)}_{(\rho,w)} : Z \in V,\ \norm{Z}_{(\rho,v)}=1\right\}.
\end{align*}
Then, if $V,W\subset\wh E(\Oc)$ are subbundles, $Q:\wh\Usf(\Gamma)\to\wh\Usf(\Gamma)$ is an isomorphism, and $T:V\to W$ is a map that restricts to a linear map $T|_{(\rho,v)}:V|_{(\rho,v)}\to W|_{(\rho,Q(v))}$ for all $\rho\in\Oc$ and $v\in\wh\Usf(\Gamma)$, we define
\[\norm{T}_X:=\sup_{\rho\in\Oc,v\in X}\norm{T|_{(\rho,v)}}_{(\rho,v)}\]
for any subset $X \subset \wh\Usf(\Gamma)$. In the case when $Q=\id$, we may view $T$ as an element of $S(\Hom(V,W))$, the vector space of sections of the bundle $\Hom(V,W)$. Note then that $\norm{\cdot}_{\wh\Usf(\Gamma)}$ defines a norm on $S(\Hom(V,W))$ whose corresponding distance is Cauchy complete. 

We may decompose the flow $\phi_t$ as 
$$\phi_{t} = \begin{pmatrix} A_t & B_t \\ C_t & D_t \end{pmatrix}$$
relative to the splitting $\wh{E}(\Oc) =\wh F^k\oplus\wh G^{d-k}$. Here, 
$A_t:\wh F^k\to\wh F^k$ is a map such that for all $\rho\in\Oc$ and $v\in\wh\Usf(\Gamma)$, $A_t$ restricts to a linear map $A_{(\rho,v,t)}:\wh F^k|_{(\rho,v)}\to \wh F^k|_{(\rho,\phi_t(v))}$, etc. In the case when this splitting happens to be $\phi_t$-invariant, then $B_t \equiv 0$ and $C_t \equiv 0$. 

Fix $\epsilon \in (0,1/2)$ so that
\begin{equation}
\label{eqn:epsilon is super small}
 \frac{1+\epsilon}{(1-\epsilon)} \leq 2, \quad \frac{1}{1-\epsilon} + \frac{\epsilon(1+\epsilon)}{(1-\epsilon)^2} \leq 2 \quad \text{and} \quad \epsilon\frac{(1+2\epsilon)^2(1+2\epsilon^2)}{(1-2\epsilon)^2(1-2\epsilon^2)} \leq \frac{1}{2}.
\end{equation}

\begin{lemma}
\label{lem:bds on operators}
Up to taking a subneighborhood of $\Oc$, there exists $T > 0$ so that if $t \in [T,2T]$, 
 then $A_t$ and $D_t$ are invertible and 
\begin{align*}
\max\left\{ \norm{A_t^{-1}B_t}, \norm{D_t^{-1}C_t}, L_t\right\} < \epsilon,
\end{align*}
where $L_t:=\sup_{\rho\in\Oc,v\in\wh\Usf(\Gamma)}\left(\norm{A_{(\rho,v,t)}}_{(\rho,v)}\norm{D_{(\rho,v,t)}^{-1}}_{(\rho,\phi_t(v))}\right)$.
\end{lemma} 

\begin{proof} 
For a subset $X \subset \wh\Usf(\Gamma)$, define 
$$
L_t(X):=\sup_{\rho\in\Oc,v\in X}\left(\norm{A_{(\rho,v,t)}}_{(\rho,v)}\norm{D_{(\rho,v,t)}^{-1}}_{(\rho,\phi_t(v))}\right).$$

First notice that for all $v\in\wh\Usf(\Gamma)$, $B_{(\rho_0,v,t)}=0$, $C_{(\rho_0,v,t)}=0$ and both $A_{(\rho_0,v,t)}$, $D_{(\rho_0,v,t)}$ are invertible since $F^k \oplus G^{d-k}|_{\rho_0} = \wh \Theta_{\rho_0}^k\oplus\wh \Xi_{\rho_0}^{d-k}$ is a flow invariant splitting. 
Since $\rho_0$ is $P_k$-Anosov, there exist $c,C > 0$ such that
\begin{align*}
\norm{A_{(\rho_0,v,t)}}_{(\rho_0,v)}\norm{D_{(\rho_0,v,t)}^{-1}}_{(\rho_0,\phi_t(v))} \le C e^{-ct}.
\end{align*}
for all $v\in\wh\Usf(\Gamma)$ and $t>0$. Choose $T$ so that $Ce^{-cT}<\epsilon$.

Consider the compact set
\begin{align*}
X= \left\{ v \in \wh{\Usf}(\Gamma)\  :\  \phi_t(v) \notin \wh{\Usf}(\Gamma)_{\Cc} \text{ for some } t \in [0,2T]\right\}.
\end{align*}
By shrinking $\Oc$ if necessary, we can ensure that on $\wh E(\Oc)|_{\Oc\times X}$, if $t \in [T,2T]$, then $A_t$ and $D_t$ are invertible, and
\begin{align*}
\max\left\{ \norm{A_t^{-1}B_t}_X, \norm{D_t^{-1}C_t}_X, L_t(X)\right\} < \epsilon.
\end{align*}

On the other hand, if $v \in \wh{\Usf}(\Gamma)-X$, then there is some $C\in\Cc$ such that $\phi_t(v)\in\wh\Usf(\Gamma)_C$ for all $t\in[0,2T]$. Then by construction,
$$\phi_{(\rho,v,t)}=\wh\Phi_\rho^C\circ \phi_{(\rho_0,v,t)}\circ(\wh\Phi_\rho^C)^{-1}$$
for all $\rho\in\Oc$ and $t \in [0,2T]$. Since $\wh\Phi_\rho^C$ is an isometry that preserves the splitting, it
follows that on $\wh E(\Oc)|_{\Oc\times (\wh{\Usf}(\Gamma)-X)}$, if $t \in [T,2T]$, then $A_t$ and $D_t$ are invertible, $B_t=0$, $C_t=0$, and
\begin{align*}
L_t\left(\wh\Usf(\Gamma)-X\right) = \sup_{v \in \wh\Usf(\Gamma)-X} \norm{A_{(\rho_0,v,t)}}_{(\rho_0,v)}\norm{D_{(\rho_0,v,t)}^{-1}}_{(\rho_0,\phi_t(v))}  < \epsilon.
\end{align*}
\end{proof}

Consider the bundle $\Hom(\wh G^{d-k},\wh F^k) \rightarrow \Oc \times \wh{\Usf}(\Gamma)$
with its induced operator norm $\norm{\cdot}$. Let $\mathcal{R}_r \subset \Hom(\wh G^{d-k}, \wh F^k)$ denote the ball bundle of radius $r$ about the zero section. 

\begin{proposition}\label{prop:disk_bundle} 
If $t \in [T,2T]$, then there is a well-defined map $\psi_t:\mathcal R_1\to\mathcal R_{2\epsilon}$ given by
\begin{align*}
\psi_t(f)= \left(B_t+ A_t f\right)\left(D_t + C_t f\right)^{-1}
\end{align*}
for all $\rho\in\Oc$, $v\in\wh\Usf(\Gamma)$ and $f\in\mathcal R_1|_{(\rho,v)}$. Furthermore:
$$
\norm{\psi_t(f_1)-\psi_t(f_2)} \leq 2\epsilon \norm{f_1-f_2}
$$
for all $\rho\in\Oc$, $v\in\wh\Usf(\Gamma)$ and $f_1,f_2\in\mathcal R_1|_{(\rho,v)}$.

\end{proposition} 

\begin{remark}\label{rem: phitpsit}
One can verify that the map $\psi_t$ has the defining property 
\[{\rm Graph}(\psi_t(f)) =\phi_t( {\rm Graph}(f))\] 
for all $t\in[T,2T]$, $\rho\in\Oc$, $v\in\wh\Usf(\Gamma)$ and $f\in\mathcal R_1|_{(\rho,v)}$. Similarly, if $B_t \equiv 0$ and $C_t \equiv 0$, then 
\begin{align*}
\psi_t(f)= A_t f D_t^{-1} = \phi_t \circ f \circ \phi_{-t}
\end{align*}
is a well defined flow on $\Hom(\wh G^{d-k},\wh F^k)$. 
\end{remark}

\begin{proof}[Proof of Proposition \ref{prop:disk_bundle} ] If $t \in [T,2T]$, then
\begin{align*}
D_t + C_t f= D_t\left(\id + D_t^{-1} C_t f\right).
\end{align*}
for all $\rho\in\Oc$, $v\in\wh\Usf(\Gamma)$ and $f\in\mathcal R_1|_{(\rho,v)}$. By Lemma~\ref{lem:bds on operators}, if $\norm{f}_{(\rho,v)} < 1$, then 
\[
\norm{D_t^{-1} C_t f}_{(\rho,v)}\leq\norm{D_t^{-1} C_t}_{\wh{\Usf}(\Gamma)}\norm{ f}_{(\rho,v)}<\epsilon
\]
 which implies that $\id + D_t^{-1} C_t f$ has trivial kernel for all $\rho\in\Oc$ and $v\in\wh\Usf(\Gamma)$. Hence, $\id + D_t^{-1} C_t f$ is invertible. Since Lemma~\ref{lem:bds on operators} also gives that $D_t$ is invertible, it follows that
\[\left(B_t+ A_t f\right)\left(D_t + C_t f\right)^{-1}\in\Hom(\wh G^{d-k},\wh F^k)|_{(\rho,\phi_t(v))}\]
 is a well-defined for all $f\in\mathcal R_1|_{(\rho,v)}$.
 
We first show that $\left(B_t+ A_t f\right)\left(D_t + C_t f\right)^{-1}\in\mathcal R_{2\epsilon}|_{(\rho,\phi_t(v))}$. By Lemma~\ref{lem:bds on operators}, 
\[
\norm{A_t^{-1}B_t+f}_{(\rho,v)}< 1 + \epsilon
\]
and 
\[\norm{ (\id+D_t^{-1} C_tf)^{-1}}_{(\rho,v)}<\frac{1}{1-\epsilon}.\]
Thus, by Equation~\eqref{eqn:epsilon is super small}, 
\begin{align*}
\norm{\left(B_t+ A_t f\right)\left(D_t + C_t f\right)^{-1}}_{(\rho,\phi_t(v))} & \leq L_t \norm{ A_t^{-1} B_t + f}_{(\rho,v)} \norm{(\id + D_t^{-1}C_t f)^{-1}}_{(\rho,v)}\\
& < \epsilon \frac{1+\epsilon}{1-\epsilon} \leq 2 \epsilon.
\end{align*}

Next, we prove our final claim. For any $\rho\in\Oc$, $v\in\wh\Usf(\Gamma)$, $f\in\mathcal R_1|_{(\rho,v)}$ and  $\eta \in \Hom(\wh G^{d-k},\wh F^k)|_{(\rho,v)}$, 
\begin{align*}
\left.\frac{d}{ds}\right|_{s=0} &(\psi_t)(f+s\eta)   = A_t \eta (D_t+ C_t f)^{-1} -(B_t + A_t f) (D_t + C_t f)^{-1} C_t \eta(D_t + C_t f)^{-1} \\
& = A_t\Big(  \eta( \id + D_t^{-1} C_t f)^{-1} -(A_t^{-1} B_t+f)(\id + D_t^{-1}C_t f)^{-1} D_t^{-1} C_t \eta (\id + D_t^{-1}C_t f)^{-1}\Big) D_t^{-1}\\
& \in \Hom(\wh G^{d-k},\wh F^k)|_{(\rho,\phi_t(v))}.
\end{align*}
Thus, by Equation~\eqref{eqn:epsilon is super small}, 
\begin{align*}
\norm{\left.\frac{d}{ds}\right|_{s=0}(\psi_t)(f+s\eta)}_{(\rho,\phi_t(v))} \leq \epsilon \left( \frac{1}{1-\epsilon} + \frac{\epsilon(1+\epsilon)}{(1-\epsilon)^2} \right)\norm{\eta}_{(\rho,v)} =2 \epsilon\norm{\eta}_{(\rho,v)}.
\end{align*}

\end{proof} 

Let $S(\mathcal R_r)$ be the space of continuous sections of $\mathcal R_r\to \Oc\times\Usf(\Gamma)$. Notice that $\psi_t$ induces a map $\psi^S_t : S(\mathcal R_1) \rightarrow S(\mathcal R_{2\epsilon})$ given by 
\begin{align*}
\psi_t^S(\sigma)(\rho,v) = \psi_t\left( \sigma(\rho, \phi_{-t}(v))\right).
\end{align*}
By Proposition \ref{prop:disk_bundle}, the map $\psi_t^S$ is a contraction mapping on $S(\mathcal R_1)$ for each $t\in[T,2T]$. We may now apply the contraction mapping theorem to conclude that 
for each $t \in [T,2T]$ there exists a unique $\psi_{t}^S$-invariant section $\sigma^{(t)}$ of the bundle $\mathcal{R}_{2\epsilon}$. 

We claim that $\sigma^{(t)}$ does not depend on $t$. If $t_1, t_2 \in \Qb \cap [T,2T]$, then there exist sequences $\{n_j\}, \{m_j\}$ with $n_j,m_j \rightarrow \infty$ and $n_j t_1 = m_j t_2$ for all $j \geq 1$. Then by the proof of the contraction mapping theorem
\begin{align*}
\sigma^{(t_1)} = \lim_{j \rightarrow \infty} \psi^S_{n_j t_1}\left(\sigma^{(t_2)}\right) = \lim_{j \rightarrow \infty} \psi^S_{m_j t_2}\left(\sigma^{(t_2)}\right)= \sigma^{(t_2)}. 
\end{align*}
So $\sigma^{(t)}$ does not depend on $t$ when $t \in \Qb \cap [T,2T]$. Then by uniqueness of invariant sections and the continuity of $\psi_t^S$, we see that  $\sigma^{d-k}: = \sigma^{(t)}$ does not depend on $t$. Then $\sigma^{d-k}$ determines a $\phi_t$-invariant $(d-k)$-dimensional subbundle  $\wh\Xi^{d-k}$ of $\wh{E}(\Oc)$ defined by 
\begin{align*}
\wh\Xi^{d-k}|_{(\rho,v)} = {\rm Graph}\, \sigma^{d-k}(\rho,v),
\end{align*}
see Remark \ref{rem: phitpsit}.

Applying a similar argument to the bundle $\Hom(\wh F^k,\wh G^{d-k})$ we obtain, by further shrinking $\Oc$ if necessary, a $\phi_t$-invariant $k$-dimensional subbundle $\wh\Theta^k$ of $\wh{E}(\Oc)$.

To show that $\wh E(\Oc)=\wh\Theta^k\oplus\wh\Xi^{d-k}$, it now suffices to show that the fibers $\wh\Xi^{d-k}|_{(\rho,v)}$ and $\wh\Theta^k|_{(\rho,v)}$ are transverse for every $(\rho,v)\in\Oc\times\wh\Usf(\Gamma)$. Suppose for contradiction that there is some non-zero $Z\in\wh\Xi^{d-k}|_{(\rho,v)}\cap\wh\Theta^k|_{(\rho,v)}$ for some $\rho\in\Oc$ and $v\in\wh\Usf(\Gamma)$. We may write $Z$ uniquely as $Z_k+Z_{d-k}$ where $Z_k\in \wh F^k|_{(\rho,v)}$ and $Z_{d-k}\in \wh G^{d-k}|_{(\rho,v)}$. Since $\wh\Xi^{d-k}$ corresponds to the section $\sigma^{d-k} \in S(\mathcal R_{2\epsilon})$, the fact that $Z\in\wh\Xi^{d-k}|_{(\rho,v)}$ is non-zero implies that $Z_{d-k}\neq 0$ and 
$$
\frac{\norm{Z_k}_{(\rho,v)}}{\norm{Z_{d-k}}_{(\rho,v)}}=\frac{\norm{\sigma^{d-k}(\rho,v)(Z_{d-k})}_{(\rho,v)}}{\norm{Z_{d-k}}_{(\rho,v)}}\leq\norm{\sigma^{d-k}(\rho,v)}_{(\rho,v)}<2\epsilon < 1.
$$ 
For the same reasons, $Z_{k}\neq 0$ and
 $$
 \frac{\norm{Z_{d-k}}_{(\rho,v)}}{\norm{Z_{k}}_{(\rho,v)}}<2\epsilon < 1,
 $$
which is a contradiction. Thus, $\wh\Theta^k$ and $\wh\Xi^{d-k}$ indeed give a $\phi_t$-invariant splitting of $\wh E(\Oc)$.

\begin{proposition}\label{prop: contracting Hom bundle}
The flow $f \mapsto \phi_t \circ f \circ \phi_{-t}$ on $\Hom(\wh\Xi^{d-k},\wh\Theta^k)$ is uniformly contracting.
\end{proposition}

\begin{proof} We start by proving the  following estimate:
\begin{equation}
\label{eqn:1/2 contraction after time T}
\frac{\norm{\phi_t(Y)}_{(\rho, \phi_t(v))}}{\norm{\phi_t(Z)}_{(\rho, \phi_t(v))}} \leq \frac{1}{2}\frac{\norm{Y}_{(\rho, v)}}{\norm{Z}_{(\rho,v)}}
\end{equation}
 if $\rho \in \Oc$, $v \in \wh{\Usf}(\Gamma)$, $t \in [T,2T]$, $Y \in \wh\Theta^k|_{(\rho,v)}$ and $Z \in \wh\Xi^{d-k}|_{(\rho,v)}$ is non-zero.
Given such $\rho$, $v$, $t$, $Y$ and $Z$, let 
 $$
 Y=Y_1 + Y_2 \quad \text{and} \quad Z=Z_1 + Z_2
 $$
 be the decomposition relative to $\wh{E}(\Oc) = \wh{F}^{k} \oplus \wh{G}^{d-k}$. Then, by the construction of $\wh{\Theta}^k$ and $\wh{\Xi}^{d-k}$, we have $\norm{Y_2}_{(\rho,v)} \leq 2 \epsilon \norm{Y_1}_{(\rho,v)}$ and $\norm{Z_1}_{(\rho,v)} \leq 2\epsilon \norm{Z_2}_{(\rho,v)}$. Further, 
 $$
 \phi_t(Y) = \Big( A_t(Y_1) + B_t(Y_2) \Big) + \Big( C_t(Y_1)+D_t(Y_2) \Big)
 $$
 and since $\phi_t(Y) \in \wh{\Theta}^k|_{(\rho, \phi_t(v)))}$ we have 
 $$
\norm{C_t(Y_1)+D_t(Y_2) }_{(\rho, \phi_t(v))} \leq 2\epsilon \norm{ A_t(Y_1) + B_t(Y_2)}_{(\rho, \phi_t(v))}.
$$
Thus by Lemma~\ref{lem:bds on operators}
\begin{align*}
\norm{\phi_t(Y)}_{(\rho, \phi_t(v))} & \leq (1+2\epsilon) \norm{ A_t(Y_1) + B_t(Y_2)}_{(\rho, \phi_t(v))} = (1+2\epsilon) \norm{A_t\left( Y_1 + A_t^{-1} B_t(Y_2) \right)}_{(\rho, \phi_t(v))} \\
& \leq(1+2\epsilon)\norm{A_t}_{(\rho, v)} \left( 1+ 2\epsilon^2 \right) \norm{Y_1}_{(\rho, v)} \leq(1+2\epsilon)^2 \left( 1+ 2\epsilon^2 \right) \norm{A_t}_{(\rho, v)} \norm{Y}_{(\rho, v)}.
\end{align*}
Similar reasoning shows that 
\begin{align*}
\norm{\phi_t(Z)}_{(\rho, \phi_t(v))} & \geq  (1-2\epsilon)^2 \left( 1- 2\epsilon^2 \right) \frac{1}{\norm{D_t^{-1}}_{(\rho, \phi_t(v))}} \norm{Z}_{(\rho, v)}.
\end{align*}
So by Lemma~\ref{lem:bds on operators} and Equation~\eqref{eqn:epsilon is super small},
$$
\frac{\norm{\phi_t(Y)}_{(\rho, \phi_t(v))}}{\norm{\phi_t(Z)}_{(\rho, \phi_t(v))}} \leq \epsilon\frac{(1+2\epsilon)^2(1+2\epsilon^2)}{(1-2\epsilon)^2(1-2\epsilon^2)} \frac{\norm{Y}_{(\rho, v)}}{\norm{Z}_{(\rho,v)}}\leq \frac{1}{2}\frac{\norm{Y}_{(\rho, v)}}{\norm{Z}_{(\rho,v)}}.
$$
This proves the estimate in Equation~\eqref{eqn:1/2 contraction after time T}. 

We then may apply  Equation~\eqref{eqn:1/2 contraction after time T} iteratively to show that, for all $n\in\mathbb N$,
\begin{equation*}
\frac{\norm{\phi_t(Y)}_{(\rho, \phi_t(v))}}{\norm{\phi_t(Z)}_{(\rho, \phi_t(v))}} \leq \left(\frac{1}{2}\right)^n\frac{\norm{Y}_{(\rho, v)}}{\norm{Z}_{(\rho,v)}}
\end{equation*}
 if $\rho \in \Oc$, $v \in \wh{\Usf}(\Gamma)$, $t \in [nT,(n+1)T]$, $Y \in \wh\Theta^k|_{(\rho,v)}$ and $Z \in \wh\Xi^{d-k}|_{(\rho,v)}$ is non-zero.

Finally
\begin{equation*}
\frac{\norm{\phi_t(Y)}_{(\rho, \phi_t(v))}}{\norm{\phi_t(Z)}_{(\rho, \phi_t(v))}} \leq C_0e^{-c_0 t} \frac{\norm{Y}_{(\rho, v)}}{\norm{Z}_{(\rho,v)}}
\end{equation*}
for all  $\rho \in \Oc$, $v \in \wh{\Usf}(\Gamma)$, $t \ge0$, $Y \in \wh\Theta^k|_{(\rho,v)}$ and non-zero $Z \in \wh\Xi^{d-k}|_{(\rho,v)}$,
where  $c_0:=\frac{\log 2}{T}$ and 
$$C_0:=2 \sup\left\{e^{c_0 t} \frac{\norm{\phi_t(Y)}_{(\rho,\phi_t(v))}}{\norm{Y}_{(\rho,v)}}:\rho\in\Oc, v\in\Usf(\Gamma),\ Y\in \wh \Xi^{d-k}|_{(\rho,v)}-0,\ t\in[0,T]\right\}$$
(notice that $C_0$ is finite by Lemma~\ref{obs:distortion_of_norms_holder}).

\end{proof}

\subsection*{Existence of limit maps.}
Next, we use the $\phi_t$-invariant splitting $\wh E(\Oc)=\wh\Theta^k\oplus\wh\Xi^{d-k}$ to define limit maps 
\[\xi=(\xi^k,\xi^{d-k}):\Lambda(\Gamma)\to\Gr_k(\Kb^d)\times\Gr_{d-k}(\Kb^d).\] 
 Lift this splitting of $\wh E(\Oc)$ to a splitting
$$E(\Oc)=\Theta^k\oplus\Xi^{d-k}.$$
The flow $\phi_t$ on $\wh E(\Oc)$ lifts to a flow, also denoted $\phi_t$, on $E(\Oc)$, under which this splitting is invariant. Then the bundle $\Hom(\wh\Theta^k,\wh\Xi^{d-k})$ lifts to the bundle $\Hom(\Theta^k,\Xi^{d-k})$. Finally, we use $\norm{\cdot}$ to denote the lifted norms on $E(\Oc)$ and $\Hom(\Theta^k,\Xi^{d-k})$.

By Proposition \ref{prop: contracting Hom bundle}, there exists $C_0, c_0 > 0$ such that
\begin{equation}
\label{eqn: consequence of cor contracting Hom bundle}
\norm{f}_{(\rho,\phi_t(v))} \leq C_0 e^{-c_0 t} \norm{f}_{(\rho,v)}
\end{equation}
for all $f \in \Hom(\Xi^{d-k}, \Theta^{k})|_{(\rho,v)}=\Hom(\Xi^{d-k}, \Theta^{k})|_{(\rho,\phi_t(v))}$ and $t \geq 0$.

Let
$$\sigma=(\sigma^k,\sigma^{d-k}):\Oc\times\Usf(\Gamma)\to\mathrm{Gr}_k(\mathbb K)\times\mathrm{Gr}_{d-k}(\mathbb K^d)$$
be the map so that  $\sigma^k(\rho,v)=\Theta^k|_{(\rho,v)}$ and $\sigma^{d-k}(\rho,v)=\Xi^{d-k}|_{(\rho,v)}$. Since $\sigma$ is $\phi_t$-invariant, $\sigma^k(\rho,v)$ and $\sigma^{d-k}(\rho,v)$ depend only on $\rho$, $v^+$ and $v^-$.
We now check that $\sigma^{d-k}$ depends only on $\rho$ and $v^-$.
Let $\gamma\in\Gamma$ be a hyperbolic element, let $v_\gamma\in\Usf(\Gamma)$ be a vector
so that $v_\gamma^+=\gamma^+$ and $v_\gamma^-=\gamma^-$. Let $\ell(\gamma)$ be the translation distance of $\gamma$ on $\mathbb H^2$.
Then $\phi_{n\ell(\gamma)}(v_\gamma)=\gamma^n(v_\gamma)$ for all $n$. Since $\sigma^{d-k}$ is equivariant and $\phi_t$-invariant, it follows that 
\[\sigma^{d-k}(\rho,v_\gamma)=\sigma^{d-k}(\rho,\phi_{-\ell(\gamma)}(\gamma( v_\gamma)))=\rho(\gamma)(\sigma^{d-k}(\rho,v_\gamma))\]
for all $\rho$. Furthermore, if $W\subset\Kb^d$ is a $(d-k)$-dimensional subspace that is transverse to $\sigma^{k}(\rho,v_\gamma)$, we may view $W$ as the graph of an element
$f_W\in\Hom(\Xi^{d-k},\Theta^k)|_{(\rho,v_\gamma)}$. 

Equation~\eqref{eqn: consequence of cor contracting Hom bundle} implies that
\[\norm{\rho(\gamma)^{-n}(f_W)}_{(\rho,v_\gamma)}=\norm{f_W}_{(\rho,\gamma^{n}(v_\gamma))}=\norm{f_W}_{(\rho,\phi_{n\ell(\gamma)}(v_\gamma))}\to 0,\]
which implies that $\rho(\gamma)^{-n}(W)\to\sigma^{d-k}(\rho,v_\gamma)$. Thus, $\sigma^{d-k}(\rho,v_\gamma)$ is the repelling fixed point of $\rho(\gamma)$ in $\Gr_{d-k}(\Kb^d)$. 
Now if $x\in\Lambda(\Gamma)\setminus\{\gamma^-\}$, 
then there exists $v\in\Usf(\Gamma)$ so that $v^+=x$, $v^-=\gamma^-$ and
\[\lim_{t\to \infty}d_{\Usf(\Gamma)}\big(\phi_{-t}(v_\gamma), \phi_{-t}(v)\big)=0,\]
so $\gamma^{n}( \phi_{-n\ell(\gamma)}(v))\to v_\gamma$. Thus,
\begin{align*}
\sigma^{d-k}(\rho,v_\gamma)=\lim_{n\to\infty}\sigma^{d-k}(\rho,\gamma^{n}(\phi_{-n\ell(\gamma)}(v))=\lim_{n\to\infty}\rho(\gamma)^{n}\big(\sigma^{d-k}(\rho,v)\big).
\end{align*}
Since $\sigma^{d-k}(\rho,v_\gamma)$ is the repelling fixed point of $\rho(\gamma)$ in $\Gr_{d-k}(\Kb^d)$, this implies that $\sigma^{d-k}(\rho,v)=\sigma^{d-k}(\rho,v_\gamma)$.
Therefore, since $\sigma^{d-k}$ is $\phi_t$-invariant, if $v^-=\gamma^-$, then $\sigma^{d-k}(\rho,v)$ is the repelling fixed point of $\rho(\gamma)$. Since
every point in $\Lambda(\Gamma)$ is a limit of repelling fixed points of hyperbolic elements of $\Gamma$, this implies that 
$\sigma^{d-k}$ depends only on $\rho$ and $v^-$.

One may similarly show that $\sigma^{k}$ depends only on $\rho$ and $v^+$, so there exists 
$$\xi=(\xi^k,\xi^{d-k}):\Oc\times\Lambda(\Gamma)\to\mathrm{Gr}_k(\mathbb K^d)\times\mathrm{Gr}_{d-k}(\mathbb K^d)$$
so that $\sigma(\rho,v)=(\xi^k(\rho,v^+),\xi^{d-k}(\rho,v^-))$.
As such, if $\rho\in\Oc$, then $\rho$ is $P_k$-Anosov. This proves (1).

\subsection*{The limits maps are H\"older} We now prove that, perhaps after shrinking our neighborhood $\Oc$ again, that the limit maps are uniformly H\"older. It is possible to establish this using Shub's $C^r$-Section theorem \cite[Thm. 5.18]{shub-book}, however setting up bundles with the correct regularity (see \cite[Cor. 5.19]{shub-book}) and verifying the admissibility condition is somewhat involved when $\wh\Usf(\Gamma)$ is non-compact. Instead we provide a direct argument based on the proof of Lemma 4.4 in~\cite{ZZ2019}. 

We will continue to work with $E(\Oc)$ and $\Hom(\Xi^{d-k}, \Theta^{k})|_{(\rho,v)}$. For $\rho \in \Oc$, $v \in \Usf(\Gamma)$ and $x \in \Lambda(\Gamma)-\{ v^+\}$, let $f_{\rho,v,x} \in \Hom(\Xi^{d-k}, \Theta^{k})|_{(\rho,v)}$ denote the unique element with 
\begin{align*}
{\rm Graph}(f_{\rho,v,x}) = \xi^{d-k}_\rho(x).
\end{align*}
Notice that $f_{\rho,\phi_t(v),x} = f_{\rho,v,x}$ for all $t \in \Rb$. 
If $v\in \Usf(\Gamma)$, we let $v^\perp\subset\partial\mathbb H^2$ denote the endpoints of the geodesic through the basepoint of $v$ which is
orthogonal to $v$. (One can use the orientation to canonically identify them as $v^{\perp,+}$ and $v^{\perp,-}$ but this will not be needed for our purposes.)

\begin{lemma}\label{lem:uniform_bd}  
Up to taking a subneighborhood of $\Oc$, there exists $C_1 > 1$ so that if $\rho \in \Oc$, $v\in\Usf(\Gamma)$ and $x\in v^\perp\cap\Lambda(\Gamma)$, then 
\begin{align*}
\frac{1}{C_1} \leq \norm{f_{\rho,v,x}}_{(\rho,v)} \leq C_1.
\end{align*}
\end{lemma} 

In the case when $\Gamma$ is convex co-compact, the lemma is a simple consequence of equivariance and compactness, but in the general geometrically finite case the proof is somewhat involved. Delaying the proof of the lemma, we first complete the proof of part (2) of Theorem \ref{thm: stability}.

Let $\Oc'\subset\Oc$ be a subneighborhood of $\rho_0$ such that the closure of $\overline{\Oc'}$ of $\Oc'$ is a compact subset of $\Oc$. Fix a compact set $K \subset \Usf(\Gamma)$ such that 
\begin{align*}
\Lambda(\Gamma) = \{ v^+ : v \in K\}.
\end{align*}
Let $d_\infty$ denote the distance induced by a Riemannian metric on $\partial \Hb$, and let $d_G$ denote the distance induced by a Riemannian metric on $\Gr_{d-k}(\Kb^d)$. Fix $\delta > 0$ such that: if $v\in K$ and $x\in\Lambda(\Gamma)-\{v^-\}$ satisfies $d_{\infty}(x,v^-)\leq\delta$, then there exists $t \geq 0$ such that $x\in\phi_{-t}(v)^\perp$ (in particular $x \neq v^+$). 

Fix $\rho\in\Oc'$ and $v\in K$. Then $d_G$ is bilipschitz to the norm $\norm{\cdot}_{(\rho,v)}$ on any compact subset of the affine chart 
\[\Hom(\Xi^{d-k},\Theta^k)|_{(\rho,v)}\simeq \{V\in\Gr_{d-k}(\Kb^d):V\text{ is transverse to }\xi^k_\rho(v^+)\},\]
where the isomorphism identifies each $f\in \Hom(\Xi^{d-k},\Theta^k)|_{(\rho,v)}$ with its graph in $\Gr_{d-k}(\Kb^d)$. By the compactness of $\overline{\Oc'}\times K$, there exists $C_2 > 0$ such that: if $\rho \in \Oc'$, $v \in K$ and $x \in \Lambda(\Gamma)$ with $d_{\infty}(x,v^-)\leq\delta$, then 
\begin{align*}
d_G\left( \xi_\rho^{d-k}(x), \xi_\rho^{d-k}(v^-)\right) \leq C_2 \norm{f_{\rho,v,x}}_{(\rho,v)}.
\end{align*}
There also exists $C_3> 0$ such that: if $v \in K$, $x \in \Lambda(\Gamma)$, $d_{\infty}(x,v^-) \leq \delta$ and $x\in\phi_{-t}(v)^\perp$, then 
\begin{align*}
\frac{1}{C_3} e^{-t} \leq d_\infty(x,v^-) \leq C_3 e^{-t}.
\end{align*}
Finally, let 
$$
C_4 ={\rm diam}\left( \Gr_{d-k}(\Kb^d), d_G \right)
$$
and $C = \max\{ \delta^{-c_0} C_4, C_0 C_1 C_2 C_3^{c_0}\}$. 

Now suppose $x,y \in \Lambda(\Gamma)$ and $\rho\in\Oc'$. If $d_\infty(x,y) > \delta$, then 
\begin{align*}
d_G\left( \xi^{d-k}_\rho(x), \xi^{d-k}_\rho(y)\right)  \leq C_4  \leq C_4\frac{d_\infty(x,y)^{c_0}}{\delta^{c_0}}  \leq C d_\infty(x,y)^{c_0}.
\end{align*}
If $d_{\infty}(x,y) \leq \delta$, then there exists $t \geq 0$  and $v\in K$ such that $y=v^-$ and $x\in\phi_{-t}(v)^\perp$. Then 
\begin{align*}
d_G\left( \xi^{d-k}_\rho(x), \xi^{d-k}_\rho(v^-)\right)& \leq C_2 \norm{f_{\rho,v,x}}_{(\rho,v)} \leq C_0 C_2 e^{-c_0 t}  \norm{f_{\rho,\phi_{-t}(v),x}}_{(\rho,\phi_{-t}(v))} \\
& \leq C_0 C_1 C_2 C_3^{c_0} d_\infty(x,v^-)^{c_0} \leq C d_\infty(x,y)^{c_0}.
\end{align*}
Therefore, $\xi_\rho^{d-k}$ is $c_0$-H\"older. 

One may similarly prove that, perhaps after passing to another sub-neighborhood, that  if $\rho\in\Oc'$, then $\xi_\rho^k$ is $c_1$-H\"older for some $c_1>0$.
So, (2) holds with $\alpha=\min\{c_0,c_1\}$.

\subsection*{Proof of Lemma~\ref{lem:uniform_bd}}
Let $\Oc'\subset\Oc$ be a subneighborhood of $\rho_0$ such that the closure of $\overline{\Oc'}$ of $\Oc'$ is a compact subset of $\Oc$. If the lemma fails for $\Oc'$, then there exist sequences $\{\rho_m\}$ in $\Oc^\prime$,  $\{v_m\}$  in $\Usf(\Gamma)$ and $\{x_m\}$ in $\Lambda(\Gamma)$ 
so that $x_m\in v_m^\perp$ for all $m$ and 
\begin{align}
\label{eqn:ratio_blow_up}
\lim_{m \rightarrow \infty} \abs{\log \norm{f_{\rho_m,v_m,x_m}}_{(\rho_m,v_m)}}=+\infty. 
\end{align}
By passing to a subsequence we can suppose that $\rho_m \rightarrow \rho \in \Oc$. 
After passing to a further subsequence and translating by elements in $\Gamma$, either  
\begin{enumerate}
\item $v_m \rightarrow v \in \Usf(\Gamma)$ and $x_m \rightarrow x \in \Lambda(\Gamma)$, or 
\item there exists an embedded cusp neighborhood $C=\ip{\alpha}\backslash H$ such that $\{v_m\}\subset \Usf(\Gamma)_H$ and $\{ v_m\}$ projects to an escaping sequence in $\wh{\Usf}(\Gamma)$. 
\end{enumerate}
In the first case, $x\in v^\perp$ and $\lim_{m \rightarrow \infty} \norm{f_{\rho_m,v_m,x_m}}_{(\rho_m,v_m)}=\norm{f_{\rho,v,x}}_{(\rho,v)}\ne 0$. Thus we must be in the second case.

By conjugating, we may assume that $\rho_m(\alpha)=\rho_0(\alpha)$ for all $m\in\Zb^+\cup\{\infty\}$, and that the restriction of the canonical norm on $E(\Oc')$ to each $E_{\rho_m}$ is with respect to the same cusp representation $\Psi$  for $\alpha$ and $\rho_0(\alpha)$. Then there exists a $\rho_0(\alpha)_{ss}$-invariant, $\Psi$-equivariant family of norms $\norm{\cdot}^\star_{v \in T^1 \Hb^2}$ such that 
\begin{equation}
\label{eqn:star norm versus regular norm 2}
\norm{\cdot}_{(\rho,v)}=\norm{\cdot}_{v}^{\star} 
\end{equation}
 for all $\rho \in \Oc'$ and $v \in \Usf(\Gamma)_H$. Also, given a linear map $f$ between subspaces of $\Kb^d$, let $\norm{f}^{\star}_v$ denote the operator norm relative to $\norm{\cdot}_v^\star$.

Let $y \in \Lambda(\Gamma)$ be the center of $H$ and fix a hyperbolic element $\gamma\in\SL(2,\Rb)$ with attracting fixed point $y$. After translating each $v_m$ by a power of $\alpha$ and passing to a subsequence, we can find $n_m \rightarrow \infty$ such that $\gamma^{-n_m} (v_m) \rightarrow v \in T^1 \Hb^2$ and 
$\gamma^{-n_m} (x_m) \rightarrow x \in \partial \Hb^2$.

Let $\eta^j: \partial \Hb^2 \rightarrow \Gr_j(\Kb^d)$ denote the boundary maps  associated to $\Psi$ for $j=k,d-k$, see Proposition~\ref{cusp rep anosov}. Let $\hat{f} \in \Hom(\eta^{d-k}(v^-),\eta^k(v^+))$ denote the unique element with 
\begin{align*}
{\rm Graph}(\hat{f}) = \eta^{d-k}(x).
\end{align*}
Notice that  $\hat{f} \neq 0$, since $\eta^k$ and $\eta^{d-k}$ are transverse and $x \in v^{\perp}$. By Proposition \ref{thm:asymptotics_of_limit_curve},
\begin{align*}
\Psi(\gamma^{-n_m})\Big( \xi_{\rho_m}^j(v_m^\pm)\Big) = \Big(\Psi(\gamma^{-n_m}) \circ \xi_{\rho_m}^j \circ \gamma^{n_m}\Big)\big( \gamma^{-n_m}(v_m^\pm)\big) \rightarrow \eta^j(v^\pm).
\end{align*}
for $j=k,d-k$. Similarly, $\Psi(\gamma^{-n_m})\big(\xi_{\rho_m}^j(x)\big) \rightarrow \eta^j(x)$, 
 for $j=k,d-k$, so
 $$\Psi(\gamma^{-n_m}) \circ f_{\rho_m,v_m,x_m} \circ \Psi(\gamma^{n_m}) \rightarrow \hat{f}.$$

Notice, that if $X \in \Xi^{d-k}|_{(\rho_m,v_m)}$, then Equation~\eqref{eqn:star norm versus regular norm 2} implies that
\begin{align*}
\norm{f_{\rho_m,v_m,x_m}(X)}_{(\rho_m,v_m)} = \norm{ \Big(\Psi(\gamma^{-n_m}) \circ f_{\rho_m,v_m,x_m} \circ \Psi(\gamma^{n_m})\Big)\big(\Psi(\gamma^{-n_m})(X)\big)}_{\gamma^{-n_m}(v_m)}^\star 
\end{align*}
and $\norm{\Psi(\gamma^{-n_m})(X)}_{\gamma^{-n_m}(v_m)}^\star =\norm{X}_{(\rho_m,v_m)}$. Thus
\begin{align*}
\lim_{m \rightarrow \infty} \norm{ f_{\rho_m,v_m,x_m}}_{(\rho_m,v_m)} = \lim_{m \rightarrow \infty} \norm{\Psi(\gamma^{-n_m}) \circ f_{\rho_m,v_m,x_m} \circ \Psi(\gamma^{n_m})}_{\gamma^{-n_m}(v_m)}^\star = \norm{\hat{f}}_{v}^\star \neq 0
\end{align*}
and we have a contradiction.  This completes the proof of Lemma \ref{lem:uniform_bd} and hence the proof of (2).

\subsection*{The limits maps vary analytically} 
It remains to prove the analytic variation of the limit maps. First suppose that $\Kb = \Rb$. The general strategy is to complexify and then exploit the fact that locally uniform limits of complex analytic functions are complex analytic. 

Suppose that $h:M \rightarrow \Hom_{\rm tp}(\rho_0)$ is a real analytic map and every representation in $h(M)$ is $P_k$-Anosov. 

If $\rho:\Gamma\to\mathsf{SL}(d,\mathbb R)$ is $P_k$-Anosov, we may compose
with the  inclusion map $\iota_2:\mathsf{SL}(d,\mathbb R)\to\mathsf{SL}(d,\mathbb C)$
to obtain a $P_k$-Anosov representation $\rho^{\mathbb C}=\iota_2\circ\rho:\Gamma \to\mathsf{SL}(d,\Cb)$. Fix generators $g_1,\dots g_N$ of $\Gamma$ and view $\Hom(\Gamma, \mathsf{SL}(d,\mathbb C))$ as a subset of $\mathsf{SL}(d,\mathbb C)^N$. We can then view $h$ as a map $h :M\rightarrow\mathsf{SL}(d,\mathbb C)^N$. We can also realize $M$ as a totally real submanifold of a complex manifold $M^{\Cb}$ and then extend $h$ to a complex analytic map $h : M^{\Cb} \rightarrow\mathsf{SL}(d,\mathbb C)^N$. Notice that $h(M)$ and $h(M^{\Cb})$ have the same Zariski closure in $\mathsf{SL}(d,\mathbb C)^N$. 

We claim, after possibly shrinking $M^{\Cb}$, that $h(M^{\Cb})\subset\Hom_{\rm tp}(\rho_0^{\Cb})$. For any $\alpha \in \SL(d,\Cb)$, the set 
$$
\{ g \in \SL(d,\Cb) : g \text{ is conjugate to } \alpha\}
$$
is locally closed (i.e. open in its closure) in the Zariski topology, see for instance~\cite[Theorem 3.6]{byrnes-gauger}. This implies that $\Hom_{\rm tp}(\rho_0^{\Cb})$ is itself a locally closed set in the Zariski topology on $ \mathsf{SL}(d,\mathbb C)^N$. Then since $h(M) \subset  \Hom_{\rm tp}(\rho_0^{\Cb})$, by shrinking $M^{\Cb}$ we may assume that $h(M^{\Cb}) \subset \Hom_{\rm tp}(\rho_0^{\Cb})$. 

Since every representation in $h(M)$ is $P_k$-Anosov, by shrinking $M^{\Cb}$ again if necessary, we may assume that every representation in $h(M^{\Cb})$ is also $P_k$-Anosov. Thus, if we can prove that for any $x\in \Lambda(\Gamma)$, the map $u\mapsto \xi^k_{h(u)}(x)$ from $M^{\Cb}$ to $\Gr_k(\Cb^d)$ is complex analytic, then its restriction to $M$ is real analytic.

If $\gamma$ is a hyperbolic element and $\rho\in h(M^{\Cb})$, then $\rho(\gamma)$ is $P_k$-proximal and $\xi^k_\rho(\gamma^+)$ is the attracting $k$-plane of $\rho(\gamma)$. Then it follows from standard results in the perturbation theory of linear operators,
see, for example, \cite[Chapter 6]{kato}, that the function from $M^{\Cb}$ to $\mathrm{Gr}_k(\mathbb C^d)$ given by $u \mapsto \xi^k_{h(u)}(\gamma^+)$ is complex analytic. If $x\in\Lambda(\Gamma)$, then there exists a sequence $\{\gamma_n\}$ of hyperbolic elements of $\Gamma$,
so that $\gamma_n^+\to x$. Then, since the map 
$$
(u,y) \in M^{\Cb}\times\Lambda(\Gamma)\mapsto \xi^k_{h(u)}(y) \in \mathrm{Gr}_k(\mathbb C^d)
$$ 
is continuous, the function $u \mapsto \xi^k_{h(u)}(x)$ is a locally uniform limit of complex analytic functions, hence complex analytic.
This completes the proof of (3) in the case when $\Kb=\Rb$. 

The  case when $\Kb = \Cb$ case follows by simply repeating the argument in the previous  paragraph. 
\end{proof}


\section{Positive representations in the sense of Fock-Goncharov}

In Fock and Goncharov's work \cite{fock-goncharov}, they define positive representations in the following way. Suppose that $\Gamma_0\subset\PSL(2,\Rb)$ is  a 
discrete group and $\Gamma_0\backslash\Hb^2$ is a non-compact, finite area hyperbolic surface. Recall that $\gamma\in\Gamma_0$ is {\em peripheral} if it is represented by a curve which may be freely homotoped
off of every compact subset of $S=\Gamma \backslash \Hb^2$. Let $\Lambda_p(\Gamma_0)$ be the set of fixed points 
of the peripheral elements in $\Gamma_0$ (in this case, all of which are parabolic). Notice that $\Lambda_p(\Gamma_0)$ inherits two natural cyclic orders as
a subset of $\partial\mathbb H^2$. A representation $\rho:\Gamma_0\to\SL(d,\Rb)$ is \emph{positive} if there is a positive, $\rho$-equivariant map 
$\zeta:\Lambda_p(\Gamma_0)\to\Fc(\Rb^d)$. Notice that, with Fock and Goncharov's definition, every Hitchin representation of a convex cocompact, but not
cocompact, Fuchsian group $\Gamma_0$ is also a positive representation of  a lattice $\Gamma$. (If $\Gamma\backslash\mathbb H^2$ is homeomorphic to
the interior of a compact surface $S$, then $\Gamma$ is
a finite area uniformization of  the interior of $S$.) In this case,  every peripheral element is mapped to a loxodromic element and 
there are many different  positive $\rho$-equivariant maps from $\Lambda_p(\Gamma)$ to $\Fc(\Rb^d)$, corresponding to a choice of fixed point for the (unique) fixed
point of each peripheral element of $\Gamma$.

Motivated by their definition, we define the notion of a positive type preserving representation in the following way. Let $\Gamma\subset\PSL(2,\Rb)$ be a geometrically finite group. 
Then $\gamma\in\Gamma$ is peripheral if it is either unipotent or it is a hyperbolic
element whose fixed points both lie in the boundary of $\Lambda(\Gamma)$. 

\begin{definition}
A representation $\rho:\Gamma\to\SL(d,\Rb)$ of a geometrically finite Fuchsian group is \emph{positive type preserving} if
\begin{itemize}
\item $\rho(\gamma)$ is weakly unipotent for every parabolic $\gamma\in\Gamma$, and
\item there is a positive, $\rho$-equivariant map $\zeta:\Lambda_p(\Gamma)\to\Fc(\Rb^d)$. 
\end{itemize}
\end{definition}

Observe that if $\rho:\Gamma\to\SL(d,\Rb)$ is a positive type preserving representation, then $\rho \circ f_*:\Gamma_0\to\SL(d,\Rb)$ is a positive representation for all homeomorphisms $f:\Gamma_0\backslash\Hb^2\to\Gamma\backslash\Hb^2$.

It is clear that every Hitchin representation from $\Gamma$ to $\SL(d,\Rb)$ is a positive type preserving representation. We show that the converse is also true. 

\begin{theorem}\label{thm: main1}
If $\Gamma \subset \PSL(2,\Rb)$ is geometrically finite group, then 
every positive type preserving representation $\rho:\Gamma\to\mathsf{SL}(d,\mathbb R)$ is Hitchin. 
\end{theorem}

Let $\rho$ be a positive, type preserving representation, and let $\zeta:\Lambda_p(\Gamma)\to\Fc(\Rb^d)$ 
denote a positive, $\rho$-equivariant map.
We make use of the following well-known fact which is implicit in Fock-Goncharov \cite{fock-goncharov}, see
Kim-Tan-Zhang \cite[Observation 3.18]{KTZ}  for details. It may be viewed as a generalization of the fact that every bounded monotone
sequence  in $\mathbb R$ is convergent.

\begin{proposition}\label{prop: limit}
Let $\{F_n\}$ be a sequence of flags in $\Fc(\Rb^d)$ and  $H_1,H_2\in\Fc(\Rb^d)$ such that $(F_1,\dots,F_n,H_1,H_2)$ is positive for all $n$. 
Then the sequence $\{F_n\}$ converges to a flag $F_\infty\in\Fc(\Rb^d)$, and $(F_1,\dots,F_n,F_\infty,H_2)$ is positive for all $n$.
\end{proposition}

We fix for the remainder of the section,  one of the two natural cyclic orders on $\partial \Hb^2$.  With this convention, it is
natural to define one-sided convergence of sequences.

\begin{definition}
A sequence $\{x_n\}$ in $\partial \Hb^2$ converges to $x\in\partial \Hb^2$ \emph{in the positive direction} (respectively, \emph{in the negative direction}) if $\{x_n\}$ converges to $x$, and there exists $N > 0$ such that 
\begin{align*}
x_N < x_{N+1} < \dots < x\,\,\, (\text{respectively, }x_N > x_{N+1} > \dots > x)
\end{align*}
For short, we write $x_n \nearrow x$ (respectively, $x_n \searrow x$) if $\{x_n\}$ converges to $x$ in the positive direction (respectively, in the negative direction). 
If $\{x_n\}$ converges to $x$ in either the positive direction or the negative direction,  we say that $\{x_n\}$ converges \emph{monotonically} to $x$.
\end{definition}

\begin{definition}
Let $x\in\Lambda(\Gamma)$. 
\begin{itemize}
\item If there are sequences in $\Lambda_p(\Gamma)$ that converge to $x$ from the positive direction (respectively, from the negative direction), set 
\begin{align*}
\xi_+(x) := \lim_{y \in \Lambda_p(\Gamma), y \nearrow x} \zeta(y)\,\,\, \left(\text{respectively, }\,\xi_-(x) := \lim_{y \in \Lambda_p(\Gamma), y \searrow x} \zeta(y)\right).
\end{align*}
\item If there are no sequences in $\Lambda_p(\Gamma)$ that converge to $x$ in the positive direction (respectively, in the negative direction), then there necessarily are sequences in $\Lambda_p(\Gamma)$ that converge to $x$ in the negative direction (respectively, in the positive direction). Thus, we may set
\[\xi_+(x) := \xi_-(x)\,\,\,(\text{respectively, }\,\xi_-(x) := \xi_+(x)).\]
\end{itemize}
Since $\Lambda_p(\Gamma)$ is dense in $\Lambda(\Gamma)$, Proposition \ref{prop: limit} implies that these limit maps are well-defined.
We refer to the maps $\xi_+:\Lambda(\Gamma)\to\Fc(\Rb^d)$ (respectively, $\xi_-: \Lambda(\Gamma)\to\Fc(\Rb^d)$) as the \emph{plus limit map} 
(respectively, \emph{minus limit map}).
\end{definition}

Since $\zeta$ is $\rho$-equivariant, both $\xi_+$ and $\xi_-$ are $\rho$-equivariant. We next check that they satisfy the
following positivity property, which implies in particular, that they are both positive. 

\begin{proposition}\label{prop: positive plus minus}
Let $x_1<x_2<\dots<x_k$ be points in $\Lambda(\Gamma)$, and let $s_1,\dots,s_k\in\{+,-\}$. Then
\[(\xi_{s_1}(x_1),\xi_{s_2}(x_2),\dots,\xi_{s_k}(x_k))\]
is a positive tuple of flags. 
\end{proposition}

\begin{proof} 
For each $i=1,\dots,k$, let $q_i,q_i'\in\Lambda_p(\Gamma)$ be  points satisfying the following conditions:
\begin{itemize}
\item $q_i=x_i$ if there are no sequences in $\Lambda_p(\Gamma)$ that converge to $x_i$ in the positive direction,
\item $q_i'=x_i$ if there are no sequences in $\Lambda_p(\Gamma)$ that converge to $x_i$ in the negative direction,
\item $x_i<q_i'<q_{i+1}$ if there is a sequence in $\Lambda_p(\Gamma)$ that converges to $x_i$ in the negative direction, and
\item $q_{i-1}'<q_i<x_i$ if there is a sequence in $\Lambda_p(\Gamma)$ that converges to $x_i$ in the positive direction.
\end{itemize}
Here, arithmetic in the subscripts is done modulo $k$. 

By the definition of $\xi_\pm$, for each $i=1,\dots,k$, there is a sequence $(y_{i,n})_{n\geq 1}$ in $\Lambda_p(\Gamma)$ that converges monotonically to $x_i$, and satisfies $\lim_{n\to\infty}\zeta(y_{i,n})=\xi_{s_i}(x_i)$. By passing to the tail of the sequences $(y_{i,n})_{n\geq 1}$, we may assume that for all $n$,
\[q_1<y_{1,n}<q_1'<q_2<y_{2,n}<q_2'<\dots<q_k<y_{k,n}<q_k',\]
which implies that 
\[(\zeta(q_1),\zeta(y_{1,n}),\zeta(q_1'),\zeta(q_2),\zeta(y_{2,n}),\zeta(q_2'),\dots,\zeta(q_k),\zeta(y_{k,n}),\zeta(q_k'))\]
is a positive tuple of flags in $\Fc(\Rb^d)$. Our result then follows from repeatedly applying Proposition \ref{prop: limit}.

\end{proof}

We next prove that $\xi_\pm$ satisfy an analogue of the strongly dynamics-preserving property of  Anosov limit maps. In particular, if $\xi_+=\xi_-$,
then it implies that $\xi_\pm$ is strongly dynamics-preserving. In fact, the proof mimics the proof of Theorem ~\ref{thm:hitchin_are_str_dyn_preserving}.

\begin{proposition}\label{prop:weak_dynamics} 
Suppose $\{\gamma_n\}$ is a sequence in $\Gamma$ with $\gamma_n \rightarrow x\in\Lambda(\Gamma)$ and $\gamma_n^{-1} \rightarrow y\in\Lambda(\Gamma)$. 
Then after passing to a subsequence, there exists $s_1,s_2 \in \{+,-\}$ such that
\begin{align*}
\lim_{n\to\infty}\rho(\gamma_n)(F) =\xi_{s_1}(x)
\end{align*}
for all $F \in \Fc(\Rb^d)$ transverse to $\xi_{s_2}(y)$. 
\end{proposition}

\begin{proof}
We first suppose that $x \neq y$. By passing to the tail of the sequence, we may assume that each $\gamma_n$ is hyperbolic.
Then $\gamma_n^+\rightarrow x$, $\gamma_n^- \rightarrow y$, and $\gamma_n(z) \rightarrow x$ for all $z \in \Lambda(\Gamma) - \{y\}$.

Since $\Lambda_p(\Gamma)$ is infinite, there are points $a,b \in \Lambda_p(\Gamma)-\{x,y\}$ such that 
\begin{itemize}
\item either $x<a<b<y$ or $y<b<a<x$, and
\item up to taking subsequences, the sequences $\{\gamma_n( a)\}$ and $\{\gamma_n(b)\}$ both converge monotonically to $x$, and from the same direction. 
\end{itemize}
Therefore,
\begin{align*}
 \lim_{n \rightarrow \infty}\zeta(\gamma_n(a)) = \xi_{s_1}(x) =  \lim_{n \rightarrow \infty} \zeta(\gamma_n(b)) 
\end{align*}
for some $s_1 \in \{+,-\}$. For each $n\geq 1$, choose a point $c_n \in \Lambda_p(\Gamma)\setminus\{\gamma_n^-\}$ such that $c_n\to y$ and 
$\gamma_n( c_n)\to y$. Passing to a further subsequence, we can assume that $\{\gamma_n(c_n)\}$ converges monotonically to $y$. Thus,
\[\lim_{n\to\infty}\zeta(\gamma_n(c_n))=\xi_{s_3}(y)\] 
for some $s_3 \in \{+,-\}$. 

Now, consider the open sets (defined in Definition~\ref{defn:open_sets})
 \begin{align*}
 \Oc_n:=\Oc\big(\zeta(a), \zeta(b), \zeta(c_n)\big)
 \end{align*}
for all $n$. Since $c_n\to y$, and either $a<b<y$ or $b<a<y$, there exists $N>0$ such that either $a<b<c_n$ for all $n\geq N$, or $b<a<c_n$ for all $n\geq N$. 

We may then argue, exactly as in the proof of Theorem~\ref{thm:hitchin_are_str_dyn_preserving}, that there exists an open set $\Oc\subset\Fc(\Rb^d)$ so that 
 \begin{align*}
 \lim_{n \rightarrow \infty} \rho(\gamma_n) ( F) = \xi_{s_1}(x)
 \end{align*}
 for all $F \in \Oc$. Repeating the same argument with $\gamma_n^{-1}$, we see that there exists $s_2 \in \{+,-\}$ and an open set $\Oc' \subset \Fc(\Rb^d)$ where
  \begin{align*}
 \lim_{n \rightarrow \infty} \rho(\gamma_n^{-1})(F )= \xi_{s_2}(y)
 \end{align*}
 for all $F \in \Oc'$. Hence, we may apply Lemma \ref{lem: basic singular value} to deduce the proposition when $x\neq y$.

If $x = y$,  pick $\eta \in \Gamma$ such that $z:=\eta^{-1} (x) \neq x$. Then $\gamma_n\eta\rightarrow x$, $(\gamma_n\eta)^{-1} \rightarrow z \neq x$.
By the first part, there exists $s_1,s_2 \in \{+,-\}$ such that $\rho(\gamma_n\eta)( F) \rightarrow \xi_{s_1}(x)$ for all $F \in \Fc(\Rb^d)$ transverse to $\xi_{s_2}(z)$. 
Equivalently, $\rho(\gamma_n)(F) \rightarrow \xi_{s_1}(x)$ for all $F \in \Fc(\Rb^d)$ transverse to $\xi_{s_2}(x)$. 
\end{proof}

Given $x\in\Lambda(\Gamma)-\Lambda_p(\Gamma)$, it will be useful to construct a sequence in $\Gamma$ and a pair of points in $\Lambda(\Gamma)$,
so that the orbits of the points under the sequence of elements of $\Gamma$ approach $x$ from opposite sides.

\begin{lemma} \label{lem: conical two sides}
If $x \in \Lambda(\Gamma)-\Lambda_p(\Gamma)$, then there exists a sequence $\{\gamma_n\}$ in $\Gamma$, a point $y\in\Lambda(\Gamma)$, 
and points $a,b \in \Lambda(\Gamma) - \{y\}$ so that:
\begin{enumerate}
\item $\gamma_n \rightarrow x$, 
\item $\gamma_n^{-1} \rightarrow y$, 
\item $\gamma_n(a)\nearrow x$ and $\gamma_n(b)\searrow x$. 
\end{enumerate}
\end{lemma}

\begin{proof} 
Let  $r: \Rb \rightarrow \Hb^2$ be a geodesic joining a point $w\in\Lambda(\Gamma)-\{x\}$ to $x$. Since $x$ is a conical limit point,
there exists  sequences $\{t_n\}\subset \mathbb R$ and $\{\gamma_n\}\subset\Gamma$ so that $t_n\to\infty$ and $\{d_{\Hb^2}(r(t_n),\gamma_n(r(0))\}$ is bounded.
Therefore, after passing to a subsequence, $v_n=(\gamma_n^{-1}\circ r)'(t_n)$ converges to a vector $v\in\Usf(\Gamma)$. Notice that $\gamma_n(v_n^+)=x$ and $\gamma_n(v_n^-)=w$ for all $n$.

If there is a point $a\in\Lambda(\Gamma)$ such that $v^-<a<v^+$, then $v_n^-<a<v_n^+$ for sufficiently large $n$. 
It follows that $w<\gamma_n( a)<x$ for sufficiently large $n$. Since $\gamma_n( a)\to x$, it follows that $\gamma_n( a)\nearrow x$. 
On the other hand, if there are no points $a\in\Lambda(\Gamma)$ such that $v^-<a<v^+$, then $v^+$ and $v^-$ are the fixed points of a hyperbolic, peripheral element. 
In particular, $\gamma_n^{-1}(x) \neq v^+$ for any $n$. Since $\gamma_n^{-1}( x) \rightarrow v^+$, this implies that $\gamma_n^{-1}(x)\searrow v^+$, so $\gamma_n(v^+)\nearrow x$. 
In either case, there is some $a\in\Lambda(\Gamma)\setminus\{v^-\}$ such that $\gamma_n(a)\nearrow x$. 

Similarly, we may find a point $b\in\Lambda(\Gamma)-\{v^-\}$ so that $\gamma_n(b)\searrow x$. 
\end{proof}

Lemma \ref{lem: conical two sides} allows us to complete the proof of Theorem \ref{thm: main1}.

\begin{proof}[Proof of Theorem \ref{thm: main1}] We first notice that it suffices to prove that
$\xi_+=\xi_-$ and $\xi_+|_{\Lambda_p(\Gamma)}=\zeta$.
Indeed, if we can do so, then we may set $\xi:\Lambda(\Gamma)\to\Fc(\Rb^d)$ to be the map given by $\xi(x):=\xi_\pm(x)$. 
Proposition \ref{prop: positive plus minus} implies that $\xi$ is positive, while Proposition \ref{prop: limit} implies that $\xi$ is continuous.

The proof proceeds in three cases. In the first case, we assume that $x$ is not the fixed point of a peripheral element.  If $x$ is the fixed point of a peripheral element, 
then $x$ is the fixed point of either a parabolic element or a peripheral hyperbolic element.

\medskip\noindent
{\bf Case 1:  $x \in \Lambda(\Gamma)-\Lambda_p(\Gamma)$.} Let $a,b,y \in \Lambda(\Gamma)$ and $\{\gamma_n\}$ be as in Lemma \ref{lem: conical two sides}. 
By Proposition~\ref{prop:weak_dynamics}, we can pass to a subsequence so that there exists $s_1,s_2 \in \{+,-\}$ so that
\begin{align*}
\rho(\gamma_n)( F) \rightarrow \xi_{s_1}(x)
\end{align*}
for all $F \in \Fc(\Rb^d)$ transverse to $\xi_{s_2}(y)$. By Proposition~\ref{prop: positive plus minus}, the flags $\xi_+(a)$ and $\xi_+(b)$ are both transverse to $\xi_{s_2}(y)$. 
So 
\begin{align}\label{eqn: +-}
\xi_{s_1}(x)=\lim_{n \rightarrow \infty} \rho(\gamma_n)(\xi_+(a))=\lim_{n \rightarrow \infty} \xi_+(\gamma_n(a))=\xi_+(x)
\end{align}
and 
\begin{align}\label{eqn: +-2}
\xi_{s_1}(x)=\lim_{n \rightarrow \infty} \rho(\gamma_n)(\xi_+(b))=\lim_{n \rightarrow \infty} \xi_+(\gamma_n(b))=\xi_-(x).
\end{align}
The second equality in (\ref{eqn: +-}) and (\ref{eqn: +-2}) holds because $\gamma_n\in\Gamma$ for all $n$, and the last equality in (\ref{eqn: +-}) and (\ref{eqn: +-2}) is a consequence of Proposition \ref{prop: limit}.  Thus $\xi_+(x) = \xi_-(x)$. 

\medskip\noindent
{\bf Case 2: $x \in \Lambda_p(\Gamma)$ is the fixed point of a parabolic element $\alpha \in \Gamma$.}  
As in the proof of Theorem \ref{thm: Hitchin}, 
Lemma \ref{lem: peripheral parabolic} implies that $\rho(\alpha)$ has a unique fixed flag in $\Fc(\Rb^d)$. It follows that $\xi_+(x)=\xi_-(x)=\zeta(x)$ is this unique fixed flag. 

\medskip\noindent
{\bf Case 3: $x \in \Lambda_p(\Gamma)$ is the fixed point of a hyperbolic peripheral element $\gamma \in \Gamma$.} By replacing $\gamma$ with $\gamma^{-1}$ we can assume that $x=\gamma^+$. Then either $\gamma^+<a<\gamma^-$ for all $a\in\Lambda_p(\Gamma)-\{\gamma^+,\gamma^-\}$ or $\gamma^-<a<\gamma^+$ for all $a\in\Lambda_p(\Gamma)-\{\gamma^+,\gamma^-\}$. Also, by definition, $\xi_+(x) = \xi_-(x)$. 

We now show $\xi_+(x)=\zeta(x)$. Since $\gamma^n\to \gamma^+$ and $\gamma^{-n}\to \gamma^-$, 
Proposition~\ref{prop:weak_dynamics} implies that there is an increasing sequence $\{m_n\}$ of integers so that 
\begin{align*}
\rho(\gamma^{m_n})(F) \rightarrow \xi_+(x)
\end{align*}
for all $F \in \Fc(\Rb^d)$ transverse to $\xi_+(\gamma^-)$. It follows that $\rho(\gamma)$ is loxodromic, and that $\xi_+(\gamma^+)$  
and $\xi_+(\gamma^-)$ are respectively the attracting and repelling fixed flag of $\rho(\gamma)$. 

To finish the proof, it is sufficient to show that $\zeta(\gamma^+)$ is also the attracting fixed flag of $\rho(\gamma)$. 
Let $\{x_n\}$ be a sequence in $\Lambda_p(\Gamma)$ that converges monotonically to $\gamma^-$. Since $\zeta$ is positive, the tuple
\[(\zeta(\gamma^+),\zeta(x_1),\dots,\zeta(x_n),\zeta(\gamma^-))\]
is positive for all $n$. Since $\xi_+(\gamma^-)=\lim_{n\to\infty}\zeta(x_n)$, Proposition \ref{prop: limit} implies that 
\[(\zeta(\gamma^+),\zeta(x_1),\dots,\zeta(x_n),\xi_+(\gamma^-))\]
is positive for all $n$. In particular, $\zeta(\gamma^+)$ and $\xi_+(\gamma^-)$ are transverse. Since $\zeta(\gamma^+)$ is fixed by $\rho(\gamma)$ and 
$\xi_+(\gamma^-)$ is the repelling fixed point of $\rho(\gamma)$, it follows that $\zeta(\gamma^+)$ is the attracting fixed flag of $\rho(\gamma)$. 
\end{proof}

\appendix

\section{Constructing cusp representations}\label{appendix: cusp repn}

In this appendix, we prove Proposition~\ref{prop:building_repn}. Suppose $g \in \SL(d,\Kb)$ is weakly unipotent and let
$$\u2 = \begin{pmatrix} 1 & 1\\ 0 & 1\\ \end{pmatrix}.$$ 
We will construct a representation $\Psi: \SL(2,\Rb) \rightarrow \SL(d,\Kb)$ where $\Psi\left(\u2\right) = g_u$ and $g_{ss}$ commutes with the elements of $\Psi(\SL(2,\Rb))$.

First suppose that $\Kb=\Cb$. Then using the Jordan normal form there exist $p \in \SL(d,\Cb)$, integers $d_1\ge d_2\ge \cdots\ge d_m>0$, and complex numbers $\lambda_1, \dots, \lambda_m \in \mathbb{S}^1$ so
that  $d_1+\cdots+d_m=d$,
$$pg_up^{-1}=\oplus_{j=1}^m \tau_{d_j}\left(\u2\right) \quad \text{and} \quad p g_{ss} p^{-1} = \oplus_{j=1}^{m} \lambda_j \id_{d_j}.$$
It follows that $\Psi$ defined by
\[\Psi(\beta)= p^{-1}\left(\oplus_{i=1}^m \tau_{d_i}\left(\beta\right)\right)p\]
for all $\beta\in\SL(2,\Rb)$ has the desired properties. 

Next suppose that $\Kb=\Rb$. Given two matrices $A, B$ we will let $A \otimes B$ denote the \emph{Kronecker product}, that is if $A=[a_{ij}]$ is an $m$-by-$n$ matrix, then
$$
A \otimes B = \begin{pmatrix} a_{11} B & \cdots & a_{1n} B \\ \vdots & & \vdots \\ a_{m1} B & \cdots & a_{mn} B \end{pmatrix}.
$$
Also given $\theta \in \Rb$ let 
$$
M(\theta) = \begin{pmatrix} \cos \theta & \sin \theta \\ -\sin \theta & \cos \theta \end{pmatrix}.
$$
Then using the real Jordan normal form there exist $p \in \SL(d,\Rb)$, integers $d_1\ge d_2\ge \cdots\ge d_{m+n} >0$ and numbers $\lambda_1, \dots, \lambda_m \in \{-1,1\}$, $\theta_{m+1}, \dots, \theta_{m+n} \in \Rb$ so 
that  $d_1+\cdots+d_m+2d_{m+1} + \dots + 2d_{m+n}=d$,
$$pg_up^{-1}=\left[\oplus_{j=1}^m \tau_{d_j}\left(\u2\right)\right] \oplus \left[ \oplus_{j=m+1}^{m+n} \tau_{d_j}\left(\u2\right)\otimes \id_2 \right]$$
and 
$$
p g_{ss} p^{-1} = \left[\oplus_{j=1}^{m} \lambda_j \id_{d_j} \right] \oplus \left[ \oplus_{j=m+1}^{m+n} \id_{d_j} \otimes M(\theta_j) \right].
$$
By the multiplicative property of the Kronecker product
$$
\left[  \tau_{d_j}\left(\beta\right)\otimes \id_2 \right]\left[  \id_{d_j} \otimes M(\theta_j) \right] =\left[  \id_{d_j} \otimes M(\theta_j) \right] \left[  \tau_{d_j}\left(\beta\right)\otimes \id_2 \right] =  \tau_{d_j}\left(\beta\right) \otimes M(\theta_j)
$$
for all $\beta \in \SL(2,\Rb)$ and so it follows that $\Psi$ defined by
$$
\Psi(\beta)= p^{-1}\left( \left[\oplus_{j=1}^m \tau_{d_j}\left(\beta\right)\right] \oplus \left[ \oplus_{j=m+1}^{m+n} \tau_{d_j}\left(\beta\right)\otimes \id_2 \right]\right)p
$$
has the desired properties.

\section{Anosov representations into semisimple Lie groups}\label{app: general}

In this appendix, we develop a more general theory of Anosov representations of geometrically finite Fuchsian groups into a semisimple Lie group $\mathsf{G}$ with
respect to a parabolic subgroup $P^+$.
We extend results of Guichard-Wienhard \cite[Prop. 4.3]{guichard-wienhard} (see also Gu\'eritaud-Guichard-Wienhard \cite[Section 3]{GGKW})
to show that there exists an irreducible representation $\psi:\mathsf{G}\to\mathsf{SL}(d,\mathbb R)$
so that a representation into $\mathsf{G}$ is Anosov with respect to $P^+$ if and only if its composition with $\psi$ is $P_1$-Anosov. This will allow us
to immediately recover generalizations of all the results we obtained for linear Anosov representations.

For the rest of the section, we will assume that $\mathsf{G}$ is a semisimple Lie group of non-compact type with finite center, denote its Lie algebra by $\mathfrak{g}$, 
and let ${\rm ad}: \mathfrak{g} \rightarrow \mathfrak{sl}(\mathfrak{g})$ and ${\rm Ad}:\mathsf G\to\SL(\mathfrak{g})$ be the adjoint representations.

Fix a parabolic subgroup $P^+ \subset \mathsf G$ and an opposite parabolic subgroup $P^- \subset \mathsf G$ and let $\Fc^\pm: = \mathsf G/P^{\pm}$ be the associated flag varieties. 
If $\rho : \Gamma \rightarrow \mathsf{G}$ is a  representation, we define the bundles
\begin{align*}
\wh B_\rho^\pm = \Gamma \backslash (\Usf(\Gamma) \times \Fc^{\pm})\quad\mathrm{and}\quad
\wh V_\rho^\pm = \Gamma \backslash(\Usf(\Gamma) \times T\Fc^{\pm}),
\end{align*}
where $ T\Fc^{\pm}$ is the tangent bundle of $ \Fc^{\pm}  $.
Observe that $V_\rho^\pm$ is a vector bundle over $B_\rho^\pm$ of rank $\dim(\Fc^\pm)$. 

The geodesic flow on $\Usf(\Gamma)$ extends to flows on $\Usf(\Gamma) \times \Fc^{\pm}$ and $\Usf(\Gamma) \times T\Fc^{\pm}$ whose action is trivial on the second factor.
These in turn descends to flows on $\wh B_\rho^\pm$ and $\wh V_\rho^\pm$ which covers the geodesic flow on $\wh{\Usf}(\Gamma)$. We use $\phi_t$ to denote all of these flows. 

We say that a  map 
$$\xi=(\xi^+,\xi^-):\Lambda(\Gamma)\to\Fc^+\times\Fc^-$$
is 
\begin{itemize}
\item \emph{transverse} if whenever $x\ne y\in\Lambda(\Gamma)$, the pair $(\xi^+(x),\xi^-(y))$ lies in the unique open $\mathsf{G}$-orbit in $\Fc^+\times\Fc^-$. 
\item \emph{strongly dynamics preserving} if whenever $\{\gamma_n\}$ is a sequence in $\Gamma$, $\gamma_n \rightarrow x$ and $\gamma_n^{-1} \rightarrow y$, then 
$$
\rho(\gamma_n) F \rightarrow \xi^+(x)
$$
for any $F \in \Fc^+$ which is transverse to $\xi^-(y)$. 
\end{itemize}

Given a transverse, $\rho$-equivariant, continuous map $\xi$ we define sections 
$$\sigma_\xi^\pm:\Usf(\Gamma)\to B_\rho^\pm=\Usf(\Gamma)\times \Fc^\pm$$
given by $\sigma_\xi^\pm(v)=\big(v,\xi^\pm(v^\pm)\big)$. Since $\xi$ is $\rho$-equivariant, $\sigma_\xi^\pm$ descend to  sections
$\wh\sigma_\xi^\pm:\wh\Usf(\Gamma)\to \wh B_\rho^\pm$.

\begin{definition}\label{defn:anosov_general} A representation $\rho : \Gamma \rightarrow \mathsf G$ is \emph{$P^\pm$-Anosov} if the following hold:
\begin{enumerate}
\item There exists a $\rho$-equivariant, continuous, transverse map $\xi=(\xi^+,\xi^-) : \Lambda(\Gamma) \rightarrow \Fc^+\times\Fc^-$.
\item For some  family of norms $\norm{\cdot}$ on the fibers of $\wh V_\rho^\pm\to\wh B_\rho^\pm$, the pullback of the flow $\phi_t$, also denoted $\phi_t$,
is uniformly expanding/contracting on $(\wh\sigma_\xi^\pm)^*\wh V_\rho^\pm$.
\end{enumerate}
We refer to any such map $\xi$ as the  \emph{$P^\pm$-Anosov limit map} for $\rho$.
\end{definition}

\begin{remark} With a little more work, one can show that the pullback of the flow $\phi_t$
is uniformly expanding on $(\wh\sigma_\xi^+)^*\wh V_\rho^+$ (with respect to some norm) if an only if it is uniformly contracting on
$(\wh\sigma_\xi^-)^*\wh V_\rho^-$ (with respect to some norm).
\end{remark}

\begin{remark} To be precise, the flow $\phi_t$ is uniformly expanding on $(\wh\sigma_\xi^+)^*\wh V_\rho^+$ if there exist $c,C > 0$ such that 
$$
\norm{\phi_t(Z)}_{\phi_t(\wh{\sigma}^+_\xi(v))} \geq C e^{ct} \norm{Z}_{\wh{\sigma}^+_\xi(v)}
$$
for all $t \geq 0$, $v \in \wh{\Usf}(\Gamma)$ and $Z \in \wh{V}_{\rho}^+|_{\wh{\sigma}^+_\xi(v)}$. Likewise,  the flow $\phi_t$ is uniformly contracting on $(\wh\sigma_\xi^-)^*\wh V_\rho^-$ if there exist $c,C > 0$ such that 
$$
\norm{\phi_t(Z)}_{\phi_t(\wh{\sigma}^-_\xi(v))} \leq C e^{-ct} \norm{Z}_{\wh{\sigma}^-_\xi(v)}
$$
for all $t \geq 0$, $v \in \wh{\Usf}(\Gamma)$ and $Z \in \wh{V}_{\rho}^-|_{\wh{\sigma}^-_\xi(v)}$.
\end{remark}

\begin{remark}
In the case when $\mathsf G=\SL(d,\Kb)$ and $(P^+,P^-)=(P_k,P_k^{opp})$ where 
$$
P_k = {\rm Stab}_{\mathsf G}(\ip{e_1,\dots, e_k}) \quad \text{and} \quad P_k^{opp}= {\rm Stab}_{\mathsf G}(\ip{e_{k+1},\dots, e_d}),
$$
we can identify $\Fc^+=\Gr_k(\Kb^d)$ and $\Fc^-=\Gr_{d-k}(\Kb^d)$. Then for any transverse pair $(F,G)\in \Gr_k(\Kb^d)\times\Gr_{d-k}(\Kb^d)$, there is a natural identification $T_F\Gr_k(\Kb^d)\simeq \Hom(F,G)$ and $T_G\Gr_{d-k}(\Kb^d)\simeq \Hom(G,F)$. Thus, the pullback bundles $(\wh\sigma_\xi^+)^*(\wh V_\rho^+)$ and $(\wh\sigma_\xi^-)^*(\wh V_\rho^-)$ are canonically identified with
$\Hom(\wh\Theta^k,\wh\Xi^{d-k})$ and $\Hom(\wh\Xi^{d-k},\wh\Theta^k)$ respectively, where $\wh\Theta^k$ and $\wh\Xi^{d-k}$ are the sub-bundles of $\wh E_\rho$ that lift to sub-bundles $\Theta^k$ and 
$\Xi^{d-k}$ of $E_\rho$ with the defining property $\Theta^k|_v=\xi^+(v^+)$ and $\Xi^{d-k}|_v=\xi^-(v^-)$. Thus, in this case, 
Definition \ref{defn:anosov_general} agrees with Definition \ref{defn:anosov_specific}, see Proposition \ref{prop: Hom bundle Anosov}.
\end{remark}

If $\psi:\mathsf{G}\to\mathsf{SL}(V)$ is a finite-dimensional irreducible representation, we say that $\psi$ is {\em adapted to} $(P^+,P^-)$
if $V=L_0\oplus W_0$ where $L_0$ is a line and $W_0$ is a hyperplane and
\begin{align*}
P^+ = \{ g \in \mathsf{G} : \psi(g)(L_0) = L_0 \} \quad \text{ and } \quad P^- = \{ g \in\mathsf  G : \psi(g)(W_0) = W_0\}.
\end{align*}
In this case, one may define embeddings $\psi_+ : \Fc^+ \rightarrow \Pb(V)$ and $\psi_- : \Fc^- \rightarrow \Gr_{\mathrm{dim}(V)-1}(V)$ by letting
\begin{align*}
\psi_+(gP^+) = \psi(g)(L_0) \text{ and } \psi_-(gP^-) = \psi(g)( W_0).
\end{align*}

The following result often allows one to reduce the general study of Anosov representations to the study of $P_1$-Anosov representations into $\mathsf{SL}(d,\mathbb R)$. For Anosov representations of word hyperbolic groups the analogous result is due to Guichard-Wienhard~\cite[Proposition 4.3]{guichard-wienhard}.

\begin{theorem}\label{thm:reduction_to_linear} 
Let $\mathsf{G}$ be a semisimple Lie group with finite center and let $P^\pm$ be a pair of opposite parabolic subgroups.
Suppose that $\psi : \mathsf{G} \rightarrow \SL(V)$ is a finite dimensional irreducible representation which is adapted to $(P^+,P^-)$.
Then a representation $\rho : \Gamma \rightarrow\mathsf{G}$ of a geometrically finite Fuchsian group $\Gamma$ is $P^\pm$-Anosov if and only if
 $\psi \circ \rho$ is $P_1$-Anosov. 

Moreover, if 
$\xi_\rho=(\xi_\rho^+,\xi_\rho^-)$ is a $P^\pm$-Anosov limit map for $\rho$,
then $\xi_{\psi\circ\rho}=(\psi_+\circ\xi_\rho^-,\psi_-\circ\xi_\rho^-)$ is the $P_1$-Anosov limit map  of $\psi \circ \rho$.
In particular, the $P^\pm$-Anosov limit map of $\rho$ is unique.
\end{theorem}

Adapted representations are not hard to construct.
If $\nL^+$ is the nilpotent radical of the Lie algebra of $P^+$, and $n = \dim \nL^+$, then $\psi(g) = \wedge^n \Ad(g)$ and 
$V = \Span \{ \psi(\mathsf G) (\wedge^n \nL^+)\} \subset \wedge^n \gL$ is adapted to $(P^+,P^-)$, see \cite[Remark 4.12]{guichard-wienhard}.
We obtain the following immediate corollary.

\begin{corollary}
\label{can do}
Suppose that  $\mathsf{G}$ is a semisimple Lie group with finite center and $P^\pm$ is a pair of opposite parabolic subgroups.
There exists a finite-dimensional irreducible representation $\psi : \mathsf{G} \rightarrow \SL(V)$  so that a representation
$\rho : \Gamma \rightarrow\mathsf{G}$  of a geometrically finite Fuchsian group $\Gamma$ is $P^\pm$-Anosov if and only if
 $\psi \circ \rho$ is $P_1$-Anosov. 
 \end{corollary}
 
Corollary \ref{can do} allows us to generalize our main results about linear Anosov representations into the general setting. As in the $\SL(d,\Kb)$ case, if $\rho:\Gamma\to \mathsf{G}$ is a  representation, let
$$\mathrm{Hom}_{\rm tp}(\rho)\subset\mathrm{Hom}(\Gamma,\mathsf{G})$$ be the space of representations $\sigma:\Gamma\to \mathsf{G}$ so that  if $\alpha\in\Gamma$ is parabolic, then $\sigma(\alpha)$ is conjugate to $\rho(\alpha)$. Theorem \ref{thm: stability intro} becomes:

\begin{corollary}
Suppose that  $\mathsf{G}$ is a semisimple Lie group with finite center, $P^\pm$ is a pair of opposite parabolic subgroups of $\mathsf{G}$
and $\Gamma$ is a geometrically finite Fuchsian group.
If  $\rho:\Gamma\to\mathsf{G}$
is $P^\pm$-Anosov, then there exists an open neighborhood $\Oc$ of $\rho$ in $\mathrm{Hom}_{\rm tp}(\rho)$, so that
\begin{enumerate}
\item
If $\rho\in \Oc$, then $\rho$ is $P^\pm$-Anosov.
\item
There exists $\alpha>0$ so that if $\rho\in\Oc$, then its $P^\pm$-Anosov limit map $\xi_\rho$ is $\alpha$-H\"older.
\item
If $\{\rho_u\}_{u\in M}$ is an analytic family of representations in $\Oc$ and  $z\in\Lambda(\Gamma)$, then the map from $M$ to $\Fc^+\times\Fc^-$
given by $u\to \xi_{\rho_u}(z)$ is analytic.
\end{enumerate}
\end{corollary}

Theorem \ref{thm:equivalent_to_Anosov} yields:

\begin{corollary}
Suppose that  $\mathsf{G}$ is a semisimple Lie group with finite center, $P^\pm$ is a pair of opposite parabolic subgroups of $\mathsf{G}$
and $\Gamma$ is a geometrically finite Fuchsian group.
A representation $\rho : \Gamma \rightarrow \mathsf{G}$ is $P^\pm$-Anosov if and only if 
there exists a $\rho$-equivariant, transverse, continuous, strongly dynamics preserving map 
$\xi=(\xi^k,\xi^{d-k}):\Lambda(\Gamma)\to\Fc^+\times\Fc^-$.
Furthermore, $\xi$ is the $P^\pm$-Anosov limit map.
\end{corollary}

One can also obtain analogues of parts (1) and (2)  of  Theorem \ref{thm: intro thm 1} where the roles of singular values and eigenvalues are
played by roots acting on the Cartan and Jordan projections (see \cite[Section 3]{GGKW} for a complete discussion). Part (4) of Theorem \ref{thm: intro thm 1}
remains true if we replace $X_d(\mathbb K)$ with the symmetric space of $\mathsf{G}$.

\medskip

Theorem~\ref{thm:reduction_to_linear} will be a consequence of Theorem \ref{thm:equivalent_to_Anosov} and the following dynamical property of Anosov representations. 

\begin{lemma}\label{lem:strong_dynamics} Suppose $\rho : \Gamma \rightarrow \mathsf{G}$ is a $P^\pm$-Anosov representation of a geometrically finite Fuchsian group $\Gamma$
with $P^\pm$-Anosov limit map $\xi=(\xi^+,\xi^-)$. If $\{\gamma_n\}$ is a sequence in $\Gamma$ such that $\gamma_n \rightarrow x\in\Lambda(\Gamma)$ and 
$\gamma_n^{-1} \rightarrow y\in\Lambda(\Gamma)$, then 
\begin{align*}
\lim_{n \rightarrow \infty} \rho(\gamma_n)(F) \rightarrow \xi^+(x)
\end{align*}
for all $F \in \Fc^+$ transverse to $\xi^-(y)$. 
\end{lemma}

Delaying the proof of the lemma we prove Theorem~\ref{thm:reduction_to_linear}.

\medskip\noindent
{\em Proof of Theorem~\ref{thm:reduction_to_linear}.}
We make repeated use of the following observation which follows from \cite[Prop 3.5]{GGKW}. We sketch an alternate proof.

\begin{observation}\label{obs:transverse} $(F,H) \in \Fc^+ \times \Fc^-$ are transverse if and only if $\psi_+(F)$ and
$\psi_-(H)$ are transverse. \end{observation}

\begin{proof}[Sketch of proof] Fix maximal compact subgroups $K_1 \subset \mathsf{G}=:\mathsf G_1$ and $K_2 \subset \SL(V)=:\mathsf G_2$ so that $\psi(K_1) \subset K_2$. For $j=1,2$, we can fix $\mathsf G_j$-invariant Riemannian metrics on $X_j = \mathsf{G}_j/K_j$ so that the map $T_\psi:X_1\to X_2$ given by $T_\psi(gK_1) = \psi(g)(K_2)$ is a totally geodesic isometric embedding (see~\cite[Chapter 2]{mostow}). Then $T_\psi$ extends to an embedding $T_\psi:X_1(\infty) \hookrightarrow X_2(\infty)$ 
of the $\CAT(0)$-boundaries. 

Fix $(F,H)=(g_+P^+, g_-P^-)$. Let $W_+, W_- \subset X_1(\infty)$ be the interior of the Weyl faces associated to $g_+P^+g_+^{-1}$ and $g_-P^-g_-^{-1}$ respectively, i.e.
\begin{align*}
W_\pm = \left\{ x \in X_1(\infty) : {\rm Stab}_{\mathsf{G}}(x) = g_\pm P^\pm g_\pm^{-1}\right\}.
\end{align*}
Also, let $\hat{W}_+, \hat{W}_- \subset X_2(\infty)$ denote the interior of the Weyl faces associated to the parabolic subgroups $\psi(g_+)P_1 \psi(g_+)^{-1}$ and $\psi(g_-)P_{d-1}\psi(g_-)^{-1}$ respectively. Notice that 
\begin{align*}
W_\pm = T_\psi^{-1}(\hat{W}_\pm)
\end{align*}
since $P^+ = \psi^{-1}(P_1)$ and $P^-=\psi^{-1}(P_{d-1})$. 

Next fix a maximal flat $\Fc \subset X_1$ with $W_+, W_- \subset \Fc(\infty)$ and fix some $p \in \Fc$.  Then let $s_p : X_1 \rightarrow X_1$ denote the involutive isometry based at $p$. Then $F$ and $H$ are transverse if and only if $W_+=s_p(W_-)$. 
Since $T_\psi : X_1 \rightarrow X_2$ is a totally geodesic isometric embedding, there exists a maximal flat $\hat{\Fc} \subset X_2$ with $T_\psi(\Fc) \subset \hat{\Fc}$. 
Then $\hat{W}_+,\hat{W}_- \subset \hat{\Fc}(\infty)$ and if $s_{\hat{p}} : X_2 \rightarrow X_2$ is the involutive isometry based at $\hat{p}:=T_\psi(p)$, then 
$\psi(F)$ and  $\psi(H)$ are transverse if and only if $\hat{W}_+ = s_{\hat{p}}(\hat{W}_-)$. 

Since 
\begin{align*}
T_\psi \circ s_p = s_{\hat{p}} \circ T_\psi 
\end{align*}
and any two distinct interiors of Weyl faces have trivial intersection:
\begin{align*}
W_+ = s_p(W_-) &  \Longrightarrow T_\psi(W_+) \subset \hat{W}_+\cap s_{\hat{p}}(\hat{W}_-)\Longrightarrow \hat{W}_+= s_{\hat{p}}(\hat{W}_-) \Longrightarrow T_\psi(s_p(W_-)) \subset \hat{W}_+ \\
& \Longrightarrow s_p(W_-) \subset W_+ = T_\psi^{-1}(\hat{W}_+) \Longrightarrow W_+ = s_p(W_-). 
\end{align*}
So $F$ and $H$ are transverse if and only if $\psi_+(F)$ and  $\psi_-(H)$ are transverse. 
\end{proof} 

We first prove the reverse direction of Theorem~\ref{thm:reduction_to_linear}.

\begin{lemma} 
If $\psi \circ \rho$ is $P_1$-Anosov, then $\rho$ is $P^\pm$-Anosov. Moreover, if $\xi_\rho$ is the $P^\pm$-Anosov limit map of $\rho$ and
$\xi_{\psi\circ\rho}$ is the $ P_1$-Anosov limit map of $\psi\circ\rho$, then 
$\psi_{\pm}\circ\xi_\rho^\pm=\xi_{\psi\circ\rho}^\pm$.
\end{lemma}

\begin{proof} 
If $\gamma$ is a hyperbolic element of $\Gamma$, then, since $\psi\circ\rho$ is $P_1$-Anosov, $\psi(\rho(\gamma))$ is $P_1$-biproximal and 
$\xi_{\psi\circ\rho}^+(\gamma^+)$ is the attracting eigenline of $\psi(\rho(\gamma))$. Since $\psi$ is irreducible, there exists $F\in\Fc^+$ such that $\psi_+(F)$ is transverse to the repelling hyperplane of $\psi(\rho(\gamma))$, so
$$\xi_{\psi\circ\rho}^+(\gamma^+) = \lim_{n \rightarrow \infty} (\psi \circ \rho)(\gamma)^n (\psi_+(F))= \lim_{n \rightarrow \infty} \psi_+( \rho(\gamma)^n (F)) \in \psi_+(\Fc^+).$$
Hence,
$\xi_{\psi\circ\rho}^+(\Lambda(\Gamma))\subset \psi_+(\Fc^+)$ because attracting fixed points of hyperbolic elements are dense in $\Lambda(\Gamma)$. 

Since $\psi_+$ is a $\psi$-equivariant  embedding, $\xi^+=\psi_+^{-1}\circ \xi_{\psi\circ\rho}^+$ is well-defined, continuous and $\rho$-equivariant. Similarly,
$\xi^-=\psi_-^{-1}\circ \xi_{\psi\circ\rho}^-$ is well-defined, continuous and $\rho$-equivariant. Observation~\ref{obs:transverse} implies that $\xi=(\xi^+, \xi^-)$ is transverse. 

We consider the vector bundles $\wh V_\rho^\pm$ and $\wh V_{\psi\circ\rho}^\pm$ over $\wh B_\rho^\pm$ and $\wh B_{\psi\circ\rho}^\pm$ respectively.
Notice that the map $\psi_\pm$ induces a bundle embedding $\iota_\rho^\pm : \wh V_\rho^\pm \hookrightarrow \wh V_{\psi\circ\rho}^\pm$ which intertwines the flows on the two bundles. Since $\psi\circ\rho$ is $P_1$-Anosov, there is a continuous family of norms on the fibers of the bundle $\wh V_{\psi\circ\rho}^\pm \rightarrow \wh B_{\psi\circ\rho}^\pm$ such that $\psi_t$ is uniformly expanding/contracting on the pullback bundle $(\wh\sigma_{\xi_{\psi\circ\rho}}^\pm)^*\wh V_{\psi\circ\rho}^\pm$. Equip the bundle $\wh V_\rho^\pm \rightarrow \wh B_\rho^\pm$ with the pullback of this norm via $\iota_\rho^\pm$. Since $\iota_\rho^\pm$ intertwines the flows, we see that the  flow is uniformly expanding/contracting on $\sigma_\xi^*(V_\rho^\pm)=(\iota_\rho^\pm)^*\left((\wh\sigma_{\xi_{\psi\circ\rho}}^\pm)^*\wh V_{\psi\circ\rho}^\pm\right)$.
Therefore,
$\rho$ is $P^\pm$-Anosov with $P^\pm$-Anosov limit map $\xi=(\xi^+,\xi^-)$. 
\end{proof}

We now prove the forward direction of Theorem~\ref{thm:reduction_to_linear}.

\begin{lemma}
If $\rho$ is $P^\pm$-Anosov with $P^\pm$-Anosov limit maps $\xi_\rho$, then $\psi \circ \rho$ is $P_1$-Anosov with $P_1$-Anosov limit map 
$(\psi_+\circ\xi_\rho^+,\psi_- \circ \xi_\rho^-)$.
\end{lemma}

\begin{proof} 
Let $\eta=(\eta^+,\eta^-)=(\psi_+\circ\xi_\rho^+,\psi_- \circ \xi_\rho^-)$. Then $\eta$ is  continuous, $\psi \circ \rho$-equivariant, and transverse (by Observation~\ref{obs:transverse}). 
So, by Theorem \ref{thm:equivalent_to_Anosov}, it suffices to show that $\eta$ is strongly dynamics-preserving. 

Consider a sequence $\{\gamma_n\}$ in $\Gamma$ with $\gamma_n \rightarrow x\in\Lambda(\Gamma)$ and $\gamma_n^{-1} \rightarrow y\in\Lambda(\Gamma)$.
Then, by Lemma~\ref{lem:strong_dynamics},
\begin{align*}
\lim_{n \rightarrow \infty} (\psi \circ \rho)(\gamma_n) (\psi_+(F)) = \eta^+(x)
\end{align*}
for all $F \in \Fc^+$ transverse to $\eta^-(y)$. Similarly, 
\begin{align*}
\lim_{n \rightarrow \infty} (\psi \circ \rho)(\gamma_n) (\psi_-(F)) = \eta^-(x)
\end{align*}
for all $F \in \Fc^-$ transverse to $\eta^+(x)$.

Since $\psi$ is irreducible, $\psi_+(\Fc^+)$ spans $V$ and one can repeat the proof of Corollary~\ref{zariski dense case} to show that $\eta$ is strongly dynamics preserving. 
\end{proof}

It only remains to prove Lemma \ref{lem:strong_dynamics}.

\medskip\noindent
{\em Proof of Lemma~\ref{lem:strong_dynamics}.} Let $\pL^\pm$ be the Lie algebra of $P^\pm$. Then there exists a Cartan decomposition $\gL = \kL \oplus \pL$, a Cartan subspace $\aL \subset \pL$, and an element $H_0 \in \aL$ so that 
$$\pL^\pm = \gL_0 \oplus \bigoplus_{\alpha(\pm H_0) \geq 0} \gL_\alpha$$
where 
$$\gL = \gL_0 \oplus \bigoplus_{\alpha \in \Sigma} \gL_\alpha$$
is the root space decomposition associated to $\aL$. Let $\nL^\pm = \bigoplus_{\alpha(\pm H_0) > 0} \gL_\alpha$ and define
$$T : \nL^{-} \rightarrow \Fc^+ \quad\mathrm{where}\quad
T(X) = e^{X}P^{+}$$
The map $T$ has the following properties. 

\begin{observation}\label{obs:properties_of_the_map} \
\begin{enumerate}
\item $T(\nL^-) = \{ F\in \Fc^+ : F  \text{ is transverse to } P^-\}$.
\item If $H \in \aL$, then $e^H \circ T = T\circ \Ad(e^H)$. 
\item $d(T)_0 : \nL^{-} \rightarrow T_{P^+} \Fc^+$ is a linear isomorphism. 
\end{enumerate}
\end{observation}

\begin{proof}By definition, $F$ is transverse to $P^-$ if and only if $(F,P^-)  \in \mathsf{G} \cdot (P^+, P^-)$ if and only if $F=gP^+$ for some $g \in P^-$. By the Langlands decomposition, $P^- = N^-(P^+ \cap P^-)$ where $N^- \subset \mathsf{G}$ is the connected Lie subgroup with Lie algebra $\nL^-$.  Since $N^-$ is nilpotent, $N^- = e^{\nL^-}$. So $F$ is transverse to $P^-$ if and only if $F=e^XP^+$ for some $X \in \nL^-$. This proves part (1). 

Part (2) is an immediate consequence of the definition. Part (3) follows from the fact that $\gL = \nL^- \oplus \pL^+$ (as vector spaces) and $\pL^+$ is the Lie algebra of $P^+$. 
\end{proof}

As a consequence of (1) and (2) in Observation~\ref{obs:properties_of_the_map}, we have the following. 

\begin{lemma}\label{cor:stupid_dynamics}
If $\{H_n\}$ is a sequence in $\aL$ with $\lim_{n \rightarrow \infty} \alpha(H_n) = -\infty$ for all $\alpha \in \Sigma$ with $\alpha(H_0) < 0$, then 
\begin{align*}
\lim_{n \rightarrow \infty} e^{H_n} (F) = P^+
\end{align*}
for all $F\in \Fc^+$ transverse to $P^-$. 
\end{lemma}

\begin{proof}
By Observation~\ref{obs:properties_of_the_map} (1), $F=T(X)$ for some $X\in\mathfrak n^-$. Write $X=\sum_{\alpha(H_0)<0}X_\alpha$, 
where $X_\alpha\in\mathfrak g_\alpha$. Then by Observation~\ref{obs:properties_of_the_map} (2), 
\[e^{H_n}(F)=T\left(\Ad(e^{H_n})\left(\sum_{\alpha(H_0)<0}X_\alpha\right)\right)=T\left(\sum_{\alpha(H_0)<0}e^{\mathrm{ad}(H_n)}X_\alpha\right)=T\left(\sum_{\alpha(H_0)<0}e^{\alpha(H_n)}X_\alpha\right).\]
Since $\lim_{n \rightarrow \infty} \alpha(H_n) = -\infty$ for all $\alpha \in \Sigma$ with $\alpha(H_0) < 0$, it follows that 
\[\lim_{n\to\infty}e^{H_n}(F)=T(0)=P^+.\qedhere\]
\end{proof}

Let $K \subset \mathsf{G}$ be the maximal compact subgroup with Lie algebra $\kL$ and fix a $K$-invariant Riemannian metric on $\Fc^+$, and let $|\cdot|$ denote the induced family of norms on the fibers of $T\Fc^+$. 

Let $\xi=(\xi^+,\xi^-)$ be the $P^\pm$-Anosov limit map for $\rho$. Then let $\sigma_\xi^+(v) =( v, \xi^+(v^+))$. Since $\rho$ is Anosov, there is a $\rho$-equivariant family of norms on the fibers of $\Usf(\Gamma) \times T\Fc^{\pm} \rightarrow \Usf(\Gamma) \times \Fc^{\pm}$ and constants  $C, c> 0$ such that 
\begin{align*}
\norm{Z}_{\phi_{-t}(\sigma_\xi^+(v))} \leq Ce^{-c t} \norm{Z}_{\sigma_\xi^+(v)}
\end{align*}
for all $t > 0$, $v \in \Usf(\Gamma)$ and $Z \in T_{\xi^+(v^+)} \Fc^+$.

Consider an escaping sequence $\{\gamma_n\}$ with $\gamma_n \rightarrow x$ and $\gamma_n^{-1} \rightarrow y$.

\medskip
\noindent \textbf{Case 1:} If  $x \neq y$, then $\gamma_n$ is hyperbolic when $n$ is sufficiently large.
Furthermore, we can find a bounded sequence $\{v_n\}$ in $\Usf(\Gamma)$ such that $v_n^\pm = \gamma_n^\pm$, and a bounded sequence $\{g_n\}$ in $\mathsf{G}$ such that 
\begin{align*}
g_n (\xi^+(\gamma_n^+), \xi^-(\gamma_n^-)) = (P^+, P^-).
\end{align*}
Then 
$$g_n \rho(\gamma_n) g_n^{-1} P^\pm = P^\pm\quad\mathrm{so}\quad
g_n \rho(\gamma_n) g_n^{-1} \in L:=P^+ \cap P^-.$$

Notice that  
\begin{align*}
\gL_0 \oplus \bigoplus_{\alpha(H_0)=0} \gL_\alpha
\end{align*}
is a root space decomposition of the Lie algebra of $L$. Then, using the Cartan decomposition of the reductive group $L$, there exist $k_{n,1}, k_{n,2} \in K \cap L$ and $H_n \in \aL$ so that 
\begin{align*}
g_n \rho(\gamma_n) g_n^{-1} = k_{n,1} e^{H_n} k_{n,2}.
\end{align*}

Since $\{v_n\}$ is a bounded sequence there exists $C_1 > 1$ such that 
\begin{align*}
\frac{1}{C_1} |Z|_{\xi^+(v_n^+)} \leq \norm{Z}_{\sigma^+_\xi(v_n)} \leq C_1 |Z|_{\xi^+(v_n^+)}
\end{align*}
for all $n \geq 1$ and $Z \in T_{\xi^+(v_n^+)}\Fc^+$. Likewise, there exists $C_2 > 1$ such that 
\begin{align*}
\frac{1}{C_2} |Z|_F \leq |g_n(Z)|_{g_n(F)}  \leq C_2 |Z|_F
\end{align*}
for all $n \geq 1$, $F \in \Fc^+$ and $Z \in T_F\Fc^+$.

Notice that $\gamma_n^{-1} (v_n) = \phi_{-t_n}(v_n)$ 
for some sequence $\{t_n\}$ with $t_n \rightarrow \infty$. Since both $k_{n,1}$ and $k_{n,2}$ fix $P_+$ and $|\cdot|$ is a $K$-invariant family of norms, it follows that for any $Z \in T_{P^+} \Fc^+$, we have 
\begin{align}\label{eqn: long}
|e^{H_n} (Z)|_{P^+} & = |k_{n,1}^{-1} g_n \rho(\gamma_n) g_n^{-1} k_{n,2}^{-1}( Z)|_{P^+} \leq C_2 | \rho(\gamma_n) g_n^{-1} k_{n,2}^{-1}( Z)|_{\xi^+(v_n^+)} \nonumber\\
 &\leq C_1C_2  \norm{\rho(\gamma_n) g_n^{-1} k_{n,2}^{-1} (Z)}_{\sigma_\xi^+(v_n)}= C_1C_2  \norm{g_n^{-1} k_{n,2}^{-1}(Z)}_{\phi_{-t_n}(\sigma_\xi^+(v_n))}\nonumber\\
 &\leq C_1C_2Ce^{-c t_n} \norm{g_n^{-1} k_{n,2}^{-1} (Z)}_{\sigma_\xi^+(v_n)}\leq C_1^2C_2Ce^{-c t_n} |g_n^{-1} k_{n,2}^{-1} (Z)|_{\xi^+(v_n^+)}\\
& \leq C_1^2C_2^2Ce^{-c t_n} |Z|_{P^+}.\nonumber
\end{align}

By Observation~\ref{obs:properties_of_the_map} (3), we know that for any $\alpha\in\Sigma$ such that $\alpha(H_0)<0$, and any $X\in\mathfrak g_\alpha$, there is some $Z\in T_{P^+}\Fc^+$ such that $d(T)_0(X)=Z$. Then
\[e^{H_n}(Z)=d(e^{H_n}\circ T)_0(X)=\frac{d}{dt}\bigg|_{t=0}e^{H_n}\circ T(tX).\]
Then by Observation~\ref{obs:properties_of_the_map} (2) (see proof of Lemma \ref{cor:stupid_dynamics}),
\[\frac{d}{dt}\bigg|_{t=0}e^{H_n}\circ T(tX)=\frac{d}{dt}\bigg|_{t=0}T(te^{\alpha(H_n)}X)=e^{\alpha(H_n)}\frac{d}{dt}\bigg|_{t=0}T(tX)=e^{\alpha(H_n)}Z.\]
Thus, $e^{H_n}(Z)=e^{\alpha(H_n)}Z$, so the inequality \eqref{eqn: long} implies that 
\begin{align*}
\lim_{n \rightarrow \infty} \alpha(H_n) = -\infty
\end{align*}
whenever $\alpha(H_0) < 0$. Hence, by Lemma~\ref{cor:stupid_dynamics},
\begin{align*}
\lim_{n \rightarrow \infty} e^{H_n}( F) = P^+
\end{align*}
for all $F \in \Fc^+$ transverse to $P^-$. Since $g_n (\xi^+(x), \xi^-(y)) \rightarrow (P^+,P^-)$, $k_{n,j}P^\pm = P^\pm$ and 
$\rho(\gamma_n) = g_n^{-1} k_{n,1} e^{H_n} k_{n,2} g_n$ we then have 
\begin{align*}
\lim_{n \rightarrow \infty} \rho(\gamma_n) (F)= \xi^+(x)
\end{align*}
for all $F \in G/P^+$ transverse to $\xi^-(y)$. 

\medskip
\noindent \textbf{Case 2:}  If $x = y$, pick $\beta \in \Gamma$ so that $z:=\beta^{-1}( x) \neq x$. Then $\gamma_n\beta \rightarrow x$ and  $(\gamma_n\beta)^{-1} \rightarrow z \neq x$. 
By the first case, $\rho(\gamma_n\beta)(F) \rightarrow \xi(x)$ for all $F \in \Fc^+$ transverse to $\xi^-(z) = \rho(\beta^{-1})\xi^-(x)$. 
Equivalently, $\rho(\gamma_n)(F) \rightarrow \xi^+(x)$ for all $F \in \Fc^+$ transverse to $\xi^-(x)$. 
\qed

\end{document}